\documentstyle[11pt]{article}
\title{A singularity at the criticality  for  the free energy in  percolation
} 
\author{\small{\bf Dedicated to Harry Kesten  for his contributions to percolation theory}\\
${}$\\
Yu Zhang}
\begin{document}
\baselineskip .20in
\maketitle
\begin{abstract}
Consider percolation on the triangular  lattices. 
Let $\kappa(p)$ be the free energy at the zero field.  
 We show that
$$|\kappa'''(p)| \leq  |p-p_c|^{-1/3+o(1)} \mbox{ if } p \neq  p_c.$$
Furthermore, we show that there is a sequence  $\epsilon_n \downarrow 0$ such that
$$|\kappa'''(p_c\pm \epsilon_n)|\geq \epsilon_n^{-1/3+o(1)}.$$
The last inequality  implies that $\kappa(p)$ is not third differentiable.  This  answers  affirmatively a  conjecture, asked by  Sykes and Essam in 1964, whether $\kappa(p)$ has a singularity at the criticality.\\

\end{abstract}

\section{ Introduction and statement of results.}
Consider site percolation on the triangular lattice.
We may realize the triangular lattice with vertex set ${\bf Z}^2$. For a given $(x, y)\in {\bf Z}^2$,
its nearest six neighbors are defined as $(x\pm 1, y),(x, y\pm1), (x+1, y+1)$, and $(x-1, y-1)$.
Edges between neighboring or adjacent vertices therefore correspond to vertical or horizontal displacements of one unit, 
or diagonal displacements between the two nearest vertices along  lines making  angles $\pi/4$ and $5\pi/4$
with the positive $X$-axis.  Note that each site is the center of a hexagon in its dual graph.
Each site or hexagon is independently {\em occupied} with 
probability $p$ and {\em vacant} with probability $1-p$. 
The corresponding probability measure on the configurations of occupied
and vacant vertices is denoted by $P_p$. We also denote by $E_p$ the expectation with respect to $P_p$. 
A path from $u$ to $v$ is a sequence $(v_0,...,v_{i}, v_{i+1},...,v_n)$
with adjacent  vertices $v_i$  and $v_{i+1}$  or  hexagons 
sharing an edge ($0\leq i\leq n-1$) such that $v_0=u$ and $v_n=v$. 
A circuit is a path with distinct vertices $v_i$ ($1\leq i\leq n$) except $v_0=v_n$. 
A path is called occupied or vacant if all of its vertices or hexagons  are occupied or vacant. The occupied cluster of the vertex $x$, 
${\bf C}(x)$,
consists of all vertices or hexagons that are connected to $x$ by an occupied path.
${\bf C}(x)$ is an empty set if $x$ is vacant.
For any collection $A$ of vertices, $|A|$ 
denotes the cardinality of $A$. We choose ${\bf 0}$ as the origin.
The {\em percolation probability} is
\begin{eqnarray*}
\theta (p)= P_p(|{\bf C}({\bf 0})|=\infty),
\end{eqnarray*}
and the critical probability is 
$$p_c=\sup\{p:\theta (p)=0\}.$$
We may also consider site or bond percolation on a periodic two dimensional  lattice (see Kesten (1982) Chapters 1-3 for a precise description of terminology).
It has been proved  (see chapter 3 in Kesten (1982)) that for site percolation on the triangular lattice 
$$p_c=0.5.\eqno{(1.1)}$$

We denote the cluster distribution by
$$\theta_n(p)=P_p(|{\bf C}({\bf 0})|=n).$$
By analogy with the Ising model, we introduce the {\em magnetization function} as
$$M(p,h)=1-\sum_{n=0}^\infty \theta_n(p) e^{-nh}\mbox{ for }h\geq 0.$$
By setting $h=0$ in the magnetization function,
$$M(p,0)=\theta(p).$$
Using  term-by-term differentiation,
we also have
$$\lim_{h\rightarrow 0^+}{\partial M(p,h)\over \partial h}=E_p( |{\bf C}({\bf 0})|;|{\bf C}({\bf 0})|<\infty)=\chi^f(p).$$
$\chi^f(p)$ is called the {\em mean cluster size}. 
The {\em free energy} $F(p,h)$ is defined by
$$F(p,h)=h(1-\theta_0(p))+\sum_{n=1}^\infty {1\over n} \theta_n(p)e^{-hn}\mbox{ for } h >0.$$
If we differentiate with respect to $h$, we find
$$ {\partial F(p,h)\over \partial h}=M(p,h).$$
For $h>0$, the free energy is infinitely differentiable with respect to $p$. 
The zero-field free energy $F(p,0)$ is a more interesting subject of study. By our definition,
$$F(p,0)=E_p(|{\bf C}({\bf 0})|^{-1}; |{\bf C}({\bf 0})| >0).\eqno{(1.2)}$$

Grimmett  (1981) discovered that the zero-field  free energy  also coincides  with the {\em number of clusters per vertex}. 
Let us define the number of clusters per vertex as follows. Note that  any two vertices $x,y\in B(n)=[-n,n]^2$
are said to be connected in $B(n)$ 
if either $x=y$ or there exists an occupied  path $\gamma$  in $B(n)$ connecting $x$ 
and $y$.
Let $M_n$ be the 
number of occupied clusters in $B(n)$. By a standard ergodic theorem (see Theorem 4.2 in
Grimmett  (1999)), the limit
$$\lim_{n\rightarrow \infty}{1\over {|B(n)|}}M_n=\kappa(p) 
\mbox { a.s. and }L_1\eqno{(1.3)}$$ 
exists for all $0\leq p\leq 1$.
Let $K_n=E_p(M_n).$ Thus,
$$\lim_{n\rightarrow \infty}{1\over {|B(n)|}}K_n(p)=\kappa(p). $$
$\kappa (p)$ is called the number of clusters per vertex. Grimmett (1981) proved  
that 
$$\kappa(p)=F(p,0). \eqno{(1.4)}$$

Sykes and Essam were perhaps the first to introduce the number of clusters per vertex in 1964, and they tried to use it to compute $p_c$.  
They explored a beautiful geometric argument in their paper to show
$$\kappa(p)-\kappa(1-p)=p-3p^2+2p^3.\eqno{(1.5)} $$
Sykes and Essam argued that phase transition in percolation must be manifested by a singularity 
at the critical value $p_c$. If $p_c$ is indeed the only singularity of $\kappa(p)$, then (1.5) implies that
$p_c=0.5$.
For many years, the singularity criterion of  Sykes and Essam has offered a tantalizing approach
to the famous problem that $p_c=0.5$.  Kesten used another method in 1980 to show that
$p_c=0.5$ (see Kesten (1982)). However, until the present paper, there was no proof of a singularity at $p_c$
for the free energy function.

We would like to mention some progress for $\kappa(p)$ throughout the years.
It has been ruled out that $\kappa(p)$ has  another singularity on $p$ for  $p\neq p_c$.
In other words, $\kappa(p)$ is analytic for $p\neq p_c$
(see chapter 9 in Kesten (1982)). On the other hand, it has also been proved (see chapter 9 in
Kesten (1982)) that $\kappa(p)$ is twice differentiable at $p_c$. This tells us that $\kappa(p)$ is a very smooth function.
Indeed, the smoothness of $\kappa(p_c)$ might tell us why the singularity at $p_c$ is difficult to  prove. 
The main result obtained here is to understand the behavior of $\kappa$ at the critical point. 
If $p_c$  is indeed a singularity of $\kappa(p)$, then 
it is natural to ask about the behavior of the singularity. Physicists believe that the zero-field free energy is not third  differentiable. 
It is believed that the behavior of  percolation functions can be described in terms  of critical exponents 
as $p$ approaches  $p_c$. For $\kappa(p)$, it is conjectured that there exists an exponent $\alpha$ such that
$$\kappa'''(p)\approx |p-p_c|^{-1-\alpha}.\eqno{(1.6)}$$
It is not known  how strong one expects such an asymptotic $``\approx"$  relation to be, and it is for this reason that we shall
use the logarithmic relation. More precisely, 
$f(p)\approx g(p)$ or $f_n\approx g_n$ means  
$$\log g(p)/ \log f(p)\rightarrow 1 \mbox{ or } g(p)=f^{1+o(1)}(p)\mbox{ as }p\rightarrow p_c,\mbox{ or } \log f_n/\log g_n\rightarrow 1 \mbox{ or } f_n= g_n^{1+o(1)}\mbox{ as } n\rightarrow \infty.$$

The exponent $\alpha$ is called the {\em heat exponent}, and (1.6) is called the {\em power law} for the free energy. Numerical computations indicate
$\alpha=-2/3$. In addition to this power law, it
is also widely believed that the exponents satisfy the following so-called
{\em scaling laws}. To be more specific, we need to introduce  all the other critical exponents and power laws.
We denote the {\em correlation length} by
\begin{eqnarray*}
\xi^{-1}(p)=\lim_{n\rightarrow \infty}\{-{1\over n}\log P_p({\bf 0}\rightarrow
\partial B(n), |{\bf C}({\bf 0})|<\infty)\} \mbox{ if } p \neq  p_c,
\end{eqnarray*}
and the {\em probability on the tail} of $|{\bf C}({\bf 0})|$ at $p_c$ by
\begin{eqnarray*}
\pi(n)=P_{p_c}(n\leq |{\bf C}({\bf 0})|< \infty),
\end{eqnarray*}
where $\partial B(n)$ is the surface of the box $B(n)$
 and $A\rightarrow B$
means that there exists an occupied path from
some vertex of $A$ to some vertex of $B$ for any sets $A$ and $B$.

The power laws  are introduced as follows:
$$\theta(p)\approx (p-p_c)^{\beta} \mbox{ for } p> p_c,$$
$$\chi^f(p)
\approx (p_c-p)^{-\gamma} \mbox{ for } p\neq p_c ,$$
$$\xi(p)
\approx |p_c-p|^{-\nu} \mbox{ for } p\neq p_c , \eqno{(1.7)}$$
$$\pi(n)
\approx n^{-1/\rho}\mbox{ for }n\geq 1, \eqno{(1.8)}$$
$$\kappa''' (p)
\approx |p-p_c|^{-1-\alpha} \mbox{ for } p\neq p_c .\eqno{(1.9)}$$
Numerical computations indicate that
$$\beta={5\over 36}, \gamma={43\over 18},\nu={4\over 3}, \rho={91\over 5}.$$

In addition to the power laws, it
is also widely believed that the exponents satisfy the following so-called
scaling laws:
$$\alpha=2-2\nu, \eqno{(1.10)}$$
$$ \beta={2\nu\over \rho+1},\eqno{(1.11)}$$
$$\gamma=2\nu{\rho-1\over \rho+1}.\eqno{(1.12)}$$
In particular, (1.10) is called a {\em hyper scaling} relation. 
Moreover, let us introduce $k$-arm paths. Consider the annulus 
$$A(m,n)=\{B(n)\setminus B(m)\}\cup \{\partial B(m)\}\mbox{ for } m <n.$$
Let ${\cal Q}_k(b, m,n)$ be the event that there exist $i$ disjoint occupied paths and $j$ disjoint vacant paths with $i+j=k$ for all $i,j\geq 1$
from $b+\partial B(m)$ to $b+\partial B(n)$ inside $b+A(m,n)$.
We call them $k$-{\em arm paths}.
For simplicity, let ${\cal Q}_k(m,n)={\cal Q}_k({\bf 0}, m, n)$. If $b\in B(n)$, let ${\cal Q}_k(b, n)$ be the $k$-arm paths from $b$ to $\partial B(n)$. Let ${\cal Q}_4(n)={\cal Q}_4({\bf 0}, n)$ for simplicity.
It is believed by Aizenman, Dulpantier, and Alharony (1999) that
$$P_{0.5}({\cal Q}_k(m,n))= \left({m\over n}\right)^{(k^2-1)/12+o(1)}\eqno{(1.13)}$$
for fixed $m$ as $n\rightarrow \infty$ and for $k>1$.
In fact, it is more important to show (1.13) when $k=1$ and $k=4$. Fortunately, 
by using the Schramm-Loewner evolution (SLE) argument and 
 Smirnov's  scaling limit on the triangular lattice (2001),
(1.13) was  proved  (see Lawler, Schramm, and Werner  (2001) and  Smirnov and Werner's Theorem 4  (2001)) for $k=1$ and  $k=4$. 
More precisely, on the triangular lattice, 
$$P_{0.5}({\cal Q}_4(m,n))= \left({m\over n}\right)^{5/4+o(1)}\mbox{ and } P_{0.5}({\cal Q}_1(m,n))= \left({m\over n}\right)^{5/48+o(1)}.\eqno{(1.14)}$$
For the other two-dimensional lattices,  it has been shown that there exists $\delta >0$ such that
$$P_{0.5}({\cal Q}_4(m,n))\leq  \left({m\over n}\right)^{1+\delta}\mbox{ and } P_{0.5}({\cal Q}_1(m,n))\leq  \left({m\over n}\right)^{\delta}.\eqno{(1.15)}$$
In particular,  Kesten, Sidoravicius, and Zhang (1998) showed that (1.14) holds for $k=5$ for all the two-dimensional  lattices without  using an SLE approach.
Indeed, they showed that there exist $C_1$ and $C_2$ such that
$$C_1  \left({m\over n}\right)^{2}\leq P_{0.5}({\cal Q}_5(m,n))\leq C_2  \left({m\over n}\right)^{2}.\eqno{(1.16)}$$
In this paper, $C$ and $C_i$ are positive constants independent of  $n$, $m$, and $k$ and $p$. They also change from appearance to appearance. We sometimes use  $O(n)$ for $C_1n \leq O(n) \leq C_2 n$.
When we need to indicate that $C$ depends on a particular parameter, for example, $ \alpha$, we will write $C=C(\alpha)$.

It is interesting to point out that $(k^2-1)/12$ is  always a positive integer if $k\geq 5$ is a prime number. However, not only prime number $k$ satisfies that $(k^2-1)/12$ is an integer. For example, if $k=35$,  then $(k^2-1)/12$ is an integer. 
It might be easier to show (1.13) when $(k^2-1)/12$ is an integer as Kesten, Sidoravicius, and Zhang did in (1.16).  Note that if $k\geq 6$, then the critical exponent is larger than 2. Thus, we can barely see the six arm paths in $[-n, n]^2$ for large $n$. 
On the other hand,  two-arm path and three-arm path power laws in the half space are well understood. Let ${\cal H}_k(m,n)$ be the event that there are $i$-disjoint occupied paths and $j$ disjoint vacant paths in the upper half space with $i+j=k$ from $\partial B(m)$ to $\partial B(n)$
for $i, j \geq 1$.  It is believed that
$$P_{0.5}({\cal H}_k(m,n))= \left({m\over n}\right)^{k(k+1)/6+o(1)}.\eqno{(1.17)}$$
It has been proved (see Higuchi, Takei, and Zhang (2012)) that (1.17) holds
for $k=2$ and $k=3$ and for all the two-dimensional lattices. 
We know that (1.14) together with an argument of  Kesten (see Corollaries 1 and 2 in Kesten (1987)) imply that all the power laws and scaling  relations hold except possibly for (1.6).
 Here we will give the following theorem to discuss the power laws concerning for the free energy function.\\

{\bf Theorem.} {\em For percolation on the triangular lattice,
$$|\kappa'''(p)| \leq  |p-p_c|^{-1/3+o(1)} \mbox{ if } p \neq  p_c.$$
Furthermore, we show that there is a sequence  $\epsilon_n \downarrow 0$ such that
$$|\kappa'''(p_c\pm \epsilon_n)|\geq \epsilon_n^{-1/3+o(1)}.$$
}

{\bf Remarks.} 
1.  
For  percolation on the triangular  lattice,  the last inequality  of the Theorem  implies that $\kappa(p)$ is not third  differentiable. This  answers  affirmatively a  conjecture, asked by  Sykes and Essam in 1964, whether $\kappa(p)$ has a singularity at the criticality for the triangular lattice. 

 2.    
 (1.13) for $k=4$ was only  shown to hold for the triangular lattice.  If one can show (1.13) for the other lattices, then  the Theorem will work on  those lattices. In fact,   the topology of  the derivatives for $\kappa(p)$ in  the site  or bond percolation for the square lattice is easy to handle, since there are four neighbors for each vertex rather than six neighbors.  We are unable to prove the lower bound in the Theorem for any sequence $\epsilon_n\downarrow 0$.\\

{\bf Acknowledgments.} The author would like to thank Harry Kesten  for many fruitful conversations, and he would also like to thank Geoffrey Grimmett for his many comments.

\section{ Preliminaries.}  
In this section, we  introduce a few basic properties  and estimates of site percolation in the triangular lattice.  Most results are obtained from
Kesten (1982) and (1987).
 For any $u=(u_1, u_2)\in {\bf Z}^2$, let $\|u\|$ be the distance from the origin to $u$. For $u,v\in {\bf Z}^2$, let $d(u,v)=\|u-v\|$.
 For any two vertex sets $A$ and $B$, we denote  the {\em distance} of them  by
 $$d(A, B)= \min_{x\in A, y\in B} \|x-y\|.$$
 We can also define the distance of two edge sets if we consider the vertices of them.
 Given a finite connected graph $G$,  
 as we defined in section 1, a vertex $u\not\in G$ but  adjacent
to $G$ is called the {\em boundary vertex} of $G$. We denote by $\partial G$  the  boundary vertices of $G$.
If a vertex $v\in \partial G$, and there is an infinite path from $v$ without using $G$, then $v$ is called  the {\em exterior boundary} vertex. We denote by $\Delta G$  all the exterior boundary vertices of $G$.  
We introduce
a topology result (see Lemma 2.23 in Kesten (1982) or Proposition 11.2 in Grimmett (1999)).\\

{\bf Lemma 2.1.} {\em If $G$ is a finite cluster, then 
$\Delta  (G)$ is a circuit containing $G$ in its interior.  Furthermore,  $\Delta (G)$ is the smallest 
vacant circuit containing $G$ in its interior if $G$ is occupied.}\\

Now we define  occupied and  vacant crossings in a box. A left-right (respectively, top-bottom) occupied crossing of $B(n)$ is an occupied 
path  in $B(n)$ that joins some vertex on the left (respectively upper) side of $B(n)$ to some vertex on the right (respectively, lower)
 side of $B(n)$ but  uses no vertices  in the boundary of $B(n)$.
Similarly, we may define a vacant crossing on $B(n)$. 
We denote  the occupied and the vacant crossing probabilities of $B(n)$ by 
$$\sigma(p,n)= P_p(\exists \mbox{ a left-right occupied crossing of }B(n)) \mbox{ and }$$
$$\sigma^*(p,n)= P_p(\exists \mbox{ a top-bottom vacant crossing of }B(n)).$$
We need to show that the vacant
crossing and occupied crossing probabilities of squares 
are bounded away from zero when $p$ is
near $p_c$. To make this precise, we first define (see (1.21) in Kesten (1987))
\begin{eqnarray*}
L(p)=\min\{n: \sigma (p,n)\geq 
1- \epsilon_0\}\mbox{ for }  p> p_c\mbox{ and }
L(p)=\min\{n: \sigma^* (p,n)\geq 
\epsilon_0\}\mbox{ for } p< p_c, \hskip 1.5cm (2.1)
\end{eqnarray*}
where $\epsilon_0$ is some small but  strictly positive number whose precise
value is not important.
The important property is that $\epsilon_0$ can be chosen
such that there exists a constant $\delta$ for which
$$\sigma(p,n)\geq \delta \mbox{ and } \sigma^*(p,n)\geq \delta$$
uniformly in $n\leq L(p)$.
$L(p)$ is also called the correlation length, and it is proved
(see Corollary 2 in Kesten (1987)) that
$$ L(p)\asymp \xi(p),\eqno{(2.2)}$$
where $f(p) \asymp g(p) $ or $f_n \asymp g_n$ means
$$f(p)/g(p) =O(1) \mbox{ and }  f_n/g_n =O(1).$$
By (1.14) and Kesten's Corollaries 1 and 2, for all $p\neq 0.5$,
$$L(p)\asymp \xi(p)\approx |0.5-p|^{-4/3}.\eqno{(2.3)}$$
On the other hand,  it is known (see chapter 11 in Grimmett (1999)) that for any $n$
there exists $C >0$ such that
$$\sigma (0.5, n)\asymp \sigma^*(0.5, n)\asymp C.$$

If $p\neq 0.5$ and $n\leq L(p)$, Kesten's Theorem 1 and Lemma 8 (1987) show that $P_p({\cal Q}_l(n))$ 
has the same decay rate as 
$P_{0.5}({\cal Q}_l(n))$  for $l=1$ and $l=4$.
We want to remark that  two occupied paths are separated by two vacant paths in the four arm paths  of  ${\cal Q}_4(n)$ in Kesten's Theorem 1 and Lemma 8 (1987).  In other words, if we remove the origin, there is no occupied path 
connecting two occupied paths in ${\cal Q}_4(n)$. Later, we always consider these four arm paths. 
Together with his Lemmas 4 and  7 (1987), it shows that $P_p({\cal Q}_l(m,n))$ 
has the same decay rate as  $P_{0.5}({\cal Q}_l(m,n))$ for $l=1$ and $l=4$. 
Here we summarize his results as the following Lemma.\\

{\bf Lemma 2.2} (Kesten (1987)). {\em If $m\leq n \leq L(p)$,
$$P_{p}({\cal Q}_1(m,n))\asymp P_{0.5}({\cal Q}_1(m,n))\mbox{ and } P_{p}({\cal Q}_4(m,n))\asymp P_{0.5}({\cal Q}_4(m,n)).$$}



It is well known that the size of an occupied  cluster decays exponentially  (see Theorem 6.10 in Grimmett (1999)) when $p< p_c$: 
$$C_1 n^{-1} \exp\left(-\xi^{-1}(p)n\right)\leq P_p({\bf 0}\rightarrow \partial B(n))\leq C_2n \exp\left(-\xi^{-1}(p)n\right).$$
By using (2.1) and (2.2), we obtain   the following better estimate.\\

{\bf Lemma 2.3.} {\em If $p < p_c$, then there exist $C_i$  for $i=1,2$ such that
for any $n\geq L(p)$, }
$$C_1 {L(p)\over n}\exp\left(-L^{-1}(p)n\right)\leq P_p({\bf 0}\rightarrow \partial B(n))\leq C_2{n\over L(p)} \exp\left(-L^{-1}(p)n\right).$$

{\bf Proof.} Let $\beta(n)=P_p({\bf 0}\rightarrow \partial B(n))$. For $m \geq L(p)$, we may divide the parameter of 
$\partial B(m)$ into $m/L(p)$ segments such that each one has a length
$L(p)$. For each segment, by using (2.1) and the RSW lemma, there is an occupied circuit surrounding it with a positive probability.  Then we use the exact proof of (6.21) in Grimmett (1999) to show 
$$ \beta(n+m)\leq C_1(m/L(p)) \beta(m) \beta(n).\eqno{(2.4)}$$
The same method, together with (6.27) in Grimmett (1999), implies that
$$ \beta(n+m)\geq C_2(L(p)/m) \beta(m) \beta(n).\eqno{(2.5)}$$
So  together with (2.2), Lemma 2.3 follows from (2.4),(2.5), and the same proof of Theorem 5.10 in Grimmett (1999). $\Box$\\

By Lemmas 2.2 and 2.4, if $p \leq 0.5$, for all $m< n$, 
$$P_{p}({\cal Q}_4(m,n))\leq  \left({m\over n}\right)^{5/4+o(1)}.\eqno{(2.6)}$$
By (2.24) in Kesten (1987),  we have for a large $M$, 
$$P_p(B(L(p))\rightarrow \partial B(ML(p))) \leq C \exp(-C_1 M).\eqno{(2.7)}$$
By  using (2.7), we have the following lemma. \\

{\bf Lemma 2.4.} {\em If  $p< 0.5$, for $m \geq 1$,}
$$\sum_{i\geq mL(p)} iP_{p}(B(L(p))\rightarrow \partial B(i))\leq C L^2(p)\exp(-C_1m).\eqno{}$$

{\bf Proof.}  By (2.7),
\begin{eqnarray*}
&&\sum_{i\geq mL(p)} iP_{p}(B(L(p))\rightarrow \partial B(i))\leq \sum_{k=m}^\infty
\sum_{kL(p)\leq i\leq (k+1)L(p)}  iP_{p}(B(L(p))\rightarrow \partial B(i))\\
&\leq & 4C L^2(p)  \sum_{k=m}^\infty (k+1)^2 \exp(-C_1 k)\leq C_2 L^2(p)m^2 \exp(-C_1 m)\leq C_2 L^2(p) \exp(-C_3 m).
\end{eqnarray*}
Lemma 2.4 follows. $\Box$\\

Kesten in his Lemma 8 (1987) showed the following estimate:
$$L^2(p)P_p({\cal Q}_4(L(p)))\asymp (0.5-p)^{-1}.\eqno{(2.8)}$$
For convenience, we reselect $L(p)$ such that
$$L^2(p)P_p({\cal Q}_4(L(p)))= (0.5-p)^{-1}.\eqno{(2.9)}$$
In the following sections, we always assume that (2.9) holds for $L(p)$.
Given  four-arm paths in a square $B(n)$ and four-arm paths in an annulus  $A(n, m)$ for $n< m$, one of major estimates (see Lemma 4 and Lemma 6 in Kesten (1987)) is to reconnect them  by costing a constant probability independent of $n$ and $m$.  It is a more general argument than the RSW lemma since the RSW only reconnects the occupied or vacant paths. This connection  is called the {\em reconnection lemma}. \\

{\bf Reconnection lemma.} (Kesten (1987)). {\em If $n \leq L(p)$ and $n \leq m$ for $p \leq 0.5$,
there exists $C$ (independent of $n,m$) such that}
$$ P_p( {\cal Q}_4(n) )P_p({\cal Q}_4(n,m))\leq C P_p({\cal Q}_4(m)).$$

{\bf Remark} 3. The original proof of the reconnection lemma in Kesten's paper is involved and long. However, if we use the later developed three arm path estimate in (1.17), we should have an easier proof (see Morrow and Zhang (2005)). \\

Another lemma    (see  Kesten's Lemma 7  (2.58) (1987)) as the following lemma is useful when we estimate the pivotal vertices in  a square box.\\

{\bf Lemma 2.6.} (Kesten (1987)). {\em For $p < 0.5$ and for $j$ with $2^{j+1}\leq 2^k \leq L(p)$, there exists $\xi>0$ (independent of $p$) such that}
$$ 2^{2j} P_p(Q_4(2^j))\asymp \sum_{b\in B(2^{j-1})} P_p(Q_4(b, B(2^j)))\leq C(0.5-p)^{-1} 2^{-\xi(k-j)}.$$

With Lemma 2.4 and  Lemma 2.6,  we have the following lemma.\\

{\bf Lemma 2.7.}
{\em If $p < 0.5$, there exists $C$ (independent of $p$) such that
$$ \sum_{j=1} ^\infty jP_p(Q_4(j))\leq C(0.5-p)^{-1} .$$
Furthermore,  if $p < 0.5$, there exists $C_i$  for $i=1,2$ (independent of  $m$ and $p$) such that
$$ \sum_{j\geq m L(p)} ^\infty jP_p(Q_4(j))\leq C_1e^{-C_2 m} (0.5-p)^{-1} .$$}

\section{Derivative of $K_n(p)$ by finding its pivotal sites.}
The derivative of $K_n(p)$ in the triangular lattice is much more complicated than the derivative in the square lattice (see Zhang (2011)) since each vertex has six neighbors.
 Let us define for $v\in B(n)$, $N_n(v)$ as the number of distinct occupied clusters in $B(n)$, 
obtained after setting $v$ to be vacant, which contains a
neighbor of $v$. It follows from (5.5) in  Aizenman,  Kesten, and  Newman (1987) or Theorem 4.3 in Grimmett (1981) that
$${dK_n(p)\over dp}=\sum_{v\in B(n)} E_p(1-N_n(v)).\eqno{(3.1)}$$
Note that there is no infinite occupied cluster in $[0,0.5]$, so by (3.1) and a standard 
ergodic theorem (see Dunford and Schwartz Theorem VIII 6.9 (1958)),  
$$\kappa'(p)=\lim_{n\rightarrow \infty} |B(n)|^{-1} {dK_n(p)\over dp}=1-\lim_{n\rightarrow \infty}E_p(N_n(\bf 0)).\eqno{(3.2)}$$
Note also that $N_n(\bf 0)$ is neither increasing nor decreasing, so to take another derivative, we simply need to fix  the configurations of vertices adjacent to the origin.  
There are six vertices, denoted by $b_i$ for $i=1,2,3,4,5,6$, adjacent to the origin.
 We will try to fix the configurations of these six vertices.  Let ${\cal E}$ be the event that
 there are at least two occupied vertices in $\{b_1, b_2, b_3, b_4, b_5, b_6\}$,
 and they are separated by vacant vertices. For example, $b_1$ and $b_3$ are occupied and
 the others are vacant.  On ${\cal E}^C$, $N_n({\bf 0})$  (either zero or one)
 does not depend on the configurations outside of $\{{\bf 0}, b_1, b_2, b_3, b_4, b_5, b_6\}$.
 We may divide ${\cal E}^C$ into ${\cal E}^C_0$ or ${\cal E}^C_1$: all
 $\{b_1, b_2, b_3, b_4, b_5, b_6\}$ are vacant or only one connected occupied cluster is among them, respectively.
 Thus,
 $$E_p(N_n({\bf 0});{\cal E}^C)= P_p({\cal E}^C_1) =f(p)\eqno{(3.3)}$$
 for a polynomial $f(p)$ with a degree not more than $6$.
 
 We divide ${\cal E}$ into a few disjoint configurations. We denote by
 ${\cal E}_{AB}$ the event, for $A, B \subset \{b_1, b_2, b_3, b_4, b_5, b_6\}$, that 
 $A$ and $B$ are occupied  vertex sets separated by two vacant vertex sets $A'$ and $B'$ in  $\{b_1, b_2, b_3, b_4, b_5, b_6\}$.  In addition, $A, B, A', B'$ are adjacent individually
 in vertices $\{b_1, b_2, b_3, b_4, b_5, b_6\}$. 
 For example, $A=\{b_1\}$, $A'=\{b_2, b_3\}$, $B=\{b_4, b_5\}$, and $B'=\{b_6\}$.
 We also denote by ${\cal E}_{1,3,5}$ or ${\cal E}_{2,4, 6}$  the events that
 $b_1$, $b_3$, and $b_5$, or $b_2$, $b_4$, and $b_6$ are occupied and the others are
 vacant, respectively.
 Thus, by (3.3) and symmetry,
  \begin{eqnarray*}
&&E_p(N_n({\bf 0}))\\
&=& \sum_{A,B} E_p(N_n({\bf 0})\,\,| \,\, {\cal E}_{AB})f_{AB}(p)+ 
 E_p(N_n({\bf 0})\,\,| \,\, {\cal E}_{1,3,5})f_{1,3,5}(p)+ E_p(N_n({\bf 0})\,\,| \,\, {\cal E}_{2,4,6})f_{2,4,6}(p)+f(p)\\
&&=  \sum_{A,B} E_p(N_n({\bf 0})\,\,| \,\, {\cal E}_{AB})f_{AB}(p)+ 
 2E_p(N_n({\bf 0})\,\,| \,\, {\cal E}_{1,3,5})f_{1,3,5}(p)+f(p),
\end{eqnarray*}
 where the first sum takes all possible $A$ and $B$ and $f_{AB}(p)=P_p({\cal E}_{AB})$, $f_{1,3,5}(p)=P_p({\cal E}_{1,3,5})$,
 and $f_{2,4,6}(p)=P_p({\cal E}_{2,4,6})$ are polynomials with  degrees less than 6. Therefore, 
 \begin{eqnarray*}
 {dE_p(N_n({\bf 0}))\over dp}
 &=&\sum_{A, B} {d E_p(N_n({\bf 0})\,\,| \,\, {\cal E}_{AB})\over dp} f_{AB}(p)+\sum_{A,B} E_p(N_n({\bf 0})\,\,| \,\, {\cal E}_{AB}) f_{AB}'(p)\\
 &+&2{d E_p(N_n({\bf 0})\,\,| \,\, {\cal E}_{1,3,5})\over dp} f_{1,3,5}(p)+  2E_p(N_n({\bf 0})\,\,| \,\, {\cal E}_{1,3,5}) f_{1,3,5}'(p) +f'(p).
 \hskip 2cm (3.4)
 \end{eqnarray*}
After taking this derivative, the first sum and  the third  term in the right side of (3.4)  play the most important roles.  Let us focus on a term in the first sum in the right side of (3.4):
$dE_p(N_n({\bf 0})\,\,| \,\, {\cal E}_{AB})/ dp$ for  fixed sets $A$ and $B$.
Note that 
$$E_p(N_n({\bf 0})\,\,| \,\, {\cal E}_{AB})= P_p( N_n({\bf 0}) \geq 2\,\,| \,\, {\cal E}_{AB})+P_p( N_n({\bf 0}) \geq 1\,\,| \,\, {\cal E}_{AB})=P_p( N_n({\bf 0} )= 2\,\,| \,\, {\cal E}_{AB})+1$$
and $\{ N_n({\bf 0}) \geq 2\,\,| \,\, {\cal E}_{AB}\}=\{ N_n({\bf 0}) = 2\,\,| \,\, {\cal E}_{AB}\}$ is a decreasing event
for the configurations on $B(n)$ except for the origin and its neighbors.
By Russo's formula, if  $\{N_n({\bf 0}) =2\,\,|\,\, {\cal E}_{AB} \}(b)$ is 
the event that $b$ is a pivotal vertex (see Fig. 1.1 and the definition for a pivotal vertex in Kesten (1982)), 
then
\begin{eqnarray*}
&&{dE_p(N_n({\bf 0} )\,\,| \,\, {\cal E}_{AB})\over dp}=-  \sum_{b\in B(n),b\not\sim {\bf 0}} P_p( \{N_n({\bf 0})=2
 \,\,|\,\, {\cal E}_{AB}\}(b))\\
&=&  -  \sum_{b\in B(n), b\not\sim {\bf 0}} P_p\left( \{N_n({\bf 0})=1
\mbox{ or } 2 \mbox{ if $b$ is occupied or vacant}\}\,\,|\,\, {\cal E}_{AB}\right)  ,\hskip 4.3cm{(3.5)}
\end{eqnarray*}
where $b\not\sim {\bf 0}$ means that $b\neq {\bf 0}$ is not adjacent to ${\bf 0}$.
Let ${\bf C}_n(A)$ be the occupied cluster by assuming $A$ is occupied in 
$B(n) \setminus \{{\bf 0}\}$ on ${\cal E}_{AB}$, and let $\Delta {\bf C}_n(A)$
be the sites of the exterior boundary of ${\bf C}_n(A)$.
We can rewrite $ \{N_n({\bf 0})=1
\mbox{ or } =2 \mbox{ if $b$ is occupied or vacant}\,\,|\,\, {\cal E}_{AB}\}$ as 
$${\cal D}_{AB}(b, n)=\{ {\bf C}_n(A)\cap {\bf C}_n(B) =\emptyset, b\in\Delta {\bf C}_n(A)\cap
 \Delta {\bf C}_n(B)\}.\eqno{(3.6)}$$
 Now we need to write ${\cal D}_{AB}(b, n)$ into occupied and vacant paths. 
Let  (see Fig. 1.5 and Fig. 1.6)
\begin{eqnarray*}
&&\bar{{\cal R}}_{AB}(b,n)=\\
&&\{\exists\,\,\mbox{ disjoint occupied  paths $r_2$ and $r_4$ in $B(n)\setminus {\bf 0}$ from 
 $A$ and $B$ to two neighbors of $b$; }\\
 && \,\,\,\exists \mbox{ a vacant $r_3$ on $B(n)\setminus {\bf 0}$ inside $S(r_2, r_4, b)$
from $A'$ or $B'$ to a neighbor of $b$;  $\exists$ a vacant path}\\
&&\mbox{  $r_1$ inside $B(n)\setminus {\bf 0}$ from  $B'$ or $A'$ to a neighbor
  of $b$ but outside the closure of  $S( r_2, r_4, b)$ or }\\
&&\,\,\,\exists \mbox{ disjoint vacant paths $r_5$ and $r_6$ inside $B(n)\setminus {\bf 0}$ from $B'$  and
from a neighbor of $b$ to $\partial B(n)$,  }\\
&&\mbox{ but outside the closure of  $S( r_2, r_4, b)$, and $\exists$ an occupied path from $r_2$ or $r_4$ to $\partial B(n)$}
\},
\end{eqnarray*}
where $S(r_2, r_4, b)$ is the open set enclosed by the circuit $r_2\cup {\bf 0}\cup r_4\cup b$ and a path from $A$ means
the path starting from a neighbor of  $A$.
\begin{figure}
\begin{center}
\setlength{\unitlength}{0.0125in}%
\begin{picture}(200,500)(67,660)
\thicklines
\put(-60,770){\framebox(200,200)[br]{\mbox{$[-n,n]^2$}}}
\put(220,770){\framebox(200,200)[br]{\mbox{$[-n,n]^2$}}}
\put(94,902){\circle*{7}}

\put(85,890){\mbox{$b$}{}}
\put(87,895){*}
\put(82,895){*}
\put(76,895){*}
\put(71,895){*}
\put(69,892){*}
\put(50,882){\mbox{$r_4$}}

\put(65,890){*}
\put(60,890){*}
\put(55,890){*}
\put(50,890){*}
\put(45,885){*}
\put(40,880){*}
\put(35,875){*}
\put(30,870){*}
\put(25,865){*}
\put(20,860){*}
\put(15,855){*}

\put(95,907){\circle*{2}}
\put(95,910){\circle*{2}}
\put(90,910){\circle*{2}}
\put(85,910){\circle*{2}}
\put(80,910){\circle*{2}}
\put(75,910){\circle*{2}}
\put(70,910){\circle*{2}}
\put(65,910){\circle*{2}}
\put(60,910){\circle*{2}}
\put(55,910){\circle*{2}}
\put(50,910){\circle*{2}}
\put(45,910){\circle*{2}}
\put(40,910){\circle*{2}}
\put(35,910){\circle*{2}}
\put(30,910){\circle*{2}}
\put(25,910){\circle*{2}}
\put(20,910){\circle*{2}}
\put(15,910){\circle*{2}}
\put(15,905){\circle*{2}}
\put(15,900){\circle*{2}}
\put(15,895){\circle*{2}}
\put(15,890){\circle*{2}}
\put(15,885){\circle*{2}}
\put(15,880){\circle*{2}}
\put(15,895){\circle*{2}}
\put(15,890){\circle*{2}}
\put(15,885){\circle*{2}}
\put(15,880){\circle*{2}}
\put(15,875){\circle*{2}}
\put(15,870){\circle*{2}}
\put(15,865){\circle*{2}}
\put(15,860){\circle*{2}}
\put(15,855){\circle*{2}}
\put(35,900){\mbox{$r_3$}}

\put(20,855){\circle*{2}}
\put(25,855){\circle*{2}}
\put(30,855){\circle*{2}}
\put(35,855){\circle*{2}}
\put(40,855){\circle*{2}}
\put(45,855){\circle*{2}}
\put(50,855){\circle*{2}}
\put(55,855){\circle*{2}}
\put(60,855){\circle*{2}}
\put(65,855){\circle*{2}}
\put(70,855){\circle*{2}}
\put(75,855){\circle*{2}}
\put(80,855){\circle*{2}}
\put(85,855){\circle*{2}}
\put(90,855){\circle*{2}}
\put(95,855){\circle*{2}}
\put(95,860){\circle*{2}}
\put(95,865){\circle*{2}}
\put(95,870){\circle*{2}}
\put(95,875){\circle*{2}}
\put(95,880){\circle*{2}}
\put(95,885){\circle*{2}}
\put(95,890){\circle*{2}}
\put(95,895){\circle*{2}}
\put(95,900){\circle*{2}}
\put(80,845){\mbox{$r_1$}}

\put(98,897){*}

\put(100,900){*}
\put(105,900){*}
\put(110,900){*}
\put(110,905){*}
\put(110,910){*}
\put(110,915){*}
\put(110,920){*}
\put(105,920){*}
\put(100,920){*}
\put(95,920){*}
\put(90,920){*}
\put(85,920){*}
\put(80,920){*}
\put(75,920){*}
\put(70,920){*}
\put(65,920){*}
\put(60,920){*}
\put(55,920){*}
\put(50,920){*}
\put(45,920){*}
\put(40,920){*}
\put(35,920){*}
\put(30,920){*}
\put(25,920){*}
\put(20,920){*}
\put(15,920){*}
\put(10,920){*}
\put(5,920){*}
\put(0,920){*}
\put(-5,920){*}
\put(-10,920){*}
\put(-15,920){*}
\put(-20,920){*}
\put(-20,915){*}
\put(-20,910){*}
\put(-20,905){*}
\put(-20,900){*}
\put(-18,896){*}
\put(-15,892){*}
\put(-12,888){*}
\put(-9,884){*}
\put(-5,880){*}
\put(-5,875){*}
\put(-5,870){*}
\put(-5,865){*}
\put(-5,860){*}
\put(-5,855){*}
\put(0,855){*}
\put(5,855){*}
\put(-20,887){\mbox{$r_2$}}
\put(20,745){\mbox{$\mbox{Fig. 1.5}$}}

\put(14,859){\circle*{7}}
\put(12,845){\mbox{${\bf 0}$}{}}
\put(270,850){\mbox{${\bf 0}$}{}}
\put(283,863){\circle*{7}}

\put(383,898){\circle*{7}}
\put(287,865){\circle*{2}}
\put(287,865){\circle*{2}}
\put(290,865){\circle*{2}}
\put(295,865){\circle*{2}}
\put(300,865){\circle*{2}}
\put(305,865){\circle*{2}}
\put(310,865){\circle*{2}}
\put(315,865){\circle*{2}}
\put(320,865){\circle*{2}}
\put(325,865){\circle*{2}}
\put(330,865){\circle*{2}}
\put(335,865){\circle*{2}}
\put(340,865){\circle*{2}}
\put(340,870){\circle*{2}}
\put(350,875){\mbox{$r_3$}}
\put(340,875){\circle*{2}}
\put(340,880){\circle*{2}}
\put(340,885){\circle*{2}}
\put(340,890){\circle*{2}}
\put(340,895){\circle*{2}}
\put(340,900){\circle*{2}}
\put(345,900){\circle*{2}}
\put(350,900){\circle*{2}}
\put(355,900){\circle*{2}}
\put(360,900){\circle*{2}}
\put(365,900){\circle*{2}}
\put(370,900){\circle*{2}}
\put(375,900){\circle*{2}}
\put(390,900){\circle*{2}}
\put(395,900){\circle*{2}}
\put(400,900){\circle*{2}}
\put(405,900){\circle*{2}}
\put(410,900){\circle*{2}}
\put(415,900){\circle*{2}}
\put(420,900){\circle*{2}}

\put(400,885){\mbox{$r_6$}}

\put(275,865){\circle*{2}}
\put(270,865){\circle*{2}}
\put(265,865){\circle*{2}}
\put(260,865){\circle*{2}}
\put(240,870){\mbox{$r_5$}}
\put(255,865){\circle*{2}}
\put(250,865){\circle*{2}}
\put(245,865){\circle*{2}}
\put(240,865){\circle*{2}}
\put(235,865){\circle*{2}}
\put(230,865){\circle*{2}}
\put(225,865){\circle*{2}}
\put(220,865){\circle*{2}}

\put(280,865){*}
\put(280,870){*}
\put(280,875){*}
\put(280,880){*}
\put(280,885){*}
\put(280,890){*}
\put(280,895){*}
\put(280,900){*}
\put(280,905){*}
\put(280,910){*}
\put(280,915){*}
\put(280,920){*}
\put(280,925){*}
\put(290,925){\mbox{$r_2$}}
\put(290,745){\mbox{$\mbox{ Fig. 1.6}$}}

\put(280,930){*}
\put(280,935){*}
\put(285,935){*}
\put(290,935){*}
\put(295,935){*}
\put(300,935){*}
\put(305,935){*}
\put(310,935){*}
\put(315,935){*}
\put(320,935){*}
\put(325,935){*}
\put(330,935){*}
\put(335,935){*}
\put(340,935){*}
\put(345,935){*}
\put(350,935){*}
\put(355,935){*}
\put(360,935){*}
\put(365,935){*}
\put(370,935){*}
\put(375,935){*}
\put(380,935){*}
\put(380,930){*}
\put(380,925){*}
\put(380,920){*}
\put(380,915){*}
\put(380,910){*}
\put(380,905){*}

\put(380,900){*}

\put(380,890){*}
\put(380,885){*}
\put(380,880){*}
\put(380,875){*}
\put(380,870){*}
\put(380,865){*}
\put(380,860){*}
\put(380,855){*}
\put(380,850){*}
\put(375,850){*}
\put(370,850){*}
\put(365,850){*}
\put(360,850){*}
\put(355,850){*}
\put(350,850){*}
\put(345,850){*}
\put(340,850){*}
\put(340,845){\mbox{$r_4$}}

\put(335,850){*}
\put(330,850){*}
\put(325,850){*}
\put(320,850){*}
\put(315,850){*}
\put(310,850){*}
\put(305,850){*}
\put(300,850){*}
\put(294,850){*}
\put(290,850){*}
\put(285,850){*}
\put(280,850){*}
\put(375,883){\mbox{$b$}{}}
\put(385,860){*}
\put(390,860){*}
\put(395,860){*}
\put(400,860){*}
\put(405,860){*}
\put(410,860){*}
\put(415,860){*}

\put(-80, 1050){\line(1,0){20}}
\put(-80, 1050){\line(0,1){20}}
\put(-80, 1070){\line(1,1){20}}
\put(-60, 1050){\line(1,1){20}}
\put(-60, 1090){\line(1,0){20}}
\put(-40, 1090){\line(0,-1){20}}
\put(-60,1070){\circle*{6}}
\put(-70,1060){\mbox{${\bf 0}$}}
\put(-80, 1200){\line(1,0){20}}
\put(-80, 1200){\line(0,1){20}}
\put(-80, 1220){\line(1,1){20}}
\put(-60, 1200){\line(1,1){20}}
\put(-60, 1240){\line(1,0){20}}
\put(-40, 1240){\line(0,-1){20}}
\put(-60,1220){\circle*{6}}
\put(-70,1210){\mbox{$b$}}
\put(-85, 1039){\huge{*}}
\put(-65, 1039){\huge{*}}
\put(-65, 1079){\huge{*}}
\put(-80, 1070){\circle*{3}}
\put(-40, 1070){\circle*{3}}
\put(-40, 1090){\circle*{3}}
\put(-80, 1080){\circle*{3}}
\put(-80, 1090){\circle*{3}}
\put(-80, 1100){\circle*{3}}
\put(-80, 1110){\circle*{3}}
\put(-80, 1120){\circle*{3}}
\put(-80, 1130){\circle*{3}}
\put(-80, 1140){\circle*{3}}
\put(-80, 1150){\circle*{3}}
\put(-80, 1160){\circle*{3}}
\put(-90, 1170){\circle*{3}}
\put(-90, 1180){\circle*{3}}
\put(-90, 1190){\circle*{3}}
\put(-90, 1200){\circle*{3}}
\put(-90, 1210){\circle*{3}}
\put(-90, 1200){\circle*{3}}
\put(-90, 1210){\circle*{3}}
\put(-80, 1220){\circle*{3}}
\put(-65, 1089){\huge{*}}
\put(-65, 1099){\huge{*}}
\put(-65, 1109){\huge{*}}
\put(-65, 1119){\huge{*}}
\put(-65, 1129){\huge{*}}
\put(-65, 1139){\huge{*}}
\put(-65, 1149){\huge{*}}
\put(-65, 1159){\huge{*}}
\put(-65, 1169){\huge{*}}
\put(-65, 1179){\huge{*}}
\put(-65, 1189){\huge{*}}
\put(-40, 1100){\circle*{3}}
\put(-40, 1110){\circle*{3}}
\put(-40, 1120){\circle*{3}}
\put(-40, 1130){\circle*{3}}
\put(-40, 1140){\circle*{3}}
\put(-40, 1150){\circle*{3}}
\put(-40, 1160){\circle*{3}}
\put(-40, 1170){\circle*{3}}
\put(-40, 1180){\circle*{3}}
\put(-40, 1190){\circle*{3}}
\put(-40, 1200){\circle*{3}}
\put(-40, 1210){\circle*{3}}
\put(-40, 1220){\circle*{3}}
\put(-55, 1039){\huge{*}}
\put(-45, 1039){\huge{*}}
\put(-35, 1039){\huge{*}}
\put(-25, 1049){\huge{*}}
\put(-25, 1059){\huge{*}}
\put(-25, 1069){\huge{*}}
\put(-25, 1079){\huge{*}}
\put(-25, 1089){\huge{*}}
\put(-25, 1099){\huge{*}}
\put(-25, 1109){\huge{*}}
\put(-25, 1119){\huge{*}}
\put(-25, 1129){\huge{*}}
\put(-25, 1139){\huge{*}}
\put(-25, 1149){\huge{*}}
\put(-25, 1159){\huge{*}}
\put(-25, 1169){\huge{*}}
\put(-25, 1179){\huge{*}}
\put(-25, 1189){\huge{*}}
\put(-25, 1199){\huge{*}}
\put(-25, 1209){\huge{*}}
\put(-25, 1219){\huge{*}}
\put(-25, 1229){\huge{*}}
\put(-35, 1229){\huge{*}}
\put(-45, 1229){\huge{*}}
\put(-60, 1010){$\mbox{ Fig. 1.1}$}

\put(80, 1050){\line(1,0){20}}
\put(80, 1050){\line(0,1){20}}
\put(80, 1070){\line(1,1){20}}
\put(100, 1050){\line(1,1){20}}
\put(100, 1090){\line(1,0){20}}
\put(120, 1090){\line(0,-1){20}}
\put(100,1070){\circle*{6}}
\put(90,1060){\mbox{${\bf 0}$}}
\put(80, 1200){\line(1,0){20}}
\put(80, 1200){\line(0,1){20}}
\put(80, 1220){\line(1,1){20}}
\put(100, 1200){\line(1,1){20}}
\put(100, 1240){\line(1,0){20}}
\put(120, 1240){\line(0,-1){20}}
\put(100,1220){\circle*{6}}
\put(90,1210){\mbox{$b$}}
\put(80, 1050){\circle*{3}}
\put(76, 1060){\huge{*}}
\put(96, 1040){\huge{*}}
\put(116, 1080){\huge{*}}
\put(120, 1070){\circle*{3}}
\put(100, 1090){\circle*{3}}
\put(100, 1100){\circle*{3}}
\put(100, 1110){\circle*{3}}
\put(100, 1120){\circle*{3}}
\put(100, 1130){\circle*{3}}
\put(100, 1140){\circle*{3}}
\put(100, 1150){\circle*{3}}
\put(100, 1160){\circle*{3}}
\put(100, 1170){\circle*{3}}
\put(100, 1180){\circle*{3}}
\put(100, 1190){\circle*{3}}
\put(100, 1200){\circle*{3}}
\put(116, 1090){\huge{*}}
\put(116, 1100){\huge{*}}
\put(116, 1110){\huge{*}}
\put(116, 1100){\huge{*}}
\put(116, 1120){\huge{*}}
\put(116, 1130){\huge{*}}
\put(116, 1140){\huge{*}}
\put(116, 1150){\huge{*}}
\put(116, 1160){\huge{*}}
\put(116, 1170){\huge{*}}
\put(116, 1180){\huge{*}}
\put(116, 1190){\huge{*}}
\put(116, 1200){\huge{*}}
\put(116, 1210){\huge{*}}
\put(76, 1090){\huge{*}}
\put(76, 1100){\huge{*}}
\put(76, 1110){\huge{*}}
\put(76, 1100){\huge{*}}
\put(76, 1120){\huge{*}}
\put(76, 1130){\huge{*}}
\put(76, 1140){\huge{*}}
\put(76, 1150){\huge{*}}
\put(76, 1160){\huge{*}}
\put(76, 1170){\huge{*}}
\put(76, 1180){\huge{*}}
\put(76, 1190){\huge{*}}
\put(76, 1060){\huge{*}}
\put(76, 1070){\huge{*}}
\put(76, 1080){\huge{*}}
\put(130, 1070){\circle*{3}}
\put(140, 1070){\circle*{3}}
\put(140, 1080){\circle*{3}}
\put(140, 1080){\circle*{3}}
\put(140, 1090){\circle*{3}}
\put(140, 1100){\circle*{3}}
\put(140, 1110){\circle*{3}}
\put(140, 1120){\circle*{3}}
\put(140, 1130){\circle*{3}}
\put(140, 1070){\circle*{3}}
\put(140, 1080){\circle*{3}}
\put(140, 1080){\circle*{3}}
\put(140, 1090){\circle*{3}}
\put(140, 1100){\circle*{3}}
\put(140, 1110){\circle*{3}}
\put(140, 1120){\circle*{3}}
\put(140, 1130){\circle*{3}}
\put(140, 1140){\circle*{3}}
\put(140, 1150){\circle*{3}}
\put(140, 1160){\circle*{3}}
\put(140, 1170){\circle*{3}}
\put(140, 1180){\circle*{3}}
\put(140, 1190){\circle*{3}}
\put(140, 1200){\circle*{3}}
\put(140, 1210){\circle*{3}}
\put(140, 1220){\circle*{3}}
\put(140, 1230){\circle*{3}}
\put(140, 1240){\circle*{3}}
\put(130, 1240){\circle*{3}}
\put(120, 1240){\circle*{3}}
\put(120, 1040){\circle*{3}}
\put(110, 1040){\circle*{3}}
\put(100, 1040){\circle*{3}}
\put(90, 1040){\circle*{3}}
\put(80, 1040){\circle*{3}}
\put(130, 1050){\circle*{3}}
\put(140, 1060){\circle*{3}}
\put(90, 1010){$\mbox{ Fig. 1.2}$}

\put(220, 1050){\line(1,0){20}}
\put(220, 1050){\line(0,1){20}}
\put(220, 1070){\line(1,1){20}}
\put(240, 1050){\line(1,1){20}}
\put(240, 1090){\line(1,0){20}}
\put(260, 1090){\line(0,-1){20}}
\put(240,1070){\circle*{6}}
\put(230,1060){\mbox{${\bf 0}$}}
\put(220, 1200){\line(1,0){20}}
\put(220, 1200){\line(0,1){20}}
\put(220, 1220){\line(1,1){20}}
\put(240, 1200){\line(1,1){20}}
\put(240, 1240){\line(1,0){20}}
\put(260, 1240){\line(0,-1){20}}
\put(240,1220){\circle*{6}}
\put(230,1210){\mbox{$b$}}
\put(220, 1050){\circle*{3}}
\put(216, 1060){\huge{*}}
\put(236, 1040){\huge{*}}
\put(256, 1080){\huge{*}}
\put(260, 1070){\circle*{3}}
\put(240, 1090){\circle*{3}}
\put(116, 1080){\huge{*}}

\put(240, 1100){\circle*{3}}
\put(240, 1090){\circle*{3}}
\put(240, 1100){\circle*{3}}
\put(240, 1110){\circle*{3}}
\put(240, 1120){\circle*{3}}
\put(240, 1130){\circle*{3}}
\put(240, 1140){\circle*{3}}
\put(240, 1150){\circle*{3}}
\put(240, 1160){\circle*{3}}
\put(240, 1170){\circle*{3}}
\put(240, 1180){\circle*{3}}
\put(240, 1190){\circle*{3}}
\put(240, 1200){\circle*{3}}
\put(256, 1090){\huge{*}}
\put(256, 1100){\huge{*}}
\put(256, 1110){\huge{*}}
\put(256, 1100){\huge{*}}
\put(256, 1120){\huge{*}}
\put(256, 1130){\huge{*}}
\put(256, 1140){\huge{*}}
\put(256, 1150){\huge{*}}
\put(256, 1160){\huge{*}}
\put(256, 1170){\huge{*}}
\put(256, 1180){\huge{*}}
\put(256, 1190){\huge{*}}
\put(256, 1200){\huge{*}}
\put(256, 1210){\huge{*}}
\put(216, 1090){\huge{*}}
\put(216, 1100){\huge{*}}
\put(216, 1110){\huge{*}}
\put(216, 1100){\huge{*}}
\put(216, 1120){\huge{*}}
\put(216, 1130){\huge{*}}
\put(216, 1140){\huge{*}}
\put(216, 1150){\huge{*}}
\put(216, 1160){\huge{*}}
\put(216, 1170){\huge{*}}
\put(216, 1180){\huge{*}}
\put(216, 1190){\huge{*}}
\put(216, 1060){\huge{*}}
\put(216, 1070){\huge{*}}
\put(216, 1080){\huge{*}}
\put(270, 1070){\circle*{3}}
\put(280, 1070){\circle*{3}}
\put(280, 1080){\circle*{3}}
\put(280, 1080){\circle*{3}}
\put(280, 1090){\circle*{3}}
\put(280, 1100){\circle*{3}}
\put(280, 1110){\circle*{3}}
\put(280, 1120){\circle*{3}}
\put(280, 1130){\circle*{3}}
\put(280, 1070){\circle*{3}}
\put(280, 1080){\circle*{3}}
\put(280, 1080){\circle*{3}}
\put(280, 1090){\circle*{3}}
\put(280, 1100){\circle*{3}}
\put(280, 1110){\circle*{3}}
\put(280, 1120){\circle*{3}}
\put(280, 1130){\circle*{3}}
\put(280, 1140){\circle*{3}}
\put(280, 1150){\circle*{3}}
\put(280, 1160){\circle*{3}}
\put(280, 1170){\circle*{3}}
\put(280, 1180){\circle*{3}}
\put(280, 1190){\circle*{3}}
\put(280, 1200){\circle*{3}}
\put(280, 1210){\circle*{3}}
\put(280, 1220){\circle*{3}}
\put(280, 1230){\circle*{3}}
\put(280, 1240){\circle*{3}}
\put(270, 1240){\circle*{3}}
\put(260, 1240){\circle*{3}}
\put(236, 1030){\huge{*}}
\put(226, 1030){\huge{*}}
\put(216, 1030){\huge{*}}
\put(206, 1030){\huge{*}}
\put(206, 1040){\huge{*}}
\put(206, 1050){\huge{*}}
\put(216, 1010){$\mbox{ Fig. 1.3}$}

\put(350, 1050){\line(1,0){20}}
\put(350, 1050){\line(0,1){20}}
\put(350, 1070){\line(1,1){20}}
\put(370, 1050){\line(1,1){20}}
\put(370, 1090){\line(1,0){20}}
\put(390, 1090){\line(0,-1){20}}
\put(370,1070){\circle*{6}}
\put(360,1060){\mbox{${\bf 0}$}}
\put(350, 1200){\line(1,0){20}}
\put(350, 1200){\line(0,1){20}}
\put(350, 1220){\line(1,1){20}}
\put(370, 1200){\line(1,1){20}}
\put(370, 1240){\line(1,0){20}}
\put(390, 1240){\line(0,-1){20}}
\put(370,1220){\circle*{6}}
\put(360,1210){\mbox{$b$}}
\put(350, 1050){\circle*{3}}
\put(346, 1060){\huge{*}}
\put(366, 1040){\huge{*}}
\put(386, 1080){\huge{*}}
\put(390, 1070){\circle*{3}}
\put(370, 1090){\circle*{3}}

\put(370, 1100){\circle*{3}}
\put(370, 1090){\circle*{3}}
\put(370, 1100){\circle*{3}}
\put(370, 1110){\circle*{3}}
\put(370, 1120){\circle*{3}}
\put(370, 1130){\circle*{3}}
\put(370, 1140){\circle*{3}}
\put(370, 1150){\circle*{3}}
\put(370, 1160){\circle*{3}}
\put(370, 1170){\circle*{3}}
\put(370, 1180){\circle*{3}}
\put(370, 1190){\circle*{3}}
\put(370, 1200){\circle*{3}}
\put(386, 1090){\huge{*}}
\put(386, 1100){\huge{*}}
\put(386, 1110){\huge{*}}
\put(386, 1100){\huge{*}}
\put(386, 1120){\huge{*}}
\put(386, 1130){\huge{*}}
\put(386, 1140){\huge{*}}
\put(386, 1150){\huge{*}}
\put(386, 1160){\huge{*}}
\put(386, 1170){\huge{*}}
\put(386, 1180){\huge{*}}
\put(386, 1190){\huge{*}}
\put(386, 1200){\huge{*}}
\put(386, 1210){\huge{*}}
\put(346, 1090){\huge{*}}
\put(346, 1100){\huge{*}}
\put(346, 1110){\huge{*}}
\put(346, 1100){\huge{*}}
\put(346, 1120){\huge{*}}
\put(346, 1130){\huge{*}}
\put(346, 1140){\huge{*}}
\put(346, 1150){\huge{*}}
\put(346, 1160){\huge{*}}
\put(346, 1170){\huge{*}}
\put(346, 1180){\huge{*}}
\put(346, 1190){\huge{*}}
\put(346, 1060){\huge{*}}
\put(346, 1070){\huge{*}}
\put(346, 1080){\huge{*}}
\put(400, 1070){\circle*{3}}
\put(410, 1070){\circle*{3}}
\put(410, 1080){\circle*{3}}
\put(410, 1080){\circle*{3}}
\put(410, 1090){\circle*{3}}
\put(410, 1100){\circle*{3}}
\put(410, 1110){\circle*{3}}
\put(410, 1120){\circle*{3}}
\put(410, 1130){\circle*{3}}
\put(410, 1070){\circle*{3}}
\put(410, 1080){\circle*{3}}
\put(410, 1080){\circle*{3}}
\put(410, 1090){\circle*{3}}
\put(410, 1100){\circle*{3}}
\put(410, 1110){\circle*{3}}
\put(410, 1120){\circle*{3}}
\put(410, 1130){\circle*{3}}
\put(410, 1140){\circle*{3}}
\put(410, 1150){\circle*{3}}
\put(410, 1160){\circle*{3}}
\put(410, 1170){\circle*{3}}
\put(410, 1180){\circle*{3}}
\put(410, 1190){\circle*{3}}
\put(410, 1200){\circle*{3}}
\put(410, 1210){\circle*{3}}
\put(410, 1220){\circle*{3}}
\put(410, 1230){\circle*{3}}
\put(410, 1240){\circle*{3}}
\put(400, 1240){\circle*{3}}
\put(390, 1240){\circle*{3}}

\put(330, 1080){\circle*{3}}
\put(330, 1090){\circle*{3}}
\put(330, 1100){\circle*{3}}
\put(330, 1110){\circle*{3}}
\put(330, 1120){\circle*{3}}
\put(330, 1130){\circle*{3}}
\put(330, 1140){\circle*{3}}
\put(330, 1150){\circle*{3}}
\put(330, 1160){\circle*{3}}
\put(330, 1170){\circle*{3}}
\put(330, 1180){\circle*{3}}
\put(330, 1190){\circle*{3}}
\put(330, 1200){\circle*{3}}
\put(330, 1210){\circle*{3}}
\put(330, 1200){\circle*{3}}
\put(330, 1210){\circle*{3}}
\put(340, 1220){\circle*{3}}
\put(350, 1220){\circle*{3}}
\put(330, 1070){\circle*{3}}
\put(340, 1060){\circle*{3}}
\put(376, 1040){\huge{*}}
\put(386, 1040){\huge{*}}
\put(396, 1040){\huge{*}}
\put(406, 1040){\huge{*}}
\put(416, 1040){\huge{*}}
\put(426, 1040){\huge{*}}
\put(426, 1050){\huge{*}}
\put(426, 1060){\huge{*}}
\put(426, 1070){\huge{*}}
\put(426, 1080){\huge{*}}
\put(426, 1090){\huge{*}}
\put(426, 1100){\huge{*}}
\put(426, 1110){\huge{*}}
\put(426, 1120){\huge{*}}
\put(426, 1130){\huge{*}}
\put(426, 1140){\huge{*}}
\put(426, 1150){\huge{*}}
\put(426, 1160){\huge{*}}
\put(426, 1170){\huge{*}}
\put(426, 1180){\huge{*}}
\put(426, 1190){\huge{*}}
\put(426, 1200){\huge{*}}
\put(426, 1210){\huge{*}}
\put(426, 1220){\huge{*}}
\put(426, 1230){\huge{*}}
\put(426, 1240){\huge{*}}
\put(416, 1240){\huge{*}}
\put(406, 1240){\huge{*}}
\put(396, 1240){\huge{*}}
\put(386, 1240){\huge{*}}
\put(376, 1240){\huge{*}}
\put(366, 1240){\huge{*}}
\put(366, 1230){\huge{*}}
\put(370, 1010){$\mbox{ Fig. 1.4}$}
\thicklines
\end{picture}
\end{center}
\caption{ \em Fig. 1.5 is event ${\cal R}(b, n)$. Fig. 1.5 and  Fig. 1.6 are the two situations of $\bar{\cal R}(b, n)$.
In Fig. 1.6, since $r_1$ does not exist, we have an occupied path separating $r_5$ and $r_6$. Fig. 1.1  shows that $b$ is a pivotal vertex for $P_p(N({\bf 0})=2\,\,|\,\, {\cal E}_{AB})$. Fig. 1.2 and Fig. 1.4 show that $b$ is a pivotal vertex for 
$P_p(N({\bf 0})=3\,\,|\,\, {\cal E}_{1, 3, 5})$. Figs. 1.2, 1.3, and Fig. 1.4 show that
 $b$ is a pivotal vertex for 
$P_p(N({\bf 0})\geq 2\,\,|\,\, {\cal E}_{1,3,5})$. 
The solid dotted-paths are vacant, and the $*$-paths are occupied. }
\end{figure}

A pivotal vertex  for a positive or a negative event related to four-arm paths around the pivotal vertex is well understood (see Lemma 8 in Kesten (1987)).
The following lemma is a pure topology argument to tell the relation between the right side of (3.5) and the occupied and vacant  paths in
$\bar{\cal R}_{AB}(b, n)$. It is easy to be convinced by  graphs (see  Figs. 1.1--1.6), but it is tedious to show it. So we omit the proof.\\

{\bf Lemma 3.1.}  {\em For a fixed  disjoint $A$ and $B$ in $\{b_1,b_2, b_3, b_4, b_5, b_6\}$
and $b\not \sim {\bf 0}$, }
$$\{N_n({\bf 0})=1
\mbox{ or } 2 \mbox{ if $b$ is occupied or vacant}\,\,|\,\, {\cal E}_{AB}\}=\bar{\cal R}_{AB}(b,n).$$

Let
$$P_{p, AB}(\cdot)= P_p(\cdot\,\, |\,\, {\cal E}_{AB})\mbox{ and } P_{p, 1,3,5}(\cdot)= P_p(\cdot\,\, |\,\, {\cal E}_{1,3,5}).$$
Note that $\bar{\cal R}_{AB}(b, n)$ does not depend on the configurations
of ${\bf 0}$, $A$, $A'$, $B$, and $B'$, so it follows from Lemma 3.1 that
$${dE_p(N_n({\bf 0}) \,\,|\,\, {\cal E}_{AB} )\over dp } =-\sum_{b \in B(n),b\not\sim {\bf 0}} P_{p, AB}(\bar{\cal R}_{AB} (b, n))=-\sum_{b \in B(n),b\not\sim {\bf 0}} P_p(\bar{\cal R}_{AB} (b, n)).\eqno{(3.7)}$$
When we estimate the upper bound of the third derivative of $\kappa(p)$, we only need to deal with  
the following simpler event (see Fig. 1.1 and Fig. 1.5). For fixed vertex $b$ inside $B(n)$, let 
\begin{eqnarray*}
&&{\cal R}_{AB}(b,n)=\\
&&\{\exists\,\,\mbox{ disjoint occupied  paths $r_2$ and $r_4$ in $B(n)$ from 
 $A$ and $B$ to two neighbors of $b$;}\\
&&\exists \mbox{ a vacant path $r_3$ on $B(n)$ inside $S(r_2, r_4, b)$
from $A'$  or $B'$ to a neighbor of $b$; }\\
&&\exists \mbox{ path $r_1$ inside $B(n)$ from  $B'$ or $A'$ to a neighbor}
\mbox{   of $b$ but outside the closure of  $S( r_2, r_4, b)$}\}.
\end{eqnarray*}
For each $ n$, let
$$R_{AB}(p,n) =\sum_{b \in B(n),b\not\sim {\bf 0}} P_p({\cal R}_{AB} (b, n)).$$
Since there is no infinite occupied cluster,  the existence of an occupied path from the origin to $\partial B(n)$ is unlikely to occur.  In other words, we can use ${\cal R}_{AB}(b, n)$
to replace $\bar{\cal R}_{AB}(b, n)$ without losing too much. We may also consider ${\cal R}_{AB}(b, \infty)$ if all the paths in ${\cal R}_{AB}(b, n)$ are in ${\bf Z}^2$.
More precisely, we  show the following lemma.\\

{\bf Lemma 3.2.} 
{\em For each  $A$ and $B$,  ${dE_p(N_n({\bf 0}) \,\,|\,\, {\cal E}_{AB} )\over dp }$ and 
$  {{R}_{AB}(p, n)} $ converge uniformly on $[0, 0.5]$ such that
$$ -\lim_{n\rightarrow \infty} {dE_p(N_n({\bf 0}) \,\,|\,\, {\cal E}_{AB} )\over dp }=\lim_{n\rightarrow \infty} R_{AB}(p, n).\eqno{(3.8)}$$
Furthermore,  if $p_0< 0.5$, for $i\geq 2$
$$  -\lim_{n\rightarrow \infty}{d^iE_p(N_n({\bf 0} \,\,|\,\, {\cal E}_{AB} )\over dp^i }= \lim_{n\rightarrow \infty}{d^{i-1}{R}_{AB}(p, n)\over dp^{i-1}}
\mbox{  converges uniformly on } [0, p_0].\eqno{(3.9)}$$}

{\bf Proof.} 
We  first show that $\sum_{b\in B(n),b\not\sim {\bf 0}} P_p(\bar{{\cal R}}_n(b,n))$  converges uniformly on $[0, 0.5]$. For each $A$ and $B$,
\begin{eqnarray*}
&&\sum_{b\in B(n),b\not\sim {\bf 0}} P_p(\bar{\cal R}_{AB}( b,n))
\leq \sum_{i=1}^{n+1}\sum_{\scriptstyle{b\in B(n),} \atop{\| b\|=i} } 
P_p(\bar{{\cal R}}_{AB}(b,n))\\
&&=\sum_{i=1}^{n/2}\sum_{\scriptstyle{b\in B(n),} \atop{\| b\|=i} }
P_p(\bar{{\cal R}}_{AB}(b,n))+\sum_{i=n/2}^{n+1}\sum_{\scriptstyle{b\in B(n),} \atop{\| b\|=i} } 
P_p(\bar{{\cal R}}_{AB}(b,n)).
\end{eqnarray*}
We denote by $I$ and $II$ the two sums above.
We estimate $I$ first. By the definition of $\bar{\cal R}_{AB}(b,n)$, we know that
if $\{\bar{{\cal R}}_{AB}(b,n)\}$ occurs for $i \leq n/2$, then 
${\cal Q}_{4} ({\bf 0}, 2, i/4)$ and $ {\cal Q}_{4}(b, 0,i/4)$ occur independently (see Figs. 1.5 and 1.6). By (2.6), note that there are 
at most $Ci$ choices for $b$ when $\|b\|=i$, so
$$ I=\sum_{i=1}^{n/2}\sum_{\scriptstyle{b\in B(n),} \atop{\| b\|=i} } 
P_p(\bar{{\cal R}}_{AB}(b,n))\leq C\sum_{i=1}^{n/2} i i^{-5/2+o(1)}.\eqno{(3.10)}$$
Therefore, $I$ converges uniformly on $[0, 0.5]$. Now we estimate $II$.
If $\{\bar{{\cal R}}_{AB}(b,n)\}$ occurs for $n\geq i \geq n/2$, then 
${\cal Q}_{4} ({\bf 0},2, n/4)$ and $ {\cal Q}_{4}(b, 0,(n-i)/4)$ occur independently. By (2.6), 
$$II\leq \sum_{i=n/2}^{n}\sum_{\scriptstyle{b\in B(n),} \atop{\| b\|=i} } 
P_p(\bar{{\cal R}}_{AB}(b,n))\leq \sum_{i=n/2}^{n} n Cn^{-5/4+o(1)} (n-i)^{-5/4+o(1)}.\eqno{(3.11)}$$
Thus, by using the integral test, $II$ goes to zero as $n\rightarrow \infty$ uniformly
on $[0, 0.5]$. 
Therefore,  by (3.7) and the estimates for $I$ and $II$ above, 
$${dE_p(N_n({\bf 0}) \,\,|\,\, {\cal E}_{AB} )\over dp }\mbox{ converges uniformly on $[0, 0.5]$}. \eqno{(3.12)}$$
Note that for $p\leq 0.5$ and a fixed $b$,
$$R_{AB}(p, n) \leq \sum_{b\in B(n),b\not\sim {\bf 0}} P_p(\bar{\cal R}_{AB}(b, n)) \mbox{ and } 
\lim_{n\rightarrow \infty} P_p({\cal R}_{AB}(b, n))= \lim_{n\rightarrow \infty} P_p(\bar{\cal R}_{AB}(b, n)),\eqno{(3.13)}$$
so for $A$ and $B$, by (3.7) and the estimates for $I$ and $II$ above, 
$R_{AB}(p, n)$ converges uniformly and 
(3.8) holds.

It remains to show (3.9).
For a finite $n$,  $R_{AB}(p,n)$ is $i$-th differentiable  for any $i\geq 1$ and $p \leq p_0< 0.5$.
By  the definition of ${\cal R}_{AB}(b, n)$, 
 on ${\cal R}_{AB}(b, n)$, there is an occupied path from 
$A$ to $b$. By Lemma 2.3 and the same proof of (6.108) in Grimmett (1999),
$${d^{i-1} \sum_{b\in B(n)} P_p({\cal R}_{AB}(b, n) )\over dp^{i-1}}\mbox{ converges uniformly on }[0, p_0].\eqno{(3.14)}$$
By (3.14), $R_{AB}^{(i-1)}(p,n)$ converges uniformly on $[0, p_0]$ for $p \leq  p_0< 0.5$ and for $i\geq 2$.
Thus, (3.9)  follows. $\Box$\\

We will introduce a few  estimates for $R_{AB}(p,n)$ as the following proposition.\\

{\bf Proposition 1. } 
{\em For  any $A$ and $B$, if $p\leq 0.5$,  $m=   L(p)$ and $n\geq 2m$,  then}
$$  C L^{-2}(p) (0.5-p)^{-2} \asymp  R_{AB}(p,n)-R_{AB}(p,m).\eqno{(3.15)}$$
{\em In particular,  if $\delta L(p)= m$ and $n \geq 2m$ for a small $\delta >0$, then there exists $C$ independent  of  $\delta$, $p$, $m$, and $n$ such that
$$R_{AB}(p,n) -R_{AB}(p,m)\geq C \delta^{-1/2} L^{-2}(p)(0.5-p)^{-2}.\eqno{(3.16)}$$
Furthermore, if $m =ML(p)$ and $n\geq 2m$ for a large $M\geq 1$, then
there exist $C_1$ and $C_2$  independent  of $p$, $M$, $m$, and $n$  such that}
$$R_{AB} (p,n) -R_{AB}(p,m)\leq C_1\exp(-C_2M) L^{-2}(p)(0.5-p)^{-2}.\eqno{(3.17)}$$
\\

{\bf Proof.}  Note that event ${\cal E}_{AB}$ only depends on the configurations
of $A$, $A'$, $B$, and $B'$, so
$$ P_{p, AB}({\cal Q}_4(i)\,\,\,|\,\,\, {\cal E}_{AB})=P_{p, AB}({\cal Q}_4(i))\asymp P_p({\cal Q}_4(i)).\eqno{(3.18)}$$
By (3.18), and the reconnection lemma,   if $\|b\|=i$, then
$$P_p^2({\cal Q}_4(i/4))\asymp P_{p, AB}^2({\cal Q}_4(i/4))\asymp P_{p, AB}\left({\cal Q}_4( i/4)\cap {\cal Q}_4(b, 0,i/4)\right)\asymp P_p({\cal R}_{AB}(b,i)).\eqno{}$$
Thus, 
$$  P_p({\cal R}_{AB}(b,i))\asymp P_p^2({\cal Q}_4(i/4)).\eqno{(3.19)}$$
Note that $R_{AB}(p, n)$ converges  in $n$ uniformly for $p$, so 
by (3.19), (2.9), and, the reconnection lemma, if $m= L(p)$, then by summing all $b$ with $\|b\|=i \geq L(p)$
and all $ i \leq 2L(p)$,
$$R_{AB}(p,n)-R_{AB}(p,m) \asymp \sum_{i=L(p)}^{2L(p)} iP_p^2({\cal Q}_4(i))  \asymp 
L^2(p) P_p^2( {\cal Q}_4(L(p)))\asymp L^{-2}(p) (0.5-p)^{-2}.\eqno{(3.20)}$$
So (3.15)  in Proposition 1 is proved.

If $m=\delta L(p) $,   
$$ R_{AB}(p,n)-R_{AB}(p,m)=\sum_{m  \leq \|b\|\leq L(p)} {\cal R}_{AB}(b, n)+ \sum_{ L(p) \leq \|b\|} {\cal R}_{AB}(b, n).\eqno{(3.21)}$$
By the same estimate of (3.20) and Lemma 2.4,
$$\sum_{ L(p) \leq \|b\|} {\cal R}_{AB}(b, n)\leq CP^2_p({\cal Q}_4(L(p)))\sum_{i\geq L(p)}  iP_p(B(L(p))\rightarrow B(i))\asymp  L^{-2}(p) (0.5-p)^{-2} .\eqno{(3.22)}$$
If $m=\delta L(p)$, then  the first sum in the right side of (3.21) satisfies that
$$\sum_{m  \leq \|b\|\leq L(p)} {\cal R}_{AB}(b, n) \asymp  \sum_{k=1}^{\delta^{-1} } \sum_{km \leq  i\leq (k+1)m} iP_p^2 ( {\cal Q}_4(i))\asymp \sum_{k=1}^{\delta^{-1} } 
m^2 (k+1)  P_{p}^2 ({\cal Q}_4(k m)).\eqno{(3.23)}$$
By using the reconnection lemma to extend four-arm paths in ${\cal Q}_4(km)$ to the boundary of $B(L(p))$ in (3.23) and by (2.9),
$$\sum_{m  \leq \|b\|\leq L(p)} {\cal R}_{AB}(b, n) \asymp 
 \delta^{-1/2}L^2(p)  P_{p}^2 ({\cal Q}_4(L(p))) \sum_{k=1}^{\infty }  (k+1) k^{-5/2+o(1)} \asymp  \delta^{-1/2} L^{-2}(p) (0.5-p)^{-2}.\eqno {(3.24)}$$
If $\delta$ is small, then  (3.16) in Proposition 1 follows from (3.23) and (3.24).

Now we show (3.17) in Proposition 1. If $\|b\| \geq  M L(p)$, by the same estimate of 
(3.24),
$$ R_{AB}(p,n) -R_{AB}(p,m)\leq 
C P_p^2( {\cal Q}_4(L(p))\sum_{i=ML(p)}^\infty iP_p(\partial B(L(p)))\rightarrow \partial B(i)).\eqno{(3.25)}$$
Applying Lemma 2.4 and (2.9) in (3.25),
$$R_{AB}(p,n)-R_{AB}(p,m)\leq C_1\exp(-C_2M)L^{-2}(p)(0.5-p)^{2}.\eqno{(3.26)}$$
Therefore, (3.17)  follows from (3.26). $\Box$\\

Now we focus on the third term in the right side of (3.4).
Note that
\begin{eqnarray*}
&&{E_p(N_n({\bf 0})\,\,| \,\, {\cal E}_{1,3,5})}=P_p(N_n({\bf 0})\geq 3\,\,| \,\, {\cal E}_{1,3,5})
+P_p(N_n({\bf 0})\geq 2\,\,| \,\, {\cal E}_{1,3,5})+P_p(N_n({\bf 0})\geq 1\,\,| \,\, {\cal E}_{1,3,5}) \\
&&=P_p(N_n({\bf 0})= 3\,\,| \,\, {\cal E}_{1,3,5})
+P_p(N_n({\bf 0})\geq 2\,\,| \,\, {\cal E}_{1,3,5})+1,\hskip 6.7cm 
{(3.27)}
\end{eqnarray*}
so
$${dE_p(N_n({\bf 0})\,\,| \,\, {\cal E}_{1,3,5})\over dp}={dP_p(N_n({\bf 0})=3\,\,| \,\, {\cal E}_{1,3,5})\over dp}
+{dP_p(N_n({\bf 0})\geq 2\,\,| \,\, {\cal E}_{1,3,5})\over dp}.\eqno{(3.28)}$$

For each $b\in B(n)$ with $b\not\sim {\bf 0}$ (on $b_1, b_2, b_3$, which are occupied), let (see Fig. 1.2 and Fig. 1.4)
\begin{eqnarray*}
&&{{\cal R}}_3(b,n)=\\
&&\{\exists  \mbox{ disjoint occupied paths $r_2$ and  $r_4$  in $B(n)$ from $b_1$ and $ b_3$, or $ b_1$ and $ b_5$, or $b_3$ and $b_5$} \\
&&\mbox{ to  neighbors of $b$, } \exists  \mbox{ a vacant path $r_1$ from
$b_2$, or $b_6$, or $b_4$ to a neighbor of $b$, and }\exists\mbox{ a vacant}\\
&&\mbox{ path $r_3$ from $b_4$ and $ b_6$, or $ b_2$ and $b_4$, or $b_2$ and $b_6$ to a neighbor of $b$; or $\exists $ three occupied }  \\
&& \mbox{ and three vacant  paths from $b_1$, $b_2$, $b_3$ and $b_2$, $b_4$,  $b_6$ to neighbors of $b$, all paths in $B(n)\setminus {\bf 0}$}\},
\end{eqnarray*}
and (see Figures 1.3-1.4)
\begin{eqnarray*}
&&{{\cal R}}_2(b, n)=\\
&& \{\exists \mbox{ disjoint vacant paths $r_1$ and  $r_3$ from $b_2$ and $ b_6$, or $ b_2$ and $ b_4$, or $b_4$ and $b_6$} \\
&&\mbox{ to  neighbors of $b$, } \exists  \mbox{ an occupied path $r_2$ from
$b_1$, or $b_3$, or $b_5$ to a neighbor of $b$,  and }\\
&&\mbox{ $\exists$  an occupied path $r_4$ from $b_3$ and $ b_5$, or $ b_1$ and $b_5$, or $b_1$ and $b_3$ to a neighbor of $b$; or}\\
&&  \mbox{ $\exists$ three vacant  paths from $b_1$, $b_2$, $b_3$ and $b_2$, $b_4$,  $b_6$ to neighbors of $b$,
 all paths in $B(n)\setminus {\bf 0}$\}}.  \\
\end{eqnarray*}
We also let
$$R_2(p,n)=\sum_{b\in B(n),b\not\sim {\bf 0}} P_p({\cal R}_2(b, n)) \mbox{ and } R_3(p,n)=\sum_{b\in B(n),b\not\sim {\bf 0}} P_p({\cal R}_3(b, n)).\eqno{(3.29)}$$

Comparing ${\cal R}_{AB}(b, n)$ with ${\cal R}_2(b, n)$ and ${\cal R}_3(b, n)$ (comparing  Fig. 1.1 with Fig. 1.2,  Fig. 1.3, and  Fig. 1.4), there are at least four-arm paths from both ${\bf 0}$ and $b$ among these three events. In fact, if there are six-arm paths from 
both ${\bf 0}$ and $b$ in Fig. 1.4,  the probability estimate is much smaller than
the four-arm case when $p$ is near $0.5$. Thus,  the four-arm case  dominates
the six-arm case.  
With this observation, when we take a higher derivative for $E_p(N_n({\bf 0}))$, we can always do the same computations in the analysis of ${\cal R}_{AB}(b, n)$, ${\cal R}_2(b, n)$, and ${\cal R}_3(b, n)$. Therefore, to avoid repeating the similar proofs many times later, we would rather only deal with the  detailed estimates for the higher derivative for $P_p({\cal R}_{AB}(b, n))$, but omit the 
detailed proofs  for the others.

Using the same arguments in Lemma 3.1, (3.5), (3.7),  and in Lemma 3.2, we can show that
$$\!\!\!\!\!\!\lim_{n\rightarrow \infty} {-dP_p(N_n({\bf 0})=3\,\,| \,\, {\cal E}_{1,3,5})\over dp}= 
\lim_{n\rightarrow \infty} R_3(p, n) \mbox{ and } \lim_{n\rightarrow \infty} {-dP_p(N_n({\bf 0})\geq 2\,\,| \,\, {\cal E}_{1,3,5})\over dp}= 
\lim_{n\rightarrow \infty} R_2(p, n)
\eqno{(3.30)}$$
and  both limits in (3.30) converge uniformly on $[0, 0.5]$.
Thus, by (3.30),
$$\lim_{n\rightarrow \infty} {dE_p(N_n({\bf 0})\,\,| \,\, {\cal E}_{1,3,5})\over dp}=-\lim_{n\rightarrow \infty} [R_2(p,n)+R_3(p,n)],\eqno{(3.31)}$$
and the limit converges uniformly on $[0, 0.5]$.
By (3.2), (3.4),  Lemma 3.2, and (3.31) for all $p\in [0, 0.5]$,
\begin{eqnarray*}
 \kappa''(p)&=& -\lim_{n\rightarrow \infty}{d E_p(N_n)\over dp}\\
&=& - \lim_{n\rightarrow \infty}\sum_{A, B} {d E_p(N_n({\bf 0})\,\,| \,\, {\cal E}_{AB})\over dp} f_{AB}(p)-\lim_{n\rightarrow \infty}\sum_{A,B} E_p(N_n({\bf 0})\,\,| \,\, {\cal E}_{AB}) f_{AB}'(p)\\
 &&-\lim_{n\rightarrow \infty}2{d E_p(N_n({\bf 0})\,\,| \,\, {\cal E}_{1,3,5})\over dp} f_{1,3,5}(p)-  \lim_{n\rightarrow \infty}2E_p(N_n({\bf 0})\,\,| \,\, {\cal E}_{1,3,5}) f_{1,3,5}'(p) -f''(p)\\
 &=&\lim_{n\rightarrow \infty}\left[\sum_{A, B}  f_{AB}(p) R_{AB}(b, n)+ 2f_{1,3,5}(p)\sum_{i=2,3} R_i(p, n)\right]\\
 &&-\lim_{n\rightarrow \infty}\left[\sum_{A,B} E_p(N_n({\bf 0})\,\,|\,\, {\cal E}_{AB})f'_{AB}(p)+ 2E_p(N_n({\bf 0})\,\,|\,\, {\cal E}_{1,3,5})f'_{1,3,5}(p)\right]-f''(p).
 \hskip 2cm (3.32)
 \end{eqnarray*}
 For simplicity, let 
\begin{eqnarray*}
G_p(n)&=&\left[\sum_{A, B}  f_{AB}(p) R_{AB}(p, n)+ 2f_{1,2,3}(p)\sum_{i=2,3} R_i(p, n)\right]\\
 &&-\left[\sum_{A,B} E_p(N_n({\bf 0})\,\,|\,\, {\cal E}_{AB})f'_{AB}(p)+ 2E_p(N_n({\bf 0})\,\,|\,\, {\cal E}_{13,5})f'_{1,3,5}(p)\right]-f''(p).
 \end{eqnarray*}
 Thus, by Lemma 3.2, (3.31), and (3.32), 
 $$ \lim_{n\rightarrow \infty}  G_n(p) =\kappa''(p)\mbox{ converges uniformly on }[0, 0.5].\eqno{(3.33)}$$

Furthermore,  by the same argument of  the estimate in Lemma 3.2,
 if $p \leq p_0< 0.5$, for $i\geq 2$, then
$$  -\lim_{n\rightarrow \infty} {d^iE_p(N_n({\bf 0}) \,\,|\,\, {\cal E}_{1,3,5} )\over dp^i }=\lim_{n\rightarrow \infty} [R_2^{(i-1)} (p,n)+ R_3^{(i-1)}(p,n)]\mbox{ converges uniformly on }[0, p_0].\eqno{(3.34)}$$
By Lemma 3.2 and (3.34),
$$ \kappa^{(i)} (p)= \lim_{n\rightarrow \infty}
G^{(i-2)}_p(n) \mbox{ uniformly on }[0, p_0].\eqno{ (3.35)}$$
By the same argument of Proposition 1, if $m=L(p) $ and $n\geq 2m$,  then
for $j=2,3$
$$CL^{-2}(p) (0.5-p)^{-2} \leq R_j(p,n)-R_j(p, m).\eqno{(3.36)}$$
Similarly to (3.16) and  (3.17) for $R_{AB}(p, n)-R_{AB}(p, m)$, if $m=\delta L(p)$ and $n\geq 2m$, then for $j=2,3$
$$ R_j(p, n)-R_j(p,m)\geq C\delta^{-1/2+o(1)} L^{-2}(p)(0.5-p)^{-2},\eqno{(3.37)}$$
and if $m=ML(p)$,  then for $j=2,3$,
$$R_j(p, n)-R_j(p,m)\leq C\exp(-C_1 M)L^{-2}(p)(0.5-p)^{-2}.\eqno{(3.38)}$$

By the definition (see Fig 1.1--Fig. 1.4), if $n \geq m \asymp L(p)$ on ${\cal E}_{AB}$, then
$N_n({\bf 0}) -N_m({\bf 0})\geq 1 $ implies that there are four-arm paths from $A, A', B, B'$ to $\partial B(m)$. 
Similarly,  on ${\cal E}_{1,3,5}$, there are also four-arm paths
from the origin to $\partial B(m)$. By this observation, (2.3), and (1.14),
\begin{eqnarray*}
&&\left|\sum_{A,B} E_p(N_n({\bf 0})\,\,|\,\, {\cal E}_{AB})f'_{AB}(p) 
-\sum_{A,B} E_p(N_m({\bf 0})\,\,|\,\, {\cal E}_{AB})f'_{AB}(p)\right|\\
&& +2 \left  | E_p(N_n({\bf 0})\,\,|\,\, {\cal E}_{1,3,5})f'_{1,3,5}(p)-
E_p(N_m({\bf 0})\,\,|\,\, {\cal E}_{1,3,5})f'_{1,3,5}(p)\right|\\
&\leq &  C P_p({\cal Q}_4(L(p))) \leq   C(0.5-p)^{5/3+o(1)}.\hskip 9cm (3.39)
\end{eqnarray*}
Note that $(0.5-p)^{5/3+o(1)}$ is much smaller than $L^{-2}(p)(0.5-p)^{-2}$, so 
by using Proposition 1, (3.36), and (3.39)  for $G_n(p)$, if $m=L(p)$ for $p$ near $0.5$ from below, then for all $n\geq 2m$,
there exists a constant $C$ such that
$$ C L^{-2}(p)(0.5-p)^{-2}\leq G_n(p)-G_m(p).\eqno{(3.40)}$$
If $m=\delta L(p)$ and $n \geq 2m$,, by using Proposition 1, (3.37), and (3.39), there exists a constant $C$ independent of $\delta$ such that
$$G_p(n) -G_p(m)\geq C\delta^{-1/2}L^{-2}(p)(0.5-p)^{-2}. \eqno{(3.41)}$$
If $m=ML(p)$ and $n\geq m$,  by Proposition 1, (3.38), and (3.39), there are $C$ and $C_1$ independent of $M$ such that
$$G_p(n) -G_p(m)\leq C\exp(-C_1M)L^{-2}(p)(0.5-p)^{-2}. \eqno{(3.42)}$$

\section{ Higher derivatives of  $\kappa(p)$.}
In this section, we will estimate the third  and higher derivatives for  $\kappa(p)$. 
By (3.32) and  (3.33) for any $p < 0.5$, 
\begin{eqnarray*}
 \kappa'''(p)&=& \lim_{n\rightarrow \infty} G_p'(n)\\
 &=&\lim_{n\rightarrow \infty}\left[\sum_{A, B}  f_{AB}(p) R_{AB}'(p, n)+ 2f_{1,2,3}(p)\sum_{i=2,3}R_i'(p,n)\right]\\
 &+&\lim_{n\rightarrow \infty}\left[2\sum_{A, B}  f_{AB}'(p) R_{AB}(p, n)+ 4f_{1,2,3}'(p)\sum_{i=2,3} R_i(p,n)\right]\\
 &-&\lim_{n\rightarrow \infty}\left[\sum_{A,B} E_p(N_n({\bf 0})\,\,|\,\, {\cal E}_{AB})f''_{AB}(p)+ 2E_p(N_n({\bf 0})\,\,|\,\, {\cal E}_{1,3,5})f''_{1,3,5}(p)\right]-f'''(p).
 \hskip 2cm (4.1)
 \end{eqnarray*}
 It follows from Lemma 3.2 that for all $p \leq 0.5$ that
 \begin{eqnarray*}
 &&\left |2\sum_{A, B}  f_{AB}'(p) R_{AB}(p, n)+ 4f_{1,2,3}'(p)\sum_{i=2,3} R_i(p,n)\right |+\\
 &&\left |\sum_{A,B} E_p(N_n({\bf 0})\,\,|\,\, {\cal E}_{AB})f''_{AB}(p)+ 2E_p(N_n({\bf 0})\,\,|\,\, {\cal E}_{1,3,5})f'_{1,3,5}(p)\right |
  +f'''(p)\mbox{ is uniformly bounded}.
 \hskip .5cm (4.2)
 \end{eqnarray*}
So we only need to take care the first two sums in the right side of (4.1).
 
  We first focus on $R_{AB}'(p, n)$ for  fixed $A$ and $B$. 
Note that ${\cal R}_{AB}(b,n)$
is neither increasing nor decreasing, so we have to introduce
the following more general Russo's formula.  If ${\cal A}$ is increasing and
${\cal B}$ is decreasing, then (see Lemma 1 in Kesten (1987))
\begin{eqnarray*}
{dP_p({\cal A}\cap {\cal B})\over dp}=&&\sum_{u}P_p(u\mbox{ is pivotal for ${\cal A}$ not for
${\cal B}$,  and ${\cal B}$ occurs})\\
&&-\sum_{u}P_p(u\mbox{ is pivotal for ${\cal B}$ not for ${\cal A}$, and  ${\cal A}$ occurs}). 
\hskip 5cm (4.3)
\end{eqnarray*}
We need to divide  ${\cal R}_{AB}(b,n)$ into the intersection of  an increasing  and a decreasing 
events. 
We denote
by ${\cal R}^+_{AB}(b,n)$ the event that there are two disjoint occupied paths $r_2$  and $r_4$
from   the two vertices of $A$ and $B$  to   two neighbor vertices of $b$, respectively. 
We also denote
by ${\cal R}^-_{AB}(b, n)$ the event that there are two disjoint vacant dual paths $r_1$ and $r_3$
from the two vertices of $A'$ and $B'$ to the  two neighbor vertices of $b$, respectively. 
We decompose 
$${\cal R}_{AB}(b,n)={\cal R}_{AB}^+(b,n)\cap {\cal R}_{AB}^-(b,n).\eqno{(4.4)}$$
\begin{figure}
\begin{center}
\setlength{\unitlength}{0.0125in}%
\begin{picture}(250,200)(67,770)
\thicklines
\put(-60,770){\framebox(200,200)[br]{\mbox{$[-n,n]^2$}}}
\put(220,770){\framebox(200,200)[br]{\mbox{$[-n,n]^2$}}}
\put(95,900){\circle*{6}}

\put(85,890){\mbox{$b$}{}}
\put(85,895){*}
\put(80,895){*}
\put(75,895){*}
\put(70,895){*}
\put(68,892){*}
\put(50,882){\mbox{$r_2$}}

\put(65,890){*}
\put(60,890){*}
\put(55,890){*}
\put(50,890){*}
\put(45,885){*}
\put(40,880){*}
\put(35,875){*}
\put(30,870){*}
\put(25,865){*}
\put(20,860){*}
\put(20,855){*}

\put(95,905){\circle*{2}}
\put(95,910){\circle*{2}}
\put(90,910){\circle*{2}}
\put(85,910){\circle*{2}}
\put(80,910){\circle*{2}}
\put(75,910){\circle*{2}}
\put(70,910){\circle*{2}}
\put(65,910){\circle*{2}}
\put(60,910){\circle*{2}}
\put(55,910){\circle*{2}}
\put(50,910){\circle*{2}}
\put(45,910){\circle*{2}}
\put(40,910){\circle*{2}}
\put(35,910){\circle*{2}}
\put(30,910){\circle*{2}}
\put(25,910){\circle*{2}}
\put(20,910){\circle*{2}}
\put(15,910){\circle*{2}}
\put(15,905){\circle*{2}}
\put(15,900){\circle*{2}}
\put(15,895){\circle*{2}}
\put(15,890){\circle*{2}}
\put(15,885){\circle*{2}}
\put(15,880){\circle*{2}}
\put(15,895){\circle*{2}}
\put(15,890){\circle*{2}}
\put(15,885){\circle*{2}}
\put(15,880){\circle*{2}}
\put(15,875){\circle*{2}}
\put(15,870){\circle*{2}}
\put(15,865){\circle*{2}}
\put(15,860){\circle*{2}}
\put(15,855){\circle*{2}}
\put(35,900){\mbox{$r_3$}}

\put(20,855){\circle*{2}}
\put(25,855){\circle*{2}}
\put(30,855){\circle*{2}}
\put(35,855){\circle*{2}}
\put(40,855){\circle*{2}}
\put(45,855){\circle*{2}}
\put(50,855){\circle*{2}}
\put(55,855){\circle*{2}}
\put(60,855){\circle*{2}}
\put(65,855){\circle*{2}}
\put(70,855){\circle*{2}}
\put(75,855){\circle*{2}}
\put(80,855){\circle*{2}}
\put(85,855){\circle*{2}}
\put(90,855){\circle*{2}}
\put(95,855){\circle*{2}}
\put(95,860){\circle*{2}}
\put(95,865){\circle*{2}}
\put(95,870){\circle*{2}}
\put(95,875){\circle*{2}}
\put(95,880){\circle*{2}}
\put(95,885){\circle*{2}}
\put(95,890){\circle*{2}}
\put(95,895){\circle*{2}}
\put(95,900){\circle*{2}}
\put(80,845){\mbox{$r_1$}}
\put(45,845){\mbox{$f$}}
\put(42,855){\circle*{6}}
\put(40,855){*}
\put(40,860){*}
\put(40,865){*}
\put(40,870){*}
\put(40,840){*}
\put(40,835){*}
\put(40,830){*}
\put(40,825){*}
\put(35,825){*}
\put(30,825){*}
\put(25,825){*}
\put(20,825){*}
\put(15,825){*}
\put(10,825){*}
\put(05,825){*}
\put(05,830){*}
\put(05,835){*}
\put(05,840){*}
\put(05,845){*}
\put(05,850){*}
\put(05,855){*}
\put(05,840){*}

\put(98,897){*}

\put(100,900){*}
\put(105,900){*}
\put(110,900){*}
\put(110,905){*}
\put(110,910){*}
\put(110,915){*}
\put(110,920){*}
\put(105,920){*}
\put(100,920){*}
\put(95,920){*}
\put(90,920){*}
\put(85,920){*}
\put(80,920){*}
\put(75,920){*}
\put(70,920){*}
\put(65,920){*}
\put(60,920){*}
\put(55,920){*}
\put(50,920){*}
\put(45,920){*}
\put(40,920){*}
\put(35,920){*}
\put(30,920){*}
\put(25,920){*}
\put(20,920){*}
\put(15,920){*}
\put(10,920){*}
\put(5,920){*}
\put(0,920){*}
\put(-5,920){*}
\put(-10,920){*}
\put(-15,920){*}
\put(-20,920){*}
\put(-20,915){*}
\put(-20,910){*}
\put(-20,905){*}
\put(-20,900){*}
\put(-18,896){*}
\put(-15,892){*}
\put(-12,888){*}
\put(-9,884){*}
\put(-5,880){*}
\put(-5,875){*}
\put(-5,870){*}
\put(-5,865){*}
\put(-5,860){*}
\put(-5,855){*}
\put(0,855){*}
\put(5,855){*}
\put(-20,882){\mbox{$r_4$}}

\put(15,860){\circle*{6}}
\put(12,845){\mbox{${\bf 0}$}{}}
\put(50,870){\mbox{$$}{}}
\put(240,800){\mbox{$$}{}}
\put(270,845){\mbox{$\bf 0$}{}}
\put(280,862){\circle*{6}}
\put(380,900){\circle*{6}}
\put(283,860){*}

\put(283,862){*}
\put(287,865){*}
\put(290,865){*}
\put(295,865){*}
\put(300,865){*}
\put(315,865){*}
\put(320,865){*}
\put(325,865){*}
\put(330,865){*}
\put(335,865){*}
\put(340,870){*}
\put(345,875){*}
\put(335,888){\mbox{$r_2$}}
\put(309,872){\circle*{6}}
\put(315,875){$f$}

\put(310,875){\circle*{2}}
\put(310,880){\circle*{2}}
\put(310,885){\circle*{2}}
\put(310,890){\circle*{2}}
\put(310,895){\circle*{2}}
\put(310,900){\circle*{2}}
\put(310,905){\circle*{2}}
\put(310,910){\circle*{2}}
\put(310,915){\circle*{2}}
\put(310,920){\circle*{2}}

\put(310,870){\circle*{2}}
\put(310,865){\circle*{2}}
\put(310,860){\circle*{2}}
\put(310,855){\circle*{2}}

\put(350,880){*}
\put(355,885){*}
\put(360,890){*}
\put(365,895){*}
\put(370,895){*}
\put(387,895){*}
\put(394,895){*}
\put(400,895){*}
\put(400,890){*}
\put(400,885){*}
\put(400,880){*}
\put(400,875){*}
\put(400,870){*}
\put(400,865){*}
\put(400,860){*}
\put(400,855){*}
\put(400,850){*}
\put(400,845){*}
\put(400,840){*}
\put(400,835){*}
\put(400,830){*}
\put(395,830){*}
\put(390,830){*}
\put(385,830){*}
\put(380,830){*}
\put(375,830){*}
\put(370,830){*}
\put(365,830){*}
\put(360,830){*}
\put(355,830){*}
\put(350,830){*}
\put(345,830){*}
\put(340,830){*}
\put(335,830){*}
\put(330,830){*}
\put(325,830){*}
\put(320,830){*}
\put(315,830){*}
\put(315,820){\mbox{$r_4$}}

\put(310,830){*}
\put(305,830){*}
\put(300,830){*}
\put(295,830){*}
\put(290,830){*}
\put(285,830){*}
\put(280,830){*}
\put(275,830){*}
\put(270,830){*}
\put(265,830){*}
\put(260,830){*}
\put(260,835){*}
\put(260,840){*}
\put(260,845){*}
\put(260,850){*}
\put(260,855){*}
\put(265,855){*}
\put(270,855){*}

\put(280,870){\circle*{2}}
\put(280,875){\circle*{2}}
\put(280,880){\circle*{2}}
\put(280,885){\circle*{2}}
\put(280,890){\circle*{2}}
\put(280,895){\circle*{2}}
\put(280,900){\circle*{2}}
\put(280,905){\circle*{2}}
\put(285,905){\mbox{$r_1$}}

\put(280,910){\circle*{2}}
\put(280,915){\circle*{2}}
\put(285,915){\circle*{2}}
\put(290,915){\circle*{2}}
\put(295,915){\circle*{2}}
\put(300,915){\circle*{2}}
\put(305,920){\circle*{2}}
\put(310,925){\circle*{2}}
\put(315,930){\circle*{2}}
\put(320,935){\circle*{2}}
\put(325,935){\circle*{2}}
\put(330,935){\circle*{2}}
\put(335,935){\circle*{2}}
\put(340,935){\circle*{2}}
\put(345,935){\circle*{2}}
\put(350,935){\circle*{2}}
\put(355,935){\circle*{2}}
\put(360,935){\circle*{2}}
\put(365,935){\circle*{2}}
\put(370,935){\circle*{2}}
\put(375,935){\circle*{2}}
\put(380,935){\circle*{2}}
\put(380,930){\circle*{2}}
\put(380,925){\circle*{2}}
\put(380,920){\circle*{2}}
\put(380,915){\circle*{2}}
\put(380,910){\circle*{2}}
\put(380,905){\circle*{2}}

\put(382,893){\circle*{2}}

\put(385,890){\circle*{2}}
\put(385,885){\circle*{2}}
\put(385,880){\circle*{2}}
\put(385,875){\circle*{2}}
\put(385,870){\circle*{2}}
\put(385,865){\circle*{2}}
\put(385,860){\circle*{2}}
\put(385,855){\circle*{2}}
\put(380,855){\circle*{2}}
\put(375,855){\circle*{2}}
\put(370,855){\circle*{2}}
\put(365,855){\circle*{2}}
\put(360,855){\circle*{2}}
\put(355,855){\circle*{2}}
\put(350,855){\circle*{2}}
\put(345,855){\circle*{2}}
\put(340,855){\circle*{2}}
\put(350,860){\mbox{$r_3$}}

\put(335,855){\circle*{2}}
\put(330,855){\circle*{2}}
\put(325,855){\circle*{2}}
\put(320,855){\circle*{2}}
\put(315,855){\circle*{2}}
\put(310,855){\circle*{2}}
\put(305,855){\circle*{2}}
\put(300,855){\circle*{2}}
\put(295,855){\circle*{2}}
\put(290,855){\circle*{2}}
\put(290,855){\circle*{2}}
\put(285,855){\circle*{2}}
\put(281,855){\circle*{2}}

\put(375,883){\mbox{$b$}{}}
\thicklines
\end{picture}
\end{center}
\caption{\em The left figure shows that $f$ is a pivotal vertex for ${\cal R}_{AB}^-( b,n)$, and the right one shows that $f$ is a pivotal vertex for ${\cal R}_{AB}^+( b,n)$,
 where the $*$-paths are occupied and
the solid circle-paths are vacant.}
\end{figure}
By the new Russo's formula (see Fig. 2),  
\begin{eqnarray*}
R_{AB}'(p,n)
&=&\sum_{b}\sum_{f}
P_p(f\mbox{ is   pivotal for ${\cal R}_{AB}^+(b,n)$ not for
${\cal R}_{AB}^-(b,n)$, and ${\cal R}_{AB}^-(b,n)$ occurs})\\
 &&-\sum_{b}\sum_{f}P_p(f\mbox{  is  pivotal for ${\cal R}_{AB}^-(b,n)$ not for
${\cal R}_{AB}^+(b,n)$, and ${\cal R}_{AB}^+(b,n)$ occurs})\\
&=&\sum_{b}\sum_{f} P_p({\cal R}_{AB}^+(b, f,n))-\sum_{b}\sum_{f}P_p({\cal R}_{AB}^-( b, f,n))={R'}^+_{AB}(p, n)-{R'}^-_{AB}(p, n),\hskip 1cm (4.5)
\end{eqnarray*}
where the first sum above is taken over all sites $b\in B(n)$ for $b\not\sim {\bf 0}$, and the second one is taken over all sites $f$ for $f\neq b$ and $f\not\sim {\bf 0}$ in $B(n)$, and 
\begin{eqnarray*}
&&{\cal R}_{AB}^+( b, f,n)= \{f\mbox{ is   pivotal for ${\cal R}_{AB}^+(b,n)$ not for
${\cal R}_{AB}^-(b,n)$, and ${\cal R}_{AB}^-(b,n)$ occurs}\},\\
&&{\cal R}_{AB}^-(b, f,n)= \{f\mbox{ is   pivotal for ${\cal R}_{AB}^+(b,n)$ not for
${\cal R}_{AB}^-(b,n)$, and ${\cal R}_{AB}^-(b,n)$ occurs}\}.\hskip 1.5cm (4.6)
\end{eqnarray*}
We call them the {\em positive part} and the {\em negative part} of $R_{AB}'(p, n)$, respectively.

Let us focus on ${R'}^+_{AB}(p, n)$. For $p < 0.5$, we sum $b, f$ with 
$\|b\| \leq n/2$ and $\|f\| \leq n/2$. We denote by
$${R'}^+_{AB}(p, n/2, n)=\sum_{b, \|b\|\leq n/2}\sum_{f, \|f\|\leq n/2} P_p({\cal R}_{AB}^+(b, f,n)).$$
We need to point out that ${R'}^+_{AB}(p, n/2, n)$ is not a  derivative of $R^+_{AB}(p, n)$, but  the derivative of $R_{AB}(p, n)$ with restricted pivotal sites   of $b, f\in [-n/2, n/2]^2$.
By the same proof of Lemma 3.2, if $p < 0.5$, 
$$\lim_{n\rightarrow \infty }{R'}^+_{AB}(p, n/2, n)=\lim_{n\rightarrow \infty } \sum_{b, \|b\|\leq n/2}\sum_{f, \|f\|\leq n/2} P_p({\cal R}_{AB}^+(b, f,n))=\lim_{n\rightarrow \infty } \sum_{b}\sum_{f} P_p({\cal R}_{AB}^+(b, f,n)).\eqno{(4.7)}$$
On ${\cal R}_{AB}^+(b,f,n)$ with $\|b\| \leq n/2$ and $\|f\| \leq n/2$,  there are either two vacant  paths $r_5$ and $r_6$ 
from the two neighbors  of $f$ to $r_1$ and $r_3$ (see the right graph of Fig. 2), or there is at least an occupied  path from ${\bf 0}$ or $b$ to $\partial B(n)$. 
Similarly, on ${\cal R}_{AB}^-( b, f,n)$, there are either two occupied paths $r_7$ and $r_8$  from the two neighbors  of $f$ to $r_2$ and $r_4$ (see  the left graph of Fig. 2), or there is at least an occupied path from $f$ to $\partial B(n)$. 
We may also consider 
${\cal R}_{AB}^\pm(b, f,\infty)$ if all the paths above  are 
in ${\bf Z}^2$. By (4.7) and Lemma 3.2,  if $p < 0.5$, then
$$\lim_{n\rightarrow \infty } \sum_{b, \|b\|\leq n/2}\sum_{f, \|f\|\leq n/2} P_p({\cal R}_{AB}^\pm(b, f,n))=\lim_{n\rightarrow \infty } \sum_{b}\sum_{f} P_p({\cal R}_{AB}^\pm(b, f,n))={R_{AB}'}^\pm (p, \infty).\eqno{(4.8)}$$

Now we focus on $b, f\in B(n)$ with $\|b\|\leq n/2$ and $\|f\|\leq n/2$. 
 As we discussed above or in Lemma 3.1 (see Fig. 2),   there are four-arm paths 
 from  ${\bf 0}$ to $\partial B(r)$
if   $b, f\not\in B(r)$ for any $0< r < n$. Similarly, there are four paths around $b$ and $f$, respectively. We call the  {\em pivotal property}.  Note that if we  do not restrict  $b, f\in B(n)$ with $\|b\|\leq n/2$ and $\|f\|\leq n/2$, then there are no four-arm paths when
$b$ or $f$ is near the boundary of $B(n)$. We may still use three-arm paths in the half space to hand it (see Morrow and Zhang (2005)), but it needs more notations and different computations.
We also use the  notations  $dR_i^+(p, n/2, n)/dp$, $dR_{AB}^- (p, n/2, n)/dp$,  and   $dR_i^- (p,n/2,  n)/dp$ for the corresponding sums with $\|b\|, \|f\| \leq n/2$, respectively.
Similarly,  for $i=2,3$,
$$\lim_{n\rightarrow \infty }dR^\pm_{i}  (p, n)/dp=\lim_{n\rightarrow \infty }\sum_{b, f,\|b\|\leq n/2, \|f\|\leq n/2} P_p({\cal R}_{i}^+(b, f, \infty))=dR^\pm_{i}(p, \infty)/dp.\eqno{(4.9)}$$
We also let
\begin{eqnarray*}
d G^+_p(n) /dp
 &=&\sum_{A, B}  f_{AB}(p) dR_{AB}^+(p, n)/dp+ 2f_{1,2,3}(p)\sum_{i=2,3}d R_i^+(p,n)/dp\\
 &+&2\sum_{A, B}  f_{AB}'(p) R_{AB}(p, n)+ 4f_{1,2,3}'(p)\sum_{i=2,3} R_i(p,n)
 + 2E_p(N_n({\bf 0})\,\,|\,\, {\cal E}_{1,3,5})f''_{1,3,5}(p).
 \end{eqnarray*}
 and
 $$d G^-_p(n) /dp=\sum_{A, B}  f_{AB}(p) dR_{AB}^-(p, n)/dp+ 2f_{1,2,3}(p)\sum_{i=2,3}d R_i^-(p,n)/dp+\sum_{A,B} E_p(N_n({\bf 0})\,\,|\,\, {\cal E}_{AB})f''_{AB}(p)+f'''(p).$$
 Thus,
 $$G_p'(n)= d G^+_p(n) /dp-d G^-_p(n) /dp.$$
 Similarly, if we replace $dR_{AB}^\pm(p, n)/dp$ and $dR_i^\pm (p, n)/dp$ by $dR_{AB}^\pm(p, n/2, n)/dp$ and $dR_i^\pm (p,n/2,  n)/dp$ in $dG^\pm_p(n)/dp$, then we have $d G^\pm _p(n/2, n) /dp$.
By (3.33),  if $p< 0.5$, $\lim_{n} dG_p^\pm(n)/dp$ exists. In addition, by (4.7), (4.8),  and (4.9), 
$$\lim_{n\rightarrow \infty} dG_p^\pm(n)/dp=\lim_{n\rightarrow \infty} dG_p^\pm(n/2, n)/dp\stackrel {D}{=}d^3\kappa^\pm(p) /dp^3.\eqno{(4.10)}$$
 With these observations and (4.1), if $p < 0.5$, 
 $$  \kappa'''(p)= d^3\kappa^+(p) /dp^3-d^3\kappa^-(p) /dp^3.\eqno{(4.11)}$$

We may take the second derivative of $G_p(n)$ for $ p < 0.5$ to have
\begin{eqnarray*}
 G_p''(n)&=&\sum_{A, B}  f_{AB}(p) R_{AB}''(p, n)+ 2f_{1,2,3}(p)\sum_{i=2,3}R_i''(p, n)
 +3\sum_{A, B}  f_{AB}'(p) R_{AB}'(p, n)+ 6f_{1,2,3}'(p)\sum_{i=2,3}R'_i(p,n)\\
 &&+3\sum_{A, B}  f_{AB}''(p) R_{AB}(p, n)+ 6f_{1,2,3}''(p)\sum_{i=2,3} R_i(p,n)\\
&& -\sum_{A,B} E_p(N_n({\bf 0})\,\,|\,\, {\cal E}_{AB})f'''_{AB}(p)- 2E_p(N_n({\bf 0})\,\,|\,\, {\cal E}_{1,3,5})f'''_{1,3,5}(p)-f''''(p)\\
&=&\sum_{A, B}  f_{AB}(p) R_{AB}''(p, n)+ 2f_{1,2,3}(p)\sum_{i=2,3}R_i''(p, n)\\
 &&+3\sum_{A, B}  f_{AB}'(p) R_{AB}'(p, n)+ 6f_{1,2,3}'(p)\sum_{i=2,3}R'_i(p,n)+h_4(p, n),
\hskip 4.5cm (4.12)
 \end{eqnarray*}
 where $h_4(p, n)$ is a differentiable function on $(0, 1)$ and convergent uniformly in $n$. 
 By (3.34),  if $p < 0.5$, then
 $$\lim_{n\rightarrow \infty} G''_p(n)=\kappa''''(p) .\eqno{(4.13)}$$
We continue to take the $k$-th derivative.
Before taking the $k$-derivative, we need to simplify the notations.  By Lemma 3.1 and (4.1),  for any $p < 0.5$,
 using the  same discussion  in (4.2),
we only need to estimate $d^kR_{AB}(p, n)/ dp^k$ and $d^k R_i(p, n) / dp^k$ for all possible $A, B$ and for $i=2,3$.
Without  loss of generality,  as we did before,  we only estimate $dR^k_{AB}(p, n)/ dp^k$. 
For fixed $A$ and $B$,  by the definition of (4.5),
$$dR_{AB}(p,n) /dp= \sum_b \sum_{f_1} P_p({\cal R}^+_{AB}(b,f_1, n)) -
\sum_b \sum_{f _1}P_p({\cal R}^-_{AB}(b,f_1,n)).\eqno{(4.14)}$$
We need to make  notations simpler. Let $I_1$ be a sign variable taking $+$ or $-$.
We can rewrite (4.14) as
$$dR_{AB}(p,n) /dp=  \sum_b \sum_{I_1}\sum_{f_1^{I_1}} I_1P_p({\cal R}_{AB}(b,f_1^{I_1}, n)),\eqno{(4.15)}$$
where the above second and third sums take $I_1=+$ or $I_1=-$, and $f_1^{+}$  and $f_1^-$.
 If we continue to take another derivative, then 
$$d^2R_{AB}(p,n) /dp^2= \sum_b \sum_{I_1}\sum_{f_1^{I_1}} I_1{dP_p({\cal R}_{AB}(b,f^{I_1}_1, n))\over dp}.\eqno{(4.16)}$$

Let us work on $dP_p({\cal R}_{AB}(b,f_1^{I_1}, n))/dp$.  Since the topology is similar to Proposition 1, but  needed many   pages to prove, we just state  these topology results in (4.17)  without proofs.  These topology results can also be convinced by Fig. 3.
 We divide ${\cal R}_{AB}(b,f^{I_1}_1, n)$  into the intersection of an increasing event and a decreasing event. More precisely, we collect all the necessary occupied paths connecting ${\bf 0}$, $b$ and $f_1^{I_1}$ as an increasing event and all the necessary vacant paths connecting ${\bf 0}$, $b$ and $f^{I_1}_1$ as a decreasing event,
 where these occupied and vacant paths make ${\cal R}(b,f^{I_1}_1, n)$ to occur (see Fig. 3).
 We then use the new Russo's formula  for positive and negative events separately. 
 Besides of ${\bf 0}$, $b$ and $f^{I_1}_1$, we need to handle a pivotal vertex $f^{I_1,+}_2$ for occupied paths and a pivotal vertex $f^{I_1,-}_2$ for vacant paths. 
 Thus, we can use the sign notation $I_2$ to write
 $$d^2R_{AB}(p,n) /dp^2= \sum_b \sum_{I_1,I_2}\sum_{f^{I_1}, f^{I_1, I_2}} I_1I_2{P_p({\cal R}_{AB}(b,f^{I_1}_1, f^{I_1, I_2}_2, n))}\eqno{(4.17)}$$
to present the probabilities of these four pivotal vertices $f_1^\pm$ and $f_2^{\pm,\pm}$, where $I_1I_2=+$ if $I_1=\pm$ and $I_2=\pm$, otherwise $I_1I_2=-$. 
By using our new notations,  we have the following $k$-th derivative:
$$d^kR_{AB}(p,n)/ dp^{k}=  \sum_{I_1, \cdots, I_k}\sum_{b, f^{I_1}_1, f^{I_1,I_2}_2, \cdots, f^{I_1,I_2\cdots, I_{k}}_{k}} \!\!\!\!\!\!I_1I_2\cdots I_kP_p({\cal R}_{AB}(b,f^{I_1}_1,  f^{I_1,I_2}_2, \cdots, f^{I_1,\cdots, I_{k}}_{k}, n)),\eqno{(4.18)}$$
where the first sum takes all possible signs $I_j=\pm$ for $1\leq j\leq k$,
the second sum takes all  possible  pivotal vertices $\{b, f^{I_1}_1, f^{I_1,I_2}_2, \cdots, f^{I_1,I_2\cdots, I_{k}}_{k}\}$ in event 
${\cal R}_{AB}(b,f^{I_1}_1,  f^{I_1,I_2}_2, \cdots, f^{I_1,\cdots, I_{k}}_{k}, n)$  and 
$I_1I_2\cdots I_{k}=+$ if there are even number of $I_j$'s taking  negative signs, otherwise; it is negative.  
Similarly, we can define $d^kR_i(p,n)/dp^k$ for $i=2,3$.
If  we continue to differentiate  $G_p''(n)$ for $k-2$ times in (4.12),  then by using its definition and (3.33), we have 
\begin{eqnarray*}
&&d^k G_p(n)/dp^k \\
&=&\sum_{A, B}  f_{AB}(p) d^kR_{AB}(p, n)/dp^k+ 2f_{1,2,3}(p)\sum_{i=2,3}d^kR_i(p,n)/dp^k+\\
&&+\left[\sum_{j=1}^{k-1} \sum_{AB}[f_{AB, j}(p) d^jR_{AB}(p, n)/dp^j+\sum_{j=1}^{k-1}g_j(p) \sum_{i=2,3}d^j R_{i} (p, n)/dp^j \right]+ h_k(p, n),\hskip 2cm (4.19)
\end{eqnarray*}
where $f_{AB}(p)$, $f_{1,2,3}(p)$, $f_{AB, j}(p)$ and $g_j(p)$  are some polynomials with  degrees less than 6, and $h_k(p, n)$ is a differentiable function on $(0, 1)$ and convergent uniformly for $n$.
As we did for $k=1$, we denote by
$$d^kR_{AB}(p,n/2, n)/ dp^{k}=  \sum_{I_1, \cdots, I_k}\sum_{b, f^{I_1}_1, f^{I_1,I_2}_2, \cdots, f^{I_1,I_2\cdots, I_{k}}_{k}} \!\!\!\!\!\!I_1I_2\cdots I_kP_p({\cal R}_{AB}(b,f^{I_1}_1,  f^{I_1,I_2}_2, \cdots, f^{I_1,\cdots, I_{k}}_{k}, n))$$
for 
$$\{b,f^{I_1}_1,  f^{I_1,I_2}_2, \cdots, f^{I_1,\cdots, I_{k}}_{k}\}\subset [-n/2, n/2]^2.\eqno{}$$
As mentioned before, the above equation  is not the $k$-th derivative of $R_{AB}(p, n/2, n)$, but  the pivotal sites of $d^kR_{AB}(p, n)/dp^k$ restricted  in $[-n/2, n/2]^2$.
Similarly, we can define $d^kR_i(p,n/2, n)/dp^k$ for $i=2,3$ and $d^k G_p(n/2, n)/dp^k$  by replacing $R_{AB}(p, n)$ and $R_i(p, n)$ by
 $R_{AB}(p,n/2,  n)$ and $R_i(p, n/2, n)$.

Now we try to decompose  ${\cal R}_{AB}(b,f^{I_1}_1,  f^{I_1,I_2}_2, \cdots, f^{I_1,\cdots, I_{k}}_{k}, n)$ into four arm paths in disjoint squares (see Fig. 3). We fix signs $I_1, I_2, \cdots, I_k$.  
For simplicity, we write $f_0=b$ and $f^{I_1,\cdots, I_{j}}_{j}=f_j$.    We need to work on the topology of these pivotal vertices  as the same as we analyzed  in (4.6) by
using  occupied paths and vacant paths from these pivotal vertices  $\{{\bf 0}, f_0, f_1, f_2, \cdots, f_{k}\}$.  
For each pivotal vertex 
$$f\in \{f_0, f_1, f_2, \cdots, f_{k}\}\subset [-n/2, n/2]^2,$$  as we discussed for ${\cal R}_{AB}(b, f, n)$,  if there are no other pivotal vertices  in $f+B(r)$, then ${\cal Q}_4(f, 0, r/2)$ occurs for $0< r < n$. 
 In other words, we still  have the { pivotal property} (see Fig. 3 for $k=3$) for each vertex $f_j$ for $j=0, 1,2, \cdots, k$.  
  \begin{figure}
\begin{center}
\setlength{\unitlength}{0.0125in}%
\begin{picture}(200,400)(67,670)
\thicklines
\put(135, 1010){\framebox(20,20)[br]{}}
\put(135, 990){\framebox(20,20)[br]{}}
\put(155,1000){\mbox{${v_0}$}{}}
\put(145,1000){\circle*{8}}
\put(140,1000){\circle*{2}}
\put(135,1000){\circle*{2}}
\put(130,1000){\circle*{2}}
\put(150,995){\circle*{2}}
\put(155,995){\circle*{2}}
\put(160,995){\circle*{2}}
\put(165,995){\circle*{2}}
\put(170,995){\circle*{2}}
\put(170,1000){\circle*{2}}
\put(170,1005){\circle*{2}}
\put(170,1010){\circle*{2}}
\put(170,1015){\circle*{2}}
\put(170,1020){\circle*{2}}
\put(170,1025){\circle*{2}}
\put(170,1030){\circle*{2}}
\put(170,1035){\circle*{2}}
\put(170,1040){\circle*{2}}
\put(165,1040){\circle*{2}}
\put(160,1040){\circle*{2}}
\put(155,1040){\circle*{2}}

\put(350, 960){\framebox(40,40)[br]{}}
\put(310, 960){\framebox(40,40)[br]{}}
\put(335,980){\mbox{${v_0'}$}{}}
\put(330, 980){\circle*{8}}
\put(320,985){\circle*{2}}
\put(315,990){\circle*{2}}
\put(310,995){\circle*{2}}
\put(335,975){\circle*{2}}
\put(345,965){\circle*{2}}
\put(355,960){\circle*{2}}
\put(340,970){\circle*{2}}
\put(305,1000){\circle*{2}}

\put(-70,810){\framebox(150,150)[br]{\mbox{}}}
\put(-40,840){\framebox(90,90)[br]{\mbox{}}}
\put(0,870){\framebox(30,30)[br]{}}
\put(-25,875){\mbox{$v_1'$}{}}
\put(-13,870){\circle*{8}}
\put(-29, 858){\framebox(28,28)[br]{}}
\put(20,880){\mbox{$v_1$}{}}
\put(12,886){{\circle*{8}}}
\put(80,960){\framebox(150,150)[br]{}}
\put(380,980){\mbox{${ v_0}$}{}}
\put(370,978){{\circle*{8}}}
\put(150,740){\framebox(150,150)[br]{\mbox{}}}
\put(210,800){\mbox{$v_1'$}{}}
\put(217,815){\circle*{8}}
\put(187, 775){\framebox(80,80)[br]{}}

\put(267,795){\framebox(80,80)[br]{}}
\put(305,820){\mbox{$v_1$}{}}
\put(314,835){\circle*{8}}

\put(300,890){\framebox(150,150)[br]{}}
\put(155,1020){\mbox{${v_0'}$}{}}
\put(145,1020){\circle*{8}}

\put(140,1015){\circle*{2}}
\put(135,1010){\circle*{2}}
\put(130,1005){\circle*{2}}
\put(125,1000){\circle*{2}}
\put(120,995){\circle*{2}}
\put(115,990){\circle*{2}}
\put(110,985){\circle*{2}}
\put(105,980){\circle*{2}}
\put(100,975){\circle*{2}}
\put(95,970){\circle*{2}}
\put(90,965){\circle*{2}}
\put(85,960){\circle*{2}}
\put(80,955){\circle*{2}}
\put(75,950){\circle*{2}}
\put(70,945){\circle*{2}}
\put(65,940){\circle*{2}}
\put(60,935){\circle*{2}}
\put(55,930){\circle*{2}}
\put(50,925){\circle*{2}}
\put(45,920){\circle*{2}}
\put(40,915){\circle*{2}}
\put(35,910){\circle*{2}}
\put(30,905){\circle*{2}}
\put(25,900){\circle*{2}}
\put(20,895){\circle*{2}}
\put(20,890){\circle*{2}}
\put(5,885){\circle*{2}}
\put(5,880){\circle*{2}}
\put(0,875){\circle*{2}}
\put(-5,873){\circle*{2}}

\put(140,1020){*}
\put(135,1020){*}
\put(130,1020){*}
\put(125,1020){*}
\put(120,1020){*}
\put(115,1020){*}
\put(110,1018){*}
\put(105,1015){*}
\put(100,1010){*}
\put(95,1005){*}
\put(90,1000){*}
\put(85,995){*}
\put(80,990){*}
\put(75,985){*}
\put(70,980){*}
\put(65,975){*}
\put(60,970){*}
\put(55,965){*}
\put(50,960){*}
\put(45,955){*}
\put(40,950){*}
\put(35,945){*}
\put(30,940){*}
\put(25,935){*}
\put(20,930){*}
\put(15,925){*}
\put(10,920){*}
\put(5,915){*}
\put(0,910){*}
\put(-5,905){*}
\put(-10,900){*}
\put(-15,895){*}
\put(-15,890){*}
\put(-15,885){*}
\put(-15,880){*}
\put(-15,875){*}
\put(-15,870){*}

\put(-10,857){*}
\put(-5,857){*}
\put(0,857){*}
\put(5,857){*}
\put(10,857){*}
\put(15,857){*}
\put(20,857){*}
\put(25,857){*}
\put(30,857){*}
\put(35,857){*}
\put(40,860){*}
\put(45,865){*}
\put(50,870){*}
\put(55,875){*}
\put(60,880){*}
\put(65,885){*}
\put(70,890){*}
\put(75,895){*}
\put(80,900){*}
\put(85,905){*}
\put(90,910){*}
\put(95,915){*}
\put(100,920){*}
\put(105,925){*}
\put(110,930){*}
\put(115,935){*}
\put(120,940){*}
\put(125,945){*}
\put(130,950){*}
\put(135,955){*}
\put(140,960){*}
\put(140,965){*}
\put(140,970){*}
\put(140,975){*}
\put(140,980){*}
\put(140,985){*}
\put(141,990){*}
\put(142,995){*}
\put(143,1000){*}
\put(144,1005){*}

\put(150,1028){\circle*{2}}
\put(150,1034){\circle*{2}}
\put(150,1040){\circle*{2}}
\put(150,1045){\circle*{2}}
\put(145,1045){\circle*{2}}
\put(140,1045){\circle*{2}}
\put(135,1045){\circle*{2}}
\put(130,1045){\circle*{2}}
\put(125,1045){\circle*{2}}
\put(120,1045){\circle*{2}}
\put(115,1045){\circle*{2}}
\put(110,1045){\circle*{2}}
\put(105,1045){\circle*{2}}
\put(100,1045){\circle*{2}}
\put(95,1045){\circle*{2}}
\put(90,1045){\circle*{2}}
\put(85,1045){\circle*{2}}
\put(80,1040){\circle*{2}}
\put(75,1035){\circle*{2}}
\put(70,1030){\circle*{2}}
\put(65,1025){\circle*{2}}
\put(60,1020){\circle*{2}}
\put(55,1015){\circle*{2}}
\put(50,1010){\circle*{2}}
\put(45,1005){\circle*{2}}
\put(40,1000){\circle*{2}}
\put(35,995){\circle*{2}}
\put(30,990){\circle*{2}}
\put(25,985){\circle*{2}}
\put(20,980){\circle*{2}}
\put(15,975){\circle*{2}}
\put(10,970){\circle*{2}}
\put(5,965){\circle*{2}}
\put(0,960){\circle*{2}}
\put(-5,955){\circle*{2}}
\put(-10,950){\circle*{2}}
\put(-15,945){\circle*{2}}
\put(-20,940){\circle*{2}}
\put(-25,935){\circle*{2}}
\put(-30,930){\circle*{2}}
\put(-35,925){\circle*{2}}
\put(-40,920){\circle*{2}}
\put(-40,915){\circle*{2}}
\put(-40,910){\circle*{2}}
\put(-40,905){\circle*{2}}
\put(-40,900){\circle*{2}}
\put(-40,895){\circle*{2}}
\put(-40,890){\circle*{2}}
\put(-40,885){\circle*{2}}
\put(-40,880){\circle*{2}}
\put(-40,875){\circle*{2}}
\put(-40,870){\circle*{2}}
\put(-35,870){\circle*{2}}
\put(-30,870){\circle*{2}}
\put(-25,870){\circle*{2}}
\put(-20,870){\circle*{2}}
\put(-15,870){\circle*{2}}
\put(10,887){*}
\put(10,891){*}
\put(10,892){*}
\put(10,899){*}
\put(10,903){*}
\put(10,907){*}
\put(10,911){*}
\put(10,915){*}

\put(10,872){*}
\put(10,868){*}
\put(10,864){*}
\put(10,860){*}

\put(365,975){\circle*{2}}
\put(365,970){\circle*{2}}
\put(360,965){\circle*{2}}
\put(355,960){\circle*{2}}
\put(350,955){\circle*{2}}
\put(345,950){\circle*{2}}
\put(340,945){\circle*{2}}
\put(335,940){\circle*{2}}
\put(335,935){\circle*{2}}
\put(335,930){\circle*{2}}
\put(335,925){\circle*{2}}
\put(335,920){\circle*{2}}
\put(335,915){\circle*{2}}
\put(335,910){\circle*{2}}
\put(335,905){\circle*{2}}
\put(335,900){\circle*{2}}
\put(335,895){\circle*{2}}
\put(335,890){\circle*{2}}
\put(335,885){\circle*{2}}
\put(335,880){\circle*{2}}
\put(335,875){\circle*{2}}
\put(335,870){\circle*{2}}
\put(335,865){\circle*{2}}
\put(335,860){\circle*{2}}
\put(330,855){\circle*{2}}
\put(325,850){\circle*{2}}
\put(320,845){\circle*{2}}
\put(315,845){\circle*{2}}
\put(305,830){\circle*{2}}
\put(300,825){\circle*{2}}
\put(295,820){\circle*{2}}
\put(290,815){\circle*{2}}
\put(285,815){\circle*{2}}
\put(280,815){\circle*{2}}
\put(275,815){\circle*{2}}
\put(270,815){\circle*{2}}
\put(265,815){\circle*{2}}
\put(260,815){\circle*{2}}
\put(255,815){\circle*{2}}
\put(250,815){\circle*{2}}
\put(245,815){\circle*{2}}
\put(240,815){\circle*{2}}
\put(235,815){\circle*{2}}
\put(230,815){\circle*{2}}
\put(225,815){\circle*{2}}

\put(365,980){*}
\put(360,985){*}
\put(355,990){*}
\put(350,990){*}
\put(345,990){*}
\put(340,990){*}
\put(335,985){*}
\put(330,980){*}
\put(325,975){*}
\put(320,970){*}
\put(315,965){*}
\put(310,960){*}
\put(305,955){*}
\put(300,950){*}
\put(295,945){*}
\put(290,940){*}
\put(285,935){*}
\put(280,930){*}
\put(275,925){*}
\put(270,920){*}
\put(265,915){*}
\put(260,910){*}
\put(255,905){*}
\put(250,900){*}
\put(245,895){*}
\put(240,890){*}
\put(235,885){*}
\put(230,880){*}
\put(225,875){*}
\put(220,870){*}
\put(215,865){*}
\put(215,860){*}
\put(215,855){*}
\put(215,850){*}
\put(215,845){*}
\put(215,840){*}
\put(215,835){*}
\put(215,830){*}
\put(215,825){*}
\put(215,820){*}
\put(215,815){*}
\put(220,802){*}

\put(225,802){*}
\put(230,800){*}
\put(235,795){*}
\put(240,790){*}
\put(245,785){*}
\put(250,780){*}
\put(255,780){*}
\put(260,780){*}
\put(265,780){*}
\put(270,780){*}
\put(275,780){*}
\put(280,780){*}
\put(285,780){*}
\put(290,780){*}
\put(295,780){*}
\put(300,780){*}
\put(305,780){*}
\put(310,780){*}
\put(315,780){*}
\put(320,780){*}
\put(325,780){*}
\put(330,785){*}
\put(335,790){*}
\put(340,795){*}
\put(345,800){*}
\put(350,805){*}
\put(355,810){*}
\put(360,815){*}
\put(365,820){*}
\put(370,825){*}
\put(370,830){*}
\put(370,835){*}
\put(370,840){*}
\put(370,845){*}
\put(370,850){*}
\put(370,855){*}
\put(370,860){*}
\put(370,865){*}
\put(370,870){*}
\put(370,875){*}
\put(370,880){*}
\put(370,885){*}
\put(370,890){*}
\put(370,895){*}
\put(370,900){*}
\put(372,905){*}
\put(375,910){*}
\put(375,915){*}
\put(375,920){*}
\put(375,925){*}
\put(375,930){*}
\put(375,935){*}
\put(375,940){*}
\put(375,945){*}
\put(375,950){*}
\put(375,955){*}
\put(375,960){*}
\put(375,965){*}
\put(370,965){*}
\put(220,815){\circle*{2}}
\put(215,815){\circle*{2}}
\put(210,815){\circle*{2}}
\put(205,815){\circle*{2}}
\put(205,820){\circle*{2}}
\put(205,825){\circle*{2}}
\put(205,830){\circle*{2}}
\put(205,835){\circle*{2}}
\put(205,840){\circle*{2}}
\put(205,845){\circle*{2}}
\put(205,850){\circle*{2}}
\put(205,855){\circle*{2}}
\put(205,860){\circle*{2}}
\put(205,865){\circle*{2}}
\put(205,870){\circle*{2}}
\put(205,875){\circle*{2}}
\put(205,880){\circle*{2}}
\put(205,885){\circle*{2}}
\put(205,890){\circle*{2}}
\put(205,895){\circle*{2}}
\put(205,900){\circle*{2}}
\put(210,905){\circle*{2}}
\put(215,910){\circle*{2}}
\put(220,915){\circle*{2}}
\put(225,920){\circle*{2}}
\put(230,925){\circle*{2}}
\put(235,930){\circle*{2}}
\put(240,935){\circle*{2}}
\put(245,940){\circle*{2}}
\put(250,945){\circle*{2}}
\put(255,950){\circle*{2}}
\put(260,955){\circle*{2}}
\put(265,960){\circle*{2}}
\put(270,965){\circle*{2}}
\put(275,970){\circle*{2}}
\put(280,975){\circle*{2}}
\put(285,980){\circle*{2}}
\put(290,985){\circle*{2}}
\put(295,990){\circle*{2}}
\put(300,995){\circle*{2}}
\put(305,1000){\circle*{2}}
\put(310,1005){\circle*{2}}
\put(315,1010){\circle*{2}}
\put(320,1015){\circle*{2}}
\put(325,1020){\circle*{2}}
\put(330,1020){\circle*{2}}
\put(335,1020){\circle*{2}}
\put(340,1020){\circle*{2}}
\put(345,1020){\circle*{2}}
\put(350,1020){\circle*{2}}
\put(355,1020){\circle*{2}}
\put(360,1020){\circle*{2}}
\put(365,1020){\circle*{2}}
\put(370,1020){\circle*{2}}
\put(375,1020){\circle*{2}}
\put(375,1015){\circle*{2}}
\put(375,1010){\circle*{2}}
\put(375,1005){\circle*{2}}
\put(375,1000){\circle*{2}}
\put(375,995){\circle*{2}}
\put(375,990){\circle*{2}}
\put(375,985){\circle*{2}}
\put(375,980){\circle*{2}}

\put(320,825){*}
\put(320,820){*}
\put(325,815){*}
\put(330,810){*}
\put(335,805){*}
\put(340,800){*}

\put(305,830){*}
\put(300,835){*}
\put(295,840){*}
\put(290,845){*}
\put(285,850){*}
\put(280,855){*}
\put(275,860){*}
\put(270,865){*}
\put(265,870){*}
\put(260,875){*}
\put(255,880){*}
\put(250,885){*}
\put(245,890){*}
\put(230, 720){\vector(0,1){70}}
\put(160, 700){\mbox{${\cal Q}_4(v_1', 0,R_1/4)$ occurs}}
\put(350, 780){\vector(-1,1){30}}
\put(310, 760){\mbox{${\cal Q}_4(v_1, 0, R_1/4)$ occurs}}
\put(350, 1080){\vector(0,-1){80}}
\put(290, 1100){\mbox{${\cal Q}_4(v_0, 0, R_0/4)$ occurs}}
\put(290, 1085){\mbox{${\cal Q}_4(v_0',0,  R_0/4)$ occurs}}

\put(180, 1013){\vector(-1,0){20}}
\put(180, 1020){\mbox{${\cal Q}_4(v_0, 0, R_0/4)$ occurs}}
\put(180, 1005){\mbox{${\cal Q}_4(v_0', 0, R_0/4)$ occurs}}
\put(-40, 1000){\vector(0,-1){60}}
\put(-140, 1040){\mbox{${\cal Q}_4( v_1, R_1/2,  d(v_1, v_0)/4)$ }}
\put(-100, 1020){\mbox{occurs in annulus}}
\put(20, 750){\vector(0,1){130}}
\put(-30, 740){\mbox{${\cal Q}_4( v_1,0, R_1/4)$ occurs}}
\put(-60, 780){\vector(1,2){42}}
\put(-145, 770){\mbox{${\cal Q}_4( v_1', 0, R_1/4)$ occurs}}
\end{picture}
\end{center}
\caption{ \em The left figure shows  case $II_1$(b).  $e_0$  and $e_1$ are the edges adjacent $v_0$ and $v_0'$, and $v_1, $ and $v_1'$, respectively.  In this case, $d(e_0, e_1)> 2R_1=2\|v_1-v_1'\|$.  ${\cal Q}_4(v_0, 0, R_0/4)$,  ${\cal Q}_4(v_0', 0, R_0/4)$, ${\cal Q}_4(v_1, 0, R_1/4)$,  ${\cal Q}_4(v_1', 0, R_1/4), $
${\cal Q}_4(v_1, R_1/2, d(v_1, v_0)/4)$ occur disjointly. The right figure shows case $II_1$(a). In this case, $d(e_0, e_1) \leq  2R_1$.
 ${\cal Q}_4(v_0, 0, R_0/4),$  ${\cal Q}_4(v_0', 0, R_0/4)$, ${\cal Q}_4(v_1, 0, R_1/4)$, $ {\cal Q}_4(v_1', 0, R_1/4)$ occur disjointly. 
 The $*$-paths are occupied and
the solid dotted-paths are vacant.}
\end{figure}

For fixed $\{f_0, \cdots, f_k\}\subset [-n/2, n/2]^2$, we construct a  graph with $k$ edges by
using vertices  $\{f_0, \cdots, f_k\}$.
Let $v_0, v_0'\in \{f_0, \cdots, f_k\}$  with
$$\|v_0-v_0'\|=R_0=\min_{i,j\leq k} \|f_i-f_j\|.$$
 There might be many pairs in $\{f_0, \cdots, f_k\}$ with the same  value $\|v_0-v_0'\|$,  so we simply select one pair in a unique way.  We call   edge $e_0$ with vertices $v_0, v_0'$ and   $B(v_0, R_0/2)$ is called a $v_0$-square. Since  there are no other vertices of  $\{f_0, \cdots, f_k\}$ in $v_0$-square ($v_0'$-square) except $v_0$ ($v_0'$) itself, 
   on ${\cal R}_{AB} (f_0, \cdots, f_k, n)$,  by the pivotal property (see Fig. 3),
 $${\cal Q}_4(v_0, 0, R_0/4)\mbox{ and }   {\cal Q}_4(v_0', 0, R_0/4)\mbox{ occur disjointly on $v_0$- and $v_0'$-squares}.\eqno{(4.20)}$$
  In other words, four arm paths in  $v_0$- and $v_0'$- squares are disjoint, respectively. We denote by
  $${\cal T}_k(v_0, v_0')={\cal Q}_4(v_0, 0, R_0/4)\cap {\cal Q}_4(v_0', 0, R_0/4).$$

We then construct another edge $e_1\neq e_0$ with two vertices $v_1, v_1'$ such that
$$\|v_1-v_1'\|=R_1=\min_{i,j\leq k}\|f_i-f_j\|.$$
 There might be many pairs in $\{f_0, \cdots, f_k\}$ with the same  value $\|v_1-v_1'\|$,  so we simply select one pair in a unique way.   
  Now there are two possible cases. Case $I_1$: $e_0$ and $e_1$ are adjacent (with a common vertex), or  case $II_1$: case $I_1$ does not occur.
 In case $II_1$, on ${\cal R}_{AB} (f_0, \cdots, f_k, n)$, since  there are no other pivotal vertices in $v_i$-square ($v_i'$-square) except $v_i$ ($v_i'$) itself  for $i=0, 1$, by the pivotal property, 
   $${\cal Q}_4(v_0, 0, R_0/4),  {\cal Q}_4(v_0', 0, R_0/4), {\cal Q}_4(v_1, 0, R_1/4),  {\cal Q}_4(v_1', 0, R_1/4)\mbox{ occur disjointly } \eqno{(4.21)}$$ 
 on $v_i$- and $v_i'$-squares for $i=1,2$. In other words, four arm paths in  $v_i$- and $v_i'$- squares, for $i=1,2$, are disjoint, respectively.
We always say that a  group of four arm paths  occur disjointly later if they occur in different squares.
 If $d(e_0, e_1)\leq 2R_1$, called case $II_1(a)$ (see the right graph in Fig. 3), then we just keep the events in (4.21). If  
 $d(e_0, e_1)> 2R_1$, called case  $II_1(b)$,  then  by the pivotal property, 
 $${\cal Q}_4(v_0, 0, R_0/4),  {\cal Q}_4(v_0', 0, R_0/4), {\cal Q}_4(v_1, 0, R_1/4),  {\cal Q}_4(v_1', 0, R_1/4), 
{\cal Q}_4(v_1, R_1/2, d(e_0, e_1)/4) \eqno{ (4.22)}$$
occur disjointly (see the left graph in Fig. 3).  In both case $II_1$(a) and case $II_1$(b), we select $v_1$ such that it is nearer $e_0$ than $v_1'$.

 In case $I_1$,  we know that $e_0$ and $e_1$ are adjacent.  We select $v_1$ adjacent to $e_0$, so
  $v_1\in \{v_0, v_0'\}$ such that
 $\|v_1'-v_1\|=R_1$. Thus,
$${\cal Q}_4(v_0, 0, R_0/4),  {\cal Q}_4(v_0', 0, R_0/4), {\cal Q}_4(v_1', 0, R_1/4),  \mbox{ occur disjointly } \eqno{(4.23)}$$ 
 on $v_0$- and $v_0'$- and $v_1'$-squares, respectively.  
We denote by 
 \begin{eqnarray*}
&&\!\!\!\!\! {\cal T}_k(v_0, v_0', v_1, v_1')=\\
 && \{{\cal Q}_4(v_0, 0, R_0/4),  {\cal Q}_4(v_0', 0, R_0/4), {\cal Q}_4(v_1', 0, R_1/4)\mbox{ in case }I_1\}\\
  &&\bigcup  \{{\cal Q}_4(v_0, 0, R_0/4),  {\cal Q}_4(v_0', 0, R_0/4), {\cal Q}_4(v_1, 0, R_1/4),  {\cal Q}_4(v_1', 0, R_1/4) \mbox{ in case }II_1(a)\}\hskip 2cm (4.24)\\
 &&\!\!\!\!\bigcup \{{\cal Q}_4(v_0, 0, R_0/4),  {\cal Q}_4(v_0', 0, R_0/4), {\cal Q}_4(v_1, 0, R_1/4),  {\cal Q}_4(v_1', 0, R_1/4), 
 {\cal Q}_4(v_1, R_1/2, d(e_0, e_1)/4) \mbox{ in case }II_1(b)\}.
 \end{eqnarray*}

 
 Suppose that we  have edges $\{e_0, \cdots,e_j,\cdots,  e_i\}$ with vertices $\{(v_0, v_0'), \cdots, (v_j, v_j'),\cdots, (v_{i}, v_{i}')\}$ and ${\cal T}_k(v_0, v_0', \cdots, v_i, v_i')$ occurs.
 We continue  to find the edge $e_{i+1}\neq e_j$ for $j=0, \cdots, i$ with vertices $v_{i+1}, v_{i+1}'$ such that 
 $$\|v_{i+1}- v_{i+1}'\|=\min_{i,j\leq k} \|f_i-f_j\|=R_{i+1}.$$
 Similarly to case $I_1$ and case $II_1$ in (4.20)-(4.24), we have case $I_{i+1}$: $e_{i+1}$ is adjacent to $\{e_0, \cdots, e_i\}$ or case $II_{i+1}$:
 case $I_{i+1}$ does not occur. 
 In case $I_{i+1}$, we select $v_{i+1} \in\{v_0, v_0', \cdots, v_i, v_i'\}$. 
 Thus,
$${\cal T}_k(v_0, v_0', \cdots, v_i, v_i'), {\cal Q}_4(v_{i+1}', 0, R_{i+1}/4) \mbox{ occur disjointly}.\eqno{(4.25)}$$
Now we consider case $II_{i+1}$: $e_{i+1}$ is not adjacent to $\{e_0, \cdots, e_i\}$. If
$$\!\!\!\!d(e_{i+1}, \{v_0, v_0', \cdots, v_i, v_i'\}) = d(e_{i+1}, x_i) \leq 2R_{i+1}, \mbox{called case $II_{i+1}$(a), for  $x_i\in \{v_0, v_0', \cdots, v_i, v_i'\}$}, \eqno{(4.26)}$$
note that there are no other pivotal vertices in $v_{i+1} +B(R_{i+1}/4)$  ( $v_{i+1}' +B(R_{i+1}/4)$)  except $v_{i+1}$  ($v_{i+1}'$), so  
$${\cal T}_k(v_0, v_0', \cdots, v_i, v_i'),  {\cal Q}_4(v_{i+1}, 0, R_{i+1}/4), {\cal Q}_4(v_{i+1}', 0, R_{i+1}/4) \mbox{ occur disjointly}.\eqno{(4.27)}$$
 If  
$$\!\!\!\!d(e_{i+1}, \{v_0, v_0', \cdots, v_i, v_i'\}) = d(e_{i+1}, x_i)>  2R_{i+1}, \mbox{called case $II_{i+1}$(b),   for  $x_i\in \{v_0, v_0', \cdots, v_i, v_i'\}$}, \eqno{(4.28)}$$
then  in addition to events in (4.27), there are extra four arm paths in $v_{i+1}+ A(R_{i+1}/2, d(e_{i+1}, x_i)/4)$. Thus,  
$${\cal T}_k(v_0, v_0', \cdots, v_i, v_i'), {\cal Q}_4(v_{i+1}, 0, R_{i+1}/4),  {\cal Q}_4(v_{i+1}', 0, R_{i+1}/4), {\cal Q}_4(v_{i+1}, R_{i+1}/2, d(e_{i+1}, x_i)/4)\eqno{(4.29)}$$
occur disjointly.  In case $II_{i+1}$(a) and  $II_{i+1}$(b), we select $v_{i+1}$ such that it is nearer $\{e_0, \cdots, e_i\}$ than $v_{i+1}'$.
Note that if $\{f_0, \cdots, f_k\}$ are fixed, then $\{v_0,v_0', \cdots, v_{i+1}, v_{i+1}'\}$, and $u$ are also fixed. We denote by
\begin{eqnarray*}
&&{\cal T}_k(v_0, v_0', \cdots, v_{i+1}, v_{i+1}')= {\cal T}_k(v_0, v_0', \cdots, v_i, v_i')\cap\\
&& \{[{\cal Q}_4(v_{i+1}, 0, R_{i+1}/4)\mbox{ in case }I_{i+1}]\cup [{\cal Q}_4(v_{i+1}, 0, R_{i+1}/4),  {\cal Q}_4(v_{i+1}', 0, R_{i+1}/4)\mbox{ in case $II_{i+1}$(a)}]\\
&&\cup [{\cal Q}_4(v_{i+1}, 0, R_{i+1}/4),  {\cal Q}_4(v_{i+1}', 0, R_{i+1}/4), {\cal Q}_4(v_{i+1}, R_{i+1}/2, d(e_{i+1}, x_i)/4)\mbox{ in case $II_{i+1}$(b)}]\}.
\end{eqnarray*}

We continue this way to find all pairs $\{(v_0, v_0')\cdots, (v_{k}, v_{k}')\}$ with edges $\{e_0=(v_0, v_0'), \cdots, e_k=(v_k, v_k')\}$
such that ${\cal T}_k(v_0, v_0', \cdots, v_k, v_k')$ occurs. 
For fixed signs $I_1, I_2, \cdots, I_k$, note that if $\{f_0, \cdots, f_k\}$ are fixed, then $\{v_0, v_0', \cdots, v_k, v_k'\}$ are uniquely determined,  and
$${\cal R}_{AB}(f_0,f_1,  f_2, \cdots, f_{k}, n)\subset {\cal T}_k(v_0, v_0, \cdots, v_k, v_k' ),\eqno{}$$
so
$$\sum_{f_0, \cdots, f_k}P_p({\cal R}_{AB}(f_0,f_1,  f_2, \cdots, f_{k}, n))\leq \sum_{v_0, v_0, \cdots, v_k, v_k'}P_p({\cal T}_k(v_0, v_0, \cdots, v_k, v_k' )).\eqno{(4.30)}$$
We point out that some of  $v_i$ or $ v_i'$ is the origin. We just do not need to sum it in (4.30) if it is the origin. For convenience,  by translation invariance, if we  move in parallel the configurations in ${\cal T}_k(v_0, v_0', \cdots, v_k, v_k' )$ from $v_0$ to the origin, then
$$\sum_{v_0, v_0', \cdots, v_k, v_k'}P_p({\cal T}_k(v_0, v_0', \cdots, v_k, v_k'))\leq \sum_{v_0', \cdots, v_k, v_k'}P_p({\cal T}_k({\bf 0}, v_0', \cdots, v_k, v_k' )),\eqno{(4.31)}$$
where the sum in the right of (4.31)  takes all possible $v_i, v_i'\in B(n)$ with $v_0={\bf 0}$. 

On the other hand, we know that $v_i$ and $v_j$ may or may not be  same for $i\neq j$. To avoid this problem, we select $u_0={\bf 0}$, $u_1=v_0'$, $u_2=v_1$ if $v_1\not\in\{v_0, v_0'\}$  or $u_2=v_1'$ if
$v_1\in \{v_0, v_0'\}$. If $u_i=v_j$ for some $j$,  then $u_{i+1}=v_j'$. Continuously, if $u_i=v_j'$ for some $j$, then 
   $u_{i+1}=v_{j+1}$  when $v_{j+1}\not\in \{{\bf 0}, v_0',\cdots, v_j, v_j'\}$
or $u_{i+1}= v_{j+1}'$ when $v_{j+1}\in \{{\bf 0}, v_0',\cdots, v_j, v_j'\}$.  We continue to find all the different vertices $\{u_0, u_1, \cdots, u_k\}$.  Thus, 
${\cal T}_k({\bf 0}, v_0', \cdots, v_k, v_k' )\subset {\cal T}_k({\bf 0}, u_1, \cdots, u_k )$ with $u_i\neq u_j$ for $i\neq j$ for $\{u_1, \cdots, u_k\}\subset [-2n, 2n]^2$.  On the other hand,   $\{{\bf 0}, v_1', \cdots, v_k, v_k'\}=\{{\bf 0}, u_1, \cdots,  u_k\}$. Thus,
 $$\sum_{v_0', \cdots, v_k, v_k'}P_p({\cal T}_k({\bf 0}, v_0', \cdots, v_k, v_k' ))\leq \sum_{ u_1, \cdots, u_k}P_p({\cal T}_{k}({\bf 0}, u_1, \cdots, u_k)).\eqno{(4.32)}$$
 
 With the vertices $\{u_0, \cdots,  u_k\}$, we have the edges $\{e_0, e_1, \cdots, e_k\}$ connecting these vertices.  These edges  together with $\{u_0, \cdots,  u_k\}$ consist of a few clusters (possible only one). If there are more than one clusters,   we first consider the cluster containing the origin  with  vertices  $\{{\bf 0}, v_0', \cdots, v_i, v_i'\}$ for $i< k$.  We then find the cluster, denoted by $\{v_j, v_j', \cdots, v_l, v_l'\}$,   such that  the distance $T_1$ from its vertices and  $\{{\bf 0}, v_0', \cdots, v_i, v_i'\}$ is the shortest.   If there are two or more different clusters also with the same distance $T_1$ from $\{{\bf 0}, v_0', \cdots, v_i, v_i'\}$, we simply select one uniquely. 
  We pick the two vertices from these two clusters such that the distance between these vertices equal to  $T_1$.   We simply add an extra edge $e_{1}'$  connecting these selected vertices.   The added edge is called {\em extra edge.} It follows from this construction that
$$  R_i \leq  R_j \leq R_l  \leq T_1.\eqno{(4.33)}$$
With the new edge $e_1'$, we get a larger cluster containing the origin.  
By using the same way to expand,  all  the clusters are connected to be one cluster, denoted by ${\cal G}_k$.
For simplicity, , we still use $\{e_j\}$ to denote  all the edges of ${\cal G}_k$.
  For each vertex $u_m$,   we consider the edges of $\{e_j\}$ adjacent to $u_m$ (see Fig. 4).  By our constructions, there is at most one edge from $u_m$ adjacent to  $\{u_1, \cdots, u_{m-1}\}$. On the other hand,  we also consider the edges $\{u_j\}$  adjacent to $u_m$ for $j \geq m$. Suppose that $u_j$ is  one of these edges for $j> m$ adjacent to $u_m$.  Let $u+D(r)$ be the disk with the center at $u$ and with a radius $r$.  
  Let  the intersections of two disks $u_m+D(\|u_m-u_j\|)$ and $u_j+D(\|u_m-u_j\|)$ be  $\alpha_1$ and $\alpha_2$ (see Fig. 4). We consider the rays  from $u_m$ passing through  $\alpha_i$  to $\infty$ for $i=1,2$ and the cone $\Delta_m$ with boundaries of these two rays (see Fig. 4). 
  By our constructions for $\{u_0, \cdots, u_k\}$,  if $(u_m, u_j)$ is not an extra edge, then $\Delta_m$   contains only  edge $(u_m, u_j)$
  in its domain.   We  call the {\em isolated} property for $u_m$. If $(u_m, u_j)$ is an extra edge, by (4.33), we still have the isolated property for $u_m$. Thus, we may say that ${\cal G}_j$ has the isolated property for each of its vertices.  We consider the  8 rays  including the horizontal -axis starting  from $u_m$  to $\infty$ such that  the angle  between any two rays is $\pi/4$ (see Fig. 4). Thus, there are 8 equal cones with the boundaries of these rays. We order them by corn 1, $\cdots,$ corn 8. By the isolated property,  there are at most  one edge $(u_m, u_j)$ of ${\cal G}_k$ adjacent to $u_m$ in each corn. We call this {\em single edge} property. Since they are just  elementary geometry arguments, we do not precisely prove these isolated and  single edge properties, but give  explanations  in Fig. 4.  Thus, there are at most 8 edges adjacent to $u_m$ for each vertex $u_m$ in ${\cal G}_k$.
   On the other hand, ${\cal G}_k$ does not have a loop path by the construction of ${\cal G}_k$.
 We summarize these results as the following statement.
 $$\mbox{ ${\cal G}_k$ is connected  without a loop path and the degree of each of its vertices is less than 8}. \eqno{(4.34)}$$
\begin{figure}
\begin{center}
\setlength{\unitlength}{0.0125in}%
\begin{picture}(250,190)(67,860)
\thicklines

\put(50, 1000){\line (1,0){ 200}}
\put(150, 900){\line (0,1){ 200}}
\put(50, 900){\line (1,1){ 200}}
\put(250, 900){\line (-1,1){ 200}}
\put(150, 1000){\circle{ 100}}
\put(150, 1000){\circle{ 100}}
\put(162, 1020){\circle{ 100}}
\put(150, 1000){\line (1,2){ 10}}
\put(151, 1001){\line (1,2){ 9}}
\put(150, 1000){\vector (1,0){ 100}}
\put(150, 1000){\vector (-1,2){ 40}}
\put(150, 1000){\circle*{ 8}}
\put(160, 1020){\circle*{ 8}}
\put(146, 990){\mbox{$u_m$}}
\put(160, 1025){\mbox{$u_j$}}
\put(170, 1000){\circle*{ 4}}
\put(140, 1020){\circle*{ 4}}
\put(175, 990){\mbox{$\alpha_1$}}
\put(130, 1025){\mbox{$\alpha_2$}}

\end{picture}
\end{center}
\caption{ \em Edge $(u_m, u_j)\in \{e_j\}$.  Two disks $u_m+D(\|u_m-u_j\|)$ and $u_j+D(\|u_m-u_j\|)$ intersect at   $\alpha_1$ and $\alpha_2$, respectively. Two rays  from
  $u_m$ passing through  $\alpha_i$  to $\infty$ for $i=1,2$ consisting of cone $\Delta_m$. The angle between two boundaries of the cone is $2\pi/3$.  By the constructions of the pivotal sites, there is no edge adjacent to $u_m$  in $\Delta_m \cap (u_m+D(\|u_m-u_j\|) )$ except 
  $(u_m, u_j)$.
  On the other hand, note that any vertex  in $\Delta_m \setminus (u_m+D(\|u_m-u_j\|) )$ is nearer $u_j$ than $u_m$, so there are no 
  vertices in $\Delta_m \setminus (u_m+D(\|u_m-u_j\|) )$ adjacent to $u_m$ as the edges of $\{e_j\}$. In other words, 
  $(u_m, u_j)$ is the only edge adjacent to $u_m$ in the cone.  There are  
  8 equal cones starting from $u_m$.  Each of them contains at most only one edge in $\{e_j\}$ adjacent to $u_m$.}
\end{figure}

We consider any two  clusters  ${G}_k$,  and  $\bar{ G}_k$ with edges $\{e_j\}$ and vertices $\{u_0={\bf 0}, u_1,\cdots, u_k\}\subset {\bf Z}^2$,
and edges $\{\bar{e}_j\}$ and vertices $\{\bar{u}_0={\bf 0}, \bar{u}_1,\cdots, \bar{u}_k\}\subset {\bf Z}^2$, respectively. We want to point out that
the above vertices are in ${\bf Z}^2$ ordered from $u_0, \cdots, u_k$ for both groups.
We say that the  two clusters are {\em similar}, denoted by  ${G}_k\sim \bar{G}_k$, if 
  $u_i$ is only adjacent to $u_{j_1}, \cdots, u_{j_l}$,  and $\bar{u}_i$ is also  only adjacent to   $\bar{u}_{j_1}, \cdots, \bar{u} _{j_l}$ for each $0\leq i\leq k.$  
In other words, if ${ G}_k\sim \bar{G}_k$, then the situations of the connections of the two clusters are the same, 
but  the locations of  $u_i$  and $\bar{u}_i$   in the two clusters may not be  the same.
From now we only consider a cluster with $k$ vertices on ${\bf Z}^2$ and with the single edge property. We can decompose all  such the clusters into similar groups $\langle G_k\rangle$, where each $G_k$ is one of the  selected cluster in a group.
We would like to account how many such groups.
We consider a regular tree ${\bf T}$ with a degree $8$ and a root at the origin.
Let ${ A}_k$ be an {\em animal} defined  to be a finite connected subgraph of ${ \bf T}$ containing the root with $k$ vertices. Let $a_k$  be the number of such animals of $A_k$.
 By a standard graph estimate (see (4.24) in Grimmett (1999)),   
 $$a_k  \leq 7^{8k}.\eqno{( 4.35)}$$
 For each cluster of $G_k$, it is similar to an animal $A_k$ in ${\bf T}$ and two non similar clusters $G_k'$ and $G_k''$ are also corresponding two non similar animals $A_k'$ and $A_k''$ in ${\bf T}$, respectively.  Thus, we can find at least one-many-map from $G_k$ to $A_k$.  By (4.35), 
 $$\mbox{ the number of similar cluster groups  in $ \langle G_k\rangle $ is less than } 7^{8k}.\eqno{ (4.36)}$$
   With this decomposition,   
      $$\sum_{ u_1, \cdots, u_k}P_p({\cal T}_{k}({\bf 0}, u_1, \cdots, u_k))=\sum_{G_k} \sum_{ u_1, \cdots, u_k}P_p({\cal T}_{k}({\bf 0}, u_1, \cdots, u_k), {\cal G}_k\sim G_k),\eqno{}$$
  where the first sum in the right side  takes over all possible similar cluster groups $\langle G_k\rangle$.  
   With these definitions, we show the following lemma.\\

{\bf Lemma 4.1.} {\em If $p < 0.5$, and   $k\geq 2$,  then for each  selected $G_k\in \langle G_k\rangle$ and for all large $n>k$,
there exists $C$ independent of  $k$, $n$ and $p$ such that
$$\sum_{ u_1, \cdots, u_k}P_p({\cal T}_{k}({\bf 0}, u_1, \cdots, u_k),{\cal G}_k\sim G_k)\leq C(0.5-p)^{-1}\sum_{u_1, \cdots, u_{k-1}}P_p({\cal T}_{k-1}({\bf 0}, u_1, \cdots, u_{k-1}), {\cal G}_{k-1}\sim G_{k-1}),$$
 where $\{u_1, \cdots, u_k\}\subset [-2n, 2n]^2$. 
}\\

   {\bf Proof.}  
On ${\cal G}_k\sim G_k$ for a selected and fixed $G_k$.
Recall that we divided into two cases $I_k$,   $II_k$(a) and $II_k$(b). Thus, 
\begin{eqnarray*}
&&\sum_{ u_1, \cdots, u_k}P_p({\cal T}_{k}({\bf 0}, u_1, \cdots, u_k), {\cal G}_{k}\sim G_{k})
=\sum_{ u_1, \cdots, u_k}P_p({\cal T}_{k}({\bf 0}, u_1, \cdots, u_k), \mbox{case }I_k, {\cal G}_{k}\sim G_{k})\\
&&+ \sum_{ u_1, \cdots, u_k}P_p({\cal T}_{k}({\bf 0}, u_1, \cdots, u_k), \mbox{case $II_k$(a)}, {\cal G}_{k}\sim G_{k})\\
&&+ \sum_{ u_1, \cdots, u_k}P_p({\cal T}_{k}({\bf 0}, u_1, \cdots, u_k), \mbox{case $II_k$(b)},{\cal G}_{k}\sim G_k),\hskip 7cm (4.37)
\end{eqnarray*}
where cases $I_k$, $II_k$(a) and $II_k$(b) are  the connecting situations of edge $e_k$ from $u_k$ to the other pivotal sites.
We first focus on the first sum in the right side of (4.37).
In case $I_k$,  $e_k$ is adjacent to $\{e_0, \cdots, e_{k-1}\}$. Thus,  $u_{k-1}=v_{k}\in \{{\bf 0}, v_0', \cdots, v_{k-1}, v_{k-1}'\}$ and $u_k=v_{k}'$.
If $u_1, \cdots, u_k$ are fixed, then  these disjoint squares constructed in event ${\cal T}_{k}({\bf 0}, u_1, \cdots, u_k)$ are fixed and disjoint. Thus, these four arm paths occur independently on these disjoint squares.
By this observation and using (4.25) for $i+1=k$, we have
$${\cal T}_k({\bf 0}, u_1, \cdots, u_{k-1}), {\cal Q}_4(u_k, 0, R_{k}/4) \mbox{ occur indepdently}.\eqno{(4.38)}$$
Since ${\cal G}_k\sim G_k$ for a  fixed $G_k$, then $e_k$ is the edge from $u_k$ to $u_i\in \{u_0, \cdots, u_{k-1}\}$ for a fixed $i$.
To sum $\{u_1, \cdots, u_k\}$, we may first sum  all possible $\{u_k\}$ for fixing $\{u_1, \cdots, u_{k-1}\}$, and then sum all possible
$\{u_1, \cdots, u_{k-1}\}$.
We fix $\{u_1, \cdots, u_{k-1}\}$ and  sum all possible $u_k$. Since the above $i$ is fixed, $u_i$ is fixed.
So we sum all $u_k$ from $u_i$ with $\|u_i-u_k\|=l$ for a fixed $u_i$ and sum all $l$. Moreover, if ${\cal G}_k\sim G_k$, after removing  $u_k$ and $e_k$, ${\cal G}_k$ and $G_k$ will be ${\cal G}_{k-1}$ and $G_{k-1}$, respectively and ${\cal G}_{k-1}\sim G_{k-1}$. With these observations together with   (4.38) and Lemma 2.7, there exists $C$ such that 
\begin{eqnarray*}
&&\sum_{ u_1, \cdots, u_k}P_p({\cal T}_{k}({\bf 0}, u_1, \cdots, u_k),{\cal G}_k\sim G_k, \mbox{case }I_k)=\sum_{ u_1, \cdots, u_{k-1}}\sum_{u_k}P_p({\cal T}_{k}({\bf 0}, u_1, \cdots, u_k),{\cal G}_k\sim G_k, \mbox{case }I_k)\\
&&\leq \sum_{ u_1, \cdots, u_{k-1}}P_p({\cal T}_{k}({\bf 0}, u_1, \cdots, u_{k-1}),{\cal G}_{k-1}\sim G_{k-1})\sum_{l=1}^{\infty} l P_p({\cal Q}_4(l/4))\\
&\leq & C(0.5-p)^{-1} \sum_{ u_1, \cdots, u_{k-1}}P_p({\cal T}_{k}({\bf 0}, u_1, \cdots, u_{k-1}),{\cal G}_{k-1}\sim G_{k-1}),\hskip 5.5cm {(4.39)}
\end{eqnarray*}
where $G_{k-1}$ is a cluster with vertices $\{u_0, \cdots, u_{k-1}\}$ and with $\{e_j\}\setminus e_k$.
For each subset $\{e_0, \cdots, e_{k-1}\}$ of $\{e_0, \cdots, e_{k-1},  e_{k}\}$, it is just a subset of $\{e_0, \cdots, e_{k-1}\}$ itself.
Thus,  
\begin{eqnarray*}
&&\sum_{ u_1, \cdots, u_k}\!\!\!P_p({\cal T}_{k}({\bf 0}, u_1, \cdots, u_k), \mbox{case }I_k,{\cal G}_{k}\sim G_{k})\\
&\leq & C(0.5-p)^{-1} \!\!\!\sum_{ u_1, \cdots, u_{k-1}}\!\!\!\!P_p({\cal T}_{k-1}({\bf 0}, u_1, \cdots, u_{k-1}),{\cal G}_{k-1}\sim G_{k-1}).\hskip 5.3cm {(4.40)}
\end{eqnarray*}

In case $II_k$(a),  $e_k$ is not adjacent to $\{e_0, \cdots, e_{k-1}\}$, but there is  an extra edge $e'$ connecting $u_i$ and $e_k$ for a fixed  
$u_i\in \{u_0, \cdots, u_{k-2}\}$, since $G_k$ is fixed.   Thus, $u_k=v_k'$,  $u_{k-1}=v_{k}$,  and $v_k\not\in  \{{\bf 0}, v_0', \cdots, v_{k-1}, v_{k-1}'\}$. 
 By (4.27) for $i+1=k$ and the same argument in (4.38),  for fixed $\{u_1, \cdots, u_k\}$,
$${\cal T}_k({\bf 0}, u_1, \cdots, u_{k-2}),  {\cal Q}_4(u_{k-1}, 0, R_{k}/4), \mbox{ and }  {\cal Q}_4(u_{k}, 0, R_{k}/4) \mbox{ occur independently}.\eqno{}$$
By the assumption of case II (a),
$${\cal T}_k({\bf 0}, u_1, \cdots, u_{k-2})\cap  {\cal Q}_4(u_{k-1}, 0, R_{k}/4)\subset {\cal T}_k({\bf 0}, u_1, \cdots, u_{k-1}).\eqno{(4.41)}$$
With these observations, if we use the same proofs of  (4.39) and (4.40) to sum $u_k$, then
\begin{eqnarray*}
&&\sum_{ u_1, \cdots, u_k}P_p({\cal T}_{k}({\bf 0}, u_1, \cdots, u_k), \mbox{case $II_k$(a)},{\cal G}_{k}\sim G_{k})\\
&\leq &\sum_{ u_1, \cdots, u_{k-1} }P_p({\cal T}_{k}({\bf 0}, u_1, \cdots, u_{k-2}), {\cal Q}_4(u_{k-1}, 0, R_{k}/4), {\cal G}_{k-1}\sim G_{k-1}) \sum_{l=1}^\infty l P_p({\cal Q}_4(l/4))\\
&\leq & \sum_{u_1, \cdots, u_{k-1}}P_p({\cal T}_{k}({\bf 0}, u_1, \cdots, u_{k-1}), {\cal G}_{k-1}\sim G_{k-1}) \sum_{l=1}^\infty l P_p({\cal Q}_4(l/4))\\
&\leq & C(0.5-p)^{-1}\sum_{u_1, \cdots, u_{k-1}}P_p({\cal T}_{k-1}({\bf 0}, u_1, \cdots, u_{k-1}),{\cal G}_{k-1}\sim G_{k-1}). \hskip 5cm (4.42)
\end{eqnarray*}

Finally,  we work on case $II_k$(b).    We find the same $u_i$  and $e'$ in case $II_k$(a) for a fixed $u_i\in \{u_0, \cdots, u_{k-1}\}$. 
By (4.29) for $i+1=k$,
$$\!\!\!\!\!{\cal T}({\bf 0}, u_1, \cdots, u_{k-1}), {\cal Q}_4(u_{k-1}, 0, R_{k}/4),  {\cal Q}_4(u_k, 0, R_{k}/4), {\cal Q}_4(u_{k-1}, R_{k}/2, d(u_i, e_k)/4)\mbox{ occur disjointly.}\eqno{(4.43)}$$
Thus,  by using the disjoint property in (4.43), the  reconnection lemma,  and the same estimate in (4.42), there exists $C_1$ such that
\begin{eqnarray*}
&&\sum_{ u_1, \cdots, u_{k}}P_p({\cal T}_{k}({\bf 0}, u_1, \cdots, u_k), {\cal G}_{k}\sim G_{k}, \mbox{case $II_k$(b)})\\
&\leq &  \sum_{u_1, \cdots, u_{k-1}}P_p({\cal T}_{k}({\bf 0}, u_1, \cdots, u_{k-2}), {\cal Q}_4(u_{k-1}, 0, R_{k}/4), {\cal Q}_4(u_{i}, R_{k}/2, d(e_{k}, u_i)/4),{\cal G}_{k-1}\sim G_{k-1})\\
&&\hskip 1cm \times  \sum_{l=1} l P_p({\cal Q}_4(l/4))\\
&\leq & C (0.5-p)^{-1}  \sum_{u_1, \cdots, u_{k-1}}P_p({\cal T}_{k}({\bf 0}, u_1, \cdots, u_{k-2}), {\cal Q}_4(u_{k-1}, 0, d(e_{k}, u_i)/4), {\cal G}_{k-1}\sim G_{k-1}) \\
&\leq & C_1(0.5-p)^{-1}\sum_{u_1, \cdots, u_{k-1}}P_p({\cal T}_{k-1}({\bf 0}, u_1, \cdots, u_{k-1}), {\cal G}_{k-1}\sim G_{k-1}). \hskip 5cm (4.44)
\end{eqnarray*}
Thus,  Lemma 4.1 follows from (4.40), (4.42) and (4.44). 
$\Box$\\

With Lemma 4.1,  we show the following proposition for higher derivatives of $\kappa(p)$.\\

{\bf Proposition 2.} {\em If  $p < 0.5$,  then for any integers $n$ and  $n >  k\geq 1$, there exists $C$  independent of $k$, $n$ and $p$ such that
$$ | d^kG_p(n/2,n) /dp^k |  \leq C^k L^{-2}(p) (0.5-p)^{-k-2}$$
and
$$ |d^{k+2}\kappa(p)/ dp^{k+2}|\leq  C^k   L^{-2}(p) (0.5-p)^{-k-2}.$$}

{\bf Remark} 4. In the proof of Proposition 2,  $k$ is finite  independent of $p$ and $n$, but large in order to control the constants in Propositions 1-2.  
By using the same estimate of Morrow and Zhang  (2005), one can show that 
$$ |d^{k+2}\kappa(p)/ dp^{k+2}|\leq   C(k) L^{-2}(p) (0.5-p)^{-k-2}$$
for some constant $C(k)$ depending on $k$, but this is not good enough to show the Theorem. We need a better upper bound $C^k L^{-2}(p) (0.5-p)^{-k-2}$. \\

{\bf Proof of Proposition 2.} 
By using Lemma 4.1,  for each selected $G_k\in\langle G_k\rangle$, we iterate 
$$\sum_{ u_1, \cdots, u_k}P_p({\cal T}_{k}({\bf 0}, u_1, \cdots, u_k), {\cal G}_k\sim G_k)$$
 in $(k-1)$
 times  for fixed signs $I_1, \cdots, I_k$ and for $\{u_1, \cdots, u_k\}\subset [-2n, 2n]^2$  to show
 $$\sum_{u_1, \cdots, u_k} P_p({\cal T}_{k}({\bf 0}, u_1, \cdots, u_k), {\cal G}_k\sim G_k)\leq C^{k-1} (0.5-p)^{-k+1} \sum_{u_1} P_p({\cal T}_{2}({\bf 0}, u_1)).\eqno{(4.45)}$$
 By using Proposition 1 in the right side of (4.45),  if $\{u_1, \cdots, u_k\}\subset [-2n, 2n]^2$, then
 $$\sum_{u_1, \cdots, u_k} P_p({\cal T}_{k}({\bf 0},u_1, \cdots, u_k),{\cal G}_k\sim G_k)\leq C^{k}_1 (0.5-p)^{-k-2}L^{-2}(p).\eqno{(4.46)}$$
 By (4.35), (4.36), and (4.46), 
 \begin{eqnarray*}
&&\sum_{ u_1, \cdots, u_k}P_p({\cal T}_{k}({\bf 0}, u_1, \cdots, u_k))\\
&=&\sum_{G_k}\sum_{ u_1, \cdots, u_k}P_p({\cal T}_{k}({\bf 0}, u_1, \cdots, u_k), {\cal G}_{k}\sim G_{k})\leq  7^{8k} C_1^{k} (0.5-p)^{-k-2}L^{-2}(p)
.\hskip 3.5cm (4.47)
\end{eqnarray*}
 By (4.30)-(4.31) and (4.47), note that each sign $I_j$ only has two choices, so  there exists $C_2$ independent of $k$, $p$ and $n$ such that
 $$|d^kR_{AB}(p,n/2, n)/ dp^{k}|\leq   2^k\sum_{b, f_1, f_2, \cdots, f_{k}}P_p({\cal R}_{AB}(b,f_1,  f_2, \cdots, f_{k}, n))\leq C^k_2(0.5-p)^{-k-2}L^{-2}(p).\eqno{(4.48)}$$
Similarly, we can work on $|d^kR_{i}(p,n)/ dp^{k}|$ for $i=2, 3$ to have the same upper bound in (4.48).

Thus, by using (4.48) and the same upper bound of (4.48) for  $|d^kR_{i}(p,n)/ dp^{k}|$ in (4.19) but with $R_{AB}(p, n/2, n)$ and $R_i(p, n/2, n)$ instead,   there exists $C$ independent of $p$, $k$ and $n$ such that
$$ | d^kG_p(n/2,n) /dp^k |  \leq  C^k L^{-2}(p) (0.5-p)^{-k-2}.\eqno{(4.49)}$$
If we let $n\rightarrow \infty$ in (4.49), by (3.35), 
$$ |d^{k+2}\kappa(p)/ dp^{k+2}|=\lim_{n\rightarrow \infty} | d^kG_p(n) /dp^k |=\lim_{n\rightarrow \infty} | d^kG_p(n/2, n) /dp^k |\leq   C^kL^{-2}(p) (0.5-p)^{-k-2}.\eqno{(4.50)}$$
 Proposition 2  follows from (4.49) and (4.50). $\Box$\\

Now we show the following proposition for an estimate of a convergent rate of $d^kG_p(n) /dp^k$ to $ d^{k+2}\kappa(p)/ dp^{k+2}$.\\

{\bf Proposition 3.} {\em If  $p < 0.5$  and $n\geq 2m L(p)$ for a large $m\geq k^3$ and for a large $k$, then there exists $C$ independent of $k$, $m$ and $p$ such that
$$  |d^kG_p(n)/dp^k-d^{k+2} \kappa(p)/dp^{k+2}|\leq  e^{-Cm/k} L^{-2}(p)(0.5-p)^{-k-2}.$$}\\

{\bf Remark} 5. The estimate in Proposition 3 is not optimal, but it is good enough to show the Theorem.
We also want to point out that $m$ is large, but independent  of $p$.
In other words, the upper bound in the inequality of Proposition 3 is large  for $p$ near 0.5.\\

{\bf Proof of Proposition 3.}  We work on ${\cal T}_{k}({\bf 0}, u_1, \cdots, u_k)$ for 
$\{u_1, \cdots, u_k \}\subset [-n, n]^2.$
Let 
$Y_{j}=R_1+\cdots +R_{j}$ for $R_j=\|v_j-v_j'\|$ defined before. Note that ${\cal G}_k$ is a connected graph as we discussed  before, so  if $\{u_1, \cdots, u_k\} \not \subset [-2m L(p), 2m L(p)]^2$, then 
$$Y_k \geq m L(p).$$
For a fixed $G_k$, we decompose  ${\cal T}_{k}({\bf 0}, u_1, \cdots, u_k)$  by using $Y_j$ to have
$$\sum_{u_1, \cdots, u_k} P_p({\cal T}_{k}({\bf 0}, u_1, \cdots, u_k), {\cal G}_k\sim G_k, Y_k \geq mL(p))
=\sum_{m L(p)\leq l}\sum_{u_1, \cdots, u_k} P_p({\cal T}_{k}({\bf 0}, u_1, \cdots, u_k),  {\cal G}_k\sim G_k,Y_{k}= l).\eqno{}$$
On $Y_k=l $,  note that $R_k \geq R_{k-1}\geq  \cdots \geq R_1$, so  $R_k  \geq  l/k$.  
 By the same proof of Lemma 4.1 to discuss the cases $I_k$, $II_k$(a) and (b), if we sum all possible $u_k$ with $R_k\geq l/k$, then
 there exists $C$ in Lemma 4.1 such that
\begin{eqnarray*}
&& \sum_{mL(p)\leq l}\sum_{u_1, \cdots, u_k} P_p({\cal T}_{k}({\bf 0}, u_1, \cdots, u_k),  {\cal G}_k\sim G_k,Y_{k}= l  )\\
&\leq& C  \sum_{u_1, \cdots, u_{k-1}} P_p({\cal T}_{k}({\bf 0}, u_1, \cdots, u_{k-1}),{\cal G}_{k-1}\sim G_{k-1} )\sum_{l\geq mL(p) /k } l P_p({\cal Q}_4(l/4k)).
\hskip 3cm (4.51)
\end{eqnarray*}
By Lemma 2.7, there exist $C_i$ for $i=1,2, 3,4$ independent of  $k$ and $p$ such that
 $$\sum_{l\geq mL(p)/k} lP_p({\cal Q}_4(l /4))\leq C_1(0.5-p)^{-1}\sum_{ l\geq mL(p)/k}e^{-C_2l/k} \leq C_3 e^{C_4 m/l} (0.5-p)^{-1}.\eqno{(4.52)}$$
By using Proposition  2 for the second sum and using (4.52) for the third sum in the right side of (4.51), there exists $C_5$ independent of $p$, $k$ and $n$ such that
$$ \sum_{G_k}\sum_{m\leq l}\sum_{u_1, \cdots, u_k} P_p({\cal T}_{k}({\bf 0}, u_1, \cdots, u_k),  {\cal G}_k\sim G_k,Y_{k}= l  )\leq C^{k}_5 L^{-2} (p)(0.5-p)^{-k-2}e^{-C_4 m/k}.\eqno{(4.53)}$$
Thus, by (4.35), (4.51), and (4.53), if $ m \geq k^3$ for a large $k$, there exists $C$ independent of $p$ and $k$ such that
$$\sum_{u_1, \cdots, u_k} P_p({\cal T}_{k}({\bf 0}, u_1, \cdots, u_k), Y_k \geq mL(p))\leq  \exp(-C m/k) L^{-2} (p)(0.5-p)^{-k-2}.\eqno{(4.54)}$$

By (4.18), 
$$d^kR_{AB}(p,n/2, n)/ dp^{k}=  \sum_{I_1, \cdots, I_k}\sum_{b, f^{I_1}_1, f^{I_1,I_2}_2, \cdots, f^{I_1,I_2\cdots, I_{k}}_{k}} \!\!\!\!\!\!I_1I_2\cdots I_kP_p({\cal R}_{AB}(b,f^{I_1}_1,  f^{I_1,I_2}_2, \cdots, f^{I_1,\cdots, I_{k}}_{k}, n)).\eqno{(4.55)}$$
For fixed $I_1, \cdots, I_k$, let
$$d^kR_{AB}(p,n, I_1, \cdots, I_k)/ dp^{k}=\sum_{b, f^{I_1}_1, f^{I_1,I_2}_2, \cdots, f^{I_1,I_2\cdots, I_{k}}_{k}} P_p({\cal R}_{AB}(b,f^{I_1}_1,  f^{I_1,I_2}_2, \cdots, f^{I_1,\cdots, I_{k}}_{k}, n)).\eqno{(4.56)}$$
If  $\{b, f^{I_1}_1, f^{I_1,I_2}_2, \cdots, f^{I_1,I_2\cdots, I_{k}}_{k}\}\subset [-n/2, n/2]^2$, then  let
$$d^kR_{AB}(p,n/2, n, I_1, \cdots, I_k)/ dp^{k}=\sum_{b, f^{I_1}_1, f^{I_1,I_2}_2, \cdots, f^{I_1,I_2\cdots, I_{k}}_{k}} P_p({\cal R}_{AB}(b,f^{I_1}_1,  f^{I_1,I_2}_2, \cdots, f^{I_1,\cdots, I_{k}}_{k}, n)).\eqno{(4.57)}$$
We can also define in the same ways for $d^kR_{i}(p,n, I_1, \cdots, I_k)/ dp^{k}$ and $d^kR_{i}(p,n/2, n, I_1, \cdots, I_k)/ dp^{k}$ for $i=2,3$, respectively.
By Lemma 3.2,  if $p < 0.5$, note that there is no infinite occupied cluster, so  the following limit exists:
\begin{eqnarray*}
&&\lim_{n\rightarrow \infty} d^kR_{AB}(p,n/2, n, I_1, \cdots, I_k)/ dp^{k}=\lim_{n\rightarrow \infty} d^kR_{AB}(p, n, I_1, \cdots, I_k)/ dp^{k}\\
&=&\lim_{n\rightarrow \infty}\sum_{b, f^{I_1}_1, f^{I_1,I_2}_2, \cdots, f^{I_1,I_2\cdots, I_{k}}_{k}} \!\!\!\!\!\!P_p({\cal R}_{AB}(b,f^{I_1}_1,  f^{I_1,I_2}_2, \cdots, f^{I_1,\cdots, I_{k}}_{k}, n))=d^kR_{AB}(p, \infty, I_1, \cdots, I_k)/ dp^{k}.\hskip .1cm (4.58)
\end{eqnarray*}
If one of pivotals satisfies  $\|f_j^{I_1, \cdots, I_j}\|\geq 2m L(p)$, then $\|u_i\|\geq  2m L(p)$ for some $i$.  Thus, 
\begin{eqnarray*}
&&P_p({\cal R}_{AB}(b,f^{I_1}_1,  f^{I_1,I_2}_2, \cdots, f^{I_1,\cdots, I_{k}}_{k}, n)\mbox{ with }  \|f_j^{I_1, \cdots, I_j}\|\geq 2m L(p)\mbox{ for some j})\\
&\leq& P_p({\cal T}_{k}({\bf 0}, u_1, \cdots, u_k)\mbox{ with } \| u_i\| \geq 2m L(p)\mbox{ for some }i)\hskip 6.5cm (4.59)
\end{eqnarray*}
for  $\{b, f^{I_1}_1, f^{I_1,I_2}_2, \cdots, f^{I_1,I_2\cdots, I_{k}}_{k}\}\subset [-n/2, n/2]^2$.
By (4.58), (4.59), and (4.54), if $n \geq 2m L(p)$ for $m \geq k^3$, then there exists $C$ independent of $k$, $n$, and $p$ such that
\begin{eqnarray*}
&&d^kR_{AB}(p,\infty,  I_1, \cdots, I_k)/ dp^{k}- \exp(-C(m/k) )L^{-2}(p) (0.5-p)^{-k-2}\\
& \leq & d^kR_{AB}(p,n/2, n, I_1, \cdots, I_k)/ dp^{k}\leq d^kR_{AB}(p, \infty, I_1, \cdots, I_k)/ dp^{k}.\hskip 5cm {(4.60)}
\end{eqnarray*}
Similarly, if $n \geq 2m L(p)$ for $m\geq k^3$, then
\begin{eqnarray*}
&&d^kR_{i}(p,\infty,  I_1, \cdots, I_k)/ dp^{k}-\exp(-C (m/k))L^{-2}(p) (0.5-p)^{-k-2}\\
& \leq & d^kR_{i}(p,n/2, n, I_1, \cdots, I_k)/ dp^{k}\leq d^kR_{i}(p, \infty, I_1, \cdots, I_k)/ dp^{k}.\hskip 5.3cm {(4.61)}
\end{eqnarray*}

It follows from (4.18)  and the definition of  $R_{AB}(p, n/2, n)$ that for any fixed $I_1, \cdots, I_k$, 
$$ d^kR_{AB}(p,n/2, n, I_1, \cdots, I_k)dp^k\leq d^kR_{AB}(p,n, I_1, \cdots, I_k)/dp^k\leq d^kR_{AB}(p,\infty, I_1, \cdots, I_k)/dp^k.\eqno{(4.62)}$$
By (4.60)-(4.62), and (4.58), 
\begin{eqnarray*}
&&d^kR_{AB}(p,\infty,  I_1, \cdots, I_k)/ dp^{k}- \exp(-C(m/k) )L^{-2}(p) (0.5-p)^{-k-2}\\
& \leq & d^kR_{AB}(p,n, I_1, \cdots, I_k)/ dp^{k}\leq d^kR_{AB}(p, \infty, I_1, \cdots, I_k)/ dp^{k},\hskip 5.5cm {(4.63)}
\end{eqnarray*}
and
\begin{eqnarray*}
&&d^kR_{i}(p,\infty,  I_1, \cdots, I_k)/ dp^{k}- \exp(-C(m/k))L^{-2}(p) (0.5-p)^{-k-2}\\
& \leq & d^kR_{i}(p, n, I_1, \cdots, I_k)/ dp^{k}\leq d^kR_{i}(p, \infty, I_1, \cdots, I_k)/ dp^{k}.\hskip 6cm {(4.64)}
\end{eqnarray*}
Note that each $I_j$ only has two choices, so by (4.63)-(4.64),  so if $m \geq k^3$ for a large $k$, there exists $C$  independent of $p$ and $k$ and $n$ such that for $A,B,$
$$
|d^kR_{AB}(p,\infty)/ dp^{k}
- d^kR_{AB}(p,n)/ dp^{k}|\leq \exp(-C(m/k) )L^{-2}(p) (0.5-p)^{-k-2}.\eqno{(4.65)}$$
and for $i=2,3$
$$
|d^kR_{i}(p,\infty)/ dp^{k}
- d^kR_{i}(p,n)/ dp^{k}|\leq \exp(-C(m/k) )L^{-2}(p) (0.5-p)^{-k-2}.\eqno{(4.66)}$$

With these observations, we are ready to show Proposition 3. Note that by (4.19),
\begin{eqnarray*}
&&d^k \kappa(p)/dp^k =\lim_{n\rightarrow\infty }d^k G_p(n)/dp^k \\
&=&\sum_{j=1}^{k-1}\sum_{A, B}  f_{AB}(p) d^kR_{AB}(p, \infty)/dp^k+ 2f_{1,2,3}(p)\sum_{i=2,3}d^kR_i(p,\infty)/dp^k+\\
&&+\left[\sum_{j=1}^{k-1} [f_j(p) d^jR_{AB}(p, \infty)/dp^j+\sum_{j=1}^{k-1}g_j(p) \sum_{i=2,3}d^j R_{i} (p, \infty)/dp^j \right]+h_k(p, \infty).\hskip 3cm (4.67)
\end{eqnarray*}
If $m\geq k^3$ for a large $m$, but independent of $p$ and  $n$,  for $n \geq 2mL(p)$, by (4.19), and  (4.65)-(4.67),  there exists $C$ independent of  $k$ and $p$ and $n$ such that
\begin{eqnarray*}
&&|d^k \kappa(p)/dp^k -d^k G_p(n)/dp^k|\leq  4k\exp(-C_4m/k )L^{-2}(p) (0.5-p)^{-k-2}+h_k(p, \infty)\\
&&\leq \exp(-C(m/k) )L^{-2}(p) (0.5-p)^{-k-2}+ h_k(p, \infty).\hskip 7.5cm{(4.68)}
\end{eqnarray*}
As we mentioned in Remark 5, $h_k(p, \infty)$ is uniformly bounded and $\exp(-C(m/k) )L^{-2}(p) (0.5-p)^{-k-2}$ is large for  $p$ near $0.5$, so Proposition 3 follows from (4.68).
 $\Box$.

\section{ An  estimate for   $G_{p_2}(L(p_2))-G_{p_1}(L(p_2))$.}

Recall the definition of $G_p(n)$ in section 3. We know that $G_p(n)$ converges  uniformly to $\kappa''(p)$ on $[0,0.5]$.
We select a sequence $\{p_j\}$  for $j=1,2,\cdots$ such that $p_j\uparrow 0.5$. In addition,  $p_1$ is selected for a large $L(p_1)$ and  $p_j$  is selected such that for some large $M>2$,
$$\iota (0.5-p_j)\leq p_{j+1}-p_{j}  \mbox{ and } {L(p_{j+1})\over L(p_{j})} =  M, \eqno{(5.1)}$$
where $\iota$ is  a positive constant independent of  for all $p_j$ and $M$.
The existence of such a sequence in (5.1) can be proved as follows. 
In fact,  by the reconnection lemma,  Lemma 2.2, (1.14), and (2.9),  we have
$$(0.5-p_j)= (0.5-p_{j+1})M^{3/4+o(1)}.\eqno{(5.2)}$$
So
$$ (0.5-p_j)(1-M^{-3/4+o(1)})= (p_{j+1}-p_j).\eqno{(5.3)}$$
By (5.3),  $\iota $  in (5.1) exists  uniformly for all $p_j$ and $M$. Now we show the following proposition.\\

{\bf Proposition 4.} {\em For   $\delta =0.1$ and a fixed large $M>0$ defined in (5.1), there exists $\alpha=\alpha(M)$   independent of    $\{p_j\}$ defined in (5.1)
  such that  fo any $p_j$ either 
$$(a)\,\,\,\, d^{i+2}\kappa(p_{j-1} )/dp_{j-1}^{i+2}\leq  -\alpha L^{-2}(p_{j-1}) (0.5-p_{j-1})^{-i-2}\mbox{ for some }1\leq i\leq \log M,$$
or 
 $$ (b)\,\,\,G_{p_{{j-1}}}(  L(p_{{j}}))- G_{p_{{j}}}(L(p_{{j}}))\leq  -M^{-\delta/2+o(1)} L(p_{j-1}) (0.5-p_{j-1})^{-2}. \eqno{}$$}




{\bf Proof.} We take $\epsilon=M^{-1+\delta}$.
We decompose
\begin{eqnarray*}
&&G_{p_{t-1}}( L(p_{t})) -G_{p_{t}}( L(p_{t}))=[G_{p_{t-1}}( L(p_{t}))- G_{p_{t-1} }(  \epsilon  L(p_{t}))]
+ [G_{p_{t-1} }( \epsilon  L(p_{t}))-G_{p_t}( \epsilon  L(p_{t}))]\\
&&[-G_{p_{t}}(L(p_{t}))+G_{p_t}(\epsilon L(p_{t}))]=I+II+III
,\hskip 8.5 cm (5.4)
\end{eqnarray*}
where  $I, II, III$  are the differences in the three square brackets  in (5.4).
By  (3.42), if $\epsilon=M^{-1+\delta}$, then there exists $C_1$ independent of $M$ and $p_t$ such that
\begin{eqnarray*}
&&I= [G_{p_{t-1}}(  L(p_{t}))- G_{p_{t-1} }(  \epsilon  L(p_{t}))]
= [G_{p_{t-1}}(  L(p_{t}))- G_{p_{t-1} }(  M^{\delta}  L(p_{t-1}))]\\
&\leq & \exp(-C_1M^\delta) L^{-2}( p_{t-1})(0.5-p_{t-1})^{-2}.\hskip 9cm (5.5)
\end{eqnarray*}
By (3.41),  (5.1), and (5.2), if $\epsilon=M^{-1+\delta}$, then
 \begin{eqnarray*}
 &&III=[-G_{p_t}(L(p_{t}))+G_{p_t}( \epsilon  L(p_{t}))] \\
 &\leq & -\epsilon^{-1/2} M^{-1/2+o(1)} L^{-2}(p_{t-1})(0.5-p_{t-1})^{-2} \leq - M^{-\delta/2+o(1)} L^{-2}(p_{t-1})(0.5-p_{t-1})^{-2}.\hskip 1.5cm (5.6)
 \end{eqnarray*}

We now work on II.    For each large $M$,  let $(i_0!)^{1/4}=M$.
Thus, for a large $M$,  by using  Stirling's formula,
$$i_0\leq \log M.\eqno{(5.7)}$$
By Proposition 2, if we choose   a large $i_0$  but fixed  and for $i\geq i_0$ such that
 $$(i! )^{1/2} M^{\delta} \geq (i!)^{1/4} \geq C^i\mbox{ for the $C$ in Proposition 2},$$
 then
$$d^iG_{p_{t-1} }( \epsilon L(p_t)) /d{p_{t-1}^i}\geq- (i!)^{1/2}M^{-\delta}L^{-2}(p_{t-1})(p_t-p_{t-1})^{-i-2}.\eqno{(5.8)}$$
For $1\leq i\leq i_0\leq (\log M)$,  
 we first assume that   
 $$d^iG_{p_{t-1} }( \epsilon L(p_t)) /d{p_{t-1}^i}\geq- (i!)^{1/4}M^{-\delta+o(1)}L^{-2}(p_{t-1})(p_t-p_{t-1})^{-i-2},\eqno{(5.9)}$$
 for each  $1\leq i\leq i_0-1$ with a large $M$.
 By Taylor's theorem, there exists $\xi\in (p_{t-1}, p_t)$ such that
$$G_{p_t}( \epsilon L(p_t))- G_{p_{t-1}}(\epsilon L(p_t))=\sum_{k=1}^{K-1}(p_t-p_{t-1})^k{d^k(G_{p_{t-1} }(\epsilon L(p_t)))\over k!d p^k_{t-1}}+ (p_t-p_{t-1})^K{d^K(G_{\xi }(\epsilon L(p_t)))\over K!d \xi^K}.\eqno{(5.10)}$$
By (5.8) and the assumption  (5.9), for a large $K > M$, 
$$\sum_{k=1}^{K-1}(p_t-p_{t-1})^k{d^k(G_{p_{t-1} }(\epsilon L(p_t)))\over k!d p^k_{t-1}}\geq -M^{-\delta+o(1)} L^{-2}(p_{t-1})(0.5-p_{t-1})^{-2}\sum_{k=1}^{K-1} (k!)^{-1/2+o(1)}.$$
Thus,  
$$\sum_{k=1}^{K-1}(p_t-p_{t-1})^k{d^k(G_{p_{t-1} }(\epsilon L(p_t)))\over k!d p^k_{t-1}}\geq -M^{-\delta+o(1)} L^{-2}(p_{t-1})(0.5-p_{t-1})^{-2}.\eqno{(5.11)}$$

Let us work on the remainder in (5.10).  By Proposition 2, (5.1)-(5.3), and  Stirling's formula,  note that $\delta=0.1 \leq 1/2$,  so for large $K$ and for the $C$ in Proposition 2,
\begin{eqnarray*}
&&(p_t-p_{t-1})^K\left|{d^K(G_{\xi }(\epsilon L(p_t)))\over K!d \xi^K}\right|
\leq  C^K(p_t-p_{t-1})^K L^{-2}(\xi) (0.5-\xi)^{-K-2}/K!\\
&\leq &  C^Ke^KM^{3K/4} L^{-2}(p_{t-1}) (0.5-p_{t-1})^{-2}/K!
\leq  (e CM /K)^K L^{-2}(p_{t-1}) (0.5-p_{t-1})^{-2}.\hskip 2cm (5.12)
\end{eqnarray*}
By taking $K$ large  in (5.12), 
$$(p_t-p_{t-1})^K\left|{d^K(G_{\xi }(\epsilon L(p_t)))\over K!d \xi^K}\right|\leq  M^{-2\delta+o(1)} L(p_{t-1}) (0.5-p_{t-1})^{-2}.\eqno{(5.13)}$$
By substituting (5.12) and (5.13) into (5.10),  for a large $M$
$$G_{p_t}( \epsilon L(p_t))- G_{p_{t-1}}(\epsilon L(p_t))\geq -M^{-\delta+o(1)} L^{-2}(p_{t-1}) (0.5-p_{t-1})^{-2}. \eqno{(5.14)}$$
Thus,  by (5.14) for a large $M$, 
$$II=G_{p_{t-1}}( \epsilon L(p_{t}))- G_{p_{t}}(\epsilon L(p_t))\leq M^{-\delta+o(1)} L^{-2}(p_{t-1}) (0.5-p_{t-1})^{-2}. \eqno{(5.15)}$$
If we use the estimates in (5.5), (5.6) and (5.15) together for a large $M$ in (5.4), $III$ dominates in (5.4). So  Proposition 4 (b)  follows.

We now assume that (5.9) does not hold. In other words,  for  some $1\leq i\leq i_0\leq \log M$,
$$d^iG_{p_{t-1} }( \epsilon L(p_t)) /d{p_{t-1}^i}<- (i!)^{1/4}M^{-\delta+o(1)}L^{-2}(p_{t-1})(p_t-p_{t-1})^{-i-2}.\eqno{(5.16)}$$
It remains  to show that (5.16) implies Proposition 4 (a).  Note that
$\epsilon L(p_t)= M^{\delta+o(1)} L(p_{t-1})$   and $i\leq i_0 \leq  \log M$, so 
by Proposition  3 and (5.16) for a large $M$ with $M^\delta > i^3$,
\begin{eqnarray*}
&&d^{i+2}\kappa(p_{t-1})/dp^{i+2}_{t-1} 
\leq  d^iG_{p_{t-1} }( \epsilon L(p_t)) /d{p_{t-1}^i}+  \exp(-M^{\delta+o(1)}/i)L^{-2}(p_{t-1})(0.5-p_{t-1})^{-i-2}\\
&\leq &\left(-(i!)^{1/4+o(1)} M^{-\delta+o(1)}+ \exp(-M^{\delta/2+o(1)})\right)L^{-2}(p_{t-1})(0.5-p_{t-1})^{-i-2}\\
&\leq &-(i!)^{1/4+o(1)} M^{-\delta+o(1)} L^{-2}(p_{t-1})(0.5-p_{t-1})^{-i-2} .\hskip 7.5cm (5.17)
\end{eqnarray*}
Thus, there exist $1\leq i\leq i_0$ and  $\alpha=\alpha(M)$ such that
$$d^{i+2}\kappa(p_{t-1})/dp^{i+2}_{t-1}\leq  -\alpha L^{-2}(p_{t-1}) (0.5-p_{t-1})^{-i-2}.\eqno{(5.18)}$$
(5.18) indeed implies Proposition 4 (a). 
Therefore, Proposition 4 follows. $\Box$.

\section{ Proof of theorem.}
\subsection{ Proof of  the upper bound of the theorem.}
If $ p< 0.5$,  then  by Proposition 2 and (2.3), there exists $C$ independent of $p$ such that
$$|\kappa'''(p)|  \leq C L^{-2}(p) (0.5-p)^{-3}\leq (0.5-p)^{-1/3+o(1)}.\eqno{(6.1)}$$
If $p >0.5$, by (1.5), we still have (6.1).

\subsection{ Proof of the lower bound of the theorem.}
We use the sequence $\{p_j\}$ in (5.1) for large $M$ such that
$$\iota (0.5-p_j)\leq p_{j+1}-p_{j}  \mbox{ and } {L(p_{j+1})\over L(p_{j})} =  M. \eqno{(6.2)}$$
By Proposition 4,  for given $M$ and $K=\log M$, there exist $C=C(M)$ independent of $p_j$ such that either
$$d^{i+2}\kappa(p_{j-1})/dp_{j-1}^{i+2} \leq -C L^{-2}(p_{j-1}) (p_j-p_{j-1})^{-i-2} =-(p_j-p_{j-1})^{-i+2/3+o(1)}\mbox{ for some }1\leq i\leq K, \eqno{(6.3)}$$
or 
$$G_{p_j}(L(p_j))- G_{p_{j-1}}(L(p_j)) \geq  M^{-\delta/2+o(1)} L(p_{j-1}) (0.5-p_{j-1})^{-2}.\eqno{(6.4)}$$
We assume that (6.4) holds.
By (3.33), we take $n\geq L^2(p_j)$ large such that
$$\left|\kappa''(p_{{j}}) - G_{p_{{j}}}(n) \right|\leq \exp(-L(p_{{j}}))\mbox{ and } \left|\kappa''(p_{{j-1}}) - G_{p_{{j-1}}}(n) \right|\leq \exp(-L(p_{{j}})).\eqno{(6.5)}$$
Thus,
\begin{eqnarray*}
&&\kappa''(p_{{j}})-\kappa''(p_{{j-1}})=\\
&&\left[\kappa''(p_{{j}})-G_{p_{{j}}}(n)\right]+\left[
G_{p_{{j-1}}}(n)-\kappa''(p_{{j-1}})\right]
+\left[  G_{p_{{j}}}(n)-  G_{p_{{j}}}(L(p_{{j}}))\right]\\
&&+\left[G_{p_{{j}}}(L(p_{{j}}))-G_{p_{{j-1}}}(L(p_{{j}}))\right]
+\left[-G_{p_{{j-1}}}(n)+G_{p_{{j-1}}}(L(p_{{j}}))\right]\\
&=& I+II+III+IV+V,\hskip 11.7cm (6.6)
\end{eqnarray*}
where $I, II, III, IV, V$ are the differences in the five square brackets  in (6.6).
By (6.5) and (2.3),
$$|I|\leq (0.5-p_{{j}})^{2}\mbox{ and } |II|\leq (0.5-p_{{j}})^{2}.\eqno{(6.7)}$$
By  (3.40), (6.2) and (5.2), 
there exists $C_1$  independent of  $M$  such that
$$III= G_{p_{{j}}}(n)-G_{p_{{j}}}(L(p_{{j}}))\geq  C_1 L^{-2}(p_{{j}})(0.5-p_{{j}})^{-2}\geq  M^{-1/2+o(1)}  L^{-2}(p_{{j-1}})(0.5-p_{{j-1}})^{-2}.\eqno{(6.8)}$$
By  (6.4),  
$$IV =G_{p_{{j}}}(L(p_{{j}}))-G_{p_{{j-1}}}(L(p_{{j}}))\geq M^{-\delta/2+o(1)}  L^{-2}(p_{{j-1}})(0.5-p_{{j-1}})^{-2}.\eqno{(6.9)}$$
By (3.42) and (5.1), there exist $C_{2}$ and $C_3$ independent of  $M$  such that
$$|V|\leq C_{2}\exp(-C_3 M)L^{-2}(p_{{j-1}})(0.5-p_{{j-1}})^{-2}.\eqno{(6.10)}$$
 Thus, $III$  and $IV$ will dominate the others in (6.6) for a large $M$.
Together with (6.7)--(6.10),  if  $M$ is large, but a fixed  number, then
$$\kappa''(p_{{j}})-\kappa''(p_{{j-1}})\geq M^{-\delta/2+o(1)} L^{-2}(p_{{j-1}}) (0.5-p_{{j-1}})^{-2}.\eqno{(6.11)}$$
We want to remark that $M$ is selected large depending only on $C_1$ $C_2$,  and $C_3$ in (6.9) and (6.10) such that (6.11) holds. This selection of $M$ is uniformly for all $\{p_j\}$.
Thus, $K=\log M$ is also a fixed number uniformly for all $\{p_j\}$ in (6.3).
By the mean value theorem,  (6.11), (5.3),  and (2.3), there exist $\xi_j\in (p_{{j-1}}, p_{{j}})$ and  $C_4=C_4(M)$ independent of $\xi$ such that
such that
\begin{eqnarray*}
&& \kappa'''(\xi_j)\geq M^{-\delta/2+o(1)}L^{-2}(p_{{j-1}}) (0.5-p_{{j-1}})^{-2} (p_{{j}}-p_{{j-1}})^{-1}\\
&\geq &M^{-\delta/2+o(1)} L^{-2}(p_{{j-1}}) (0.5-p_{{j-1}})^{-3}\geq M^{-3/4-\delta/2+o(1)} L^{-2}(p_{{j}}) (0.5-p_{{j}})^{-3} \\
&\geq &  M^{-3/4-\delta/2+o(1)}(0.5-\xi_j)^{-1/3+o(1)}\geq C_4(0.5-\xi_j)^{-1/3+o(1)}.\hskip 5.5cm (6.12)
\end{eqnarray*}

Now we assume that (6.3) holds. 
For the sequence $\{p_j\}$ defined in (5.1), there exists $C_5=C_5( M)$ for the $M=\log K$ defined in Proposition 4 such that
$$|d^{i+2}\kappa(p_j) /dp_j^{i+2}|\geq  C_{5} L^{-2} (p_j)(0.5-p_j)^{-i-2} \mbox{ for some }1\leq i\leq K.\eqno{(6.13)}$$
(6.12) and (6.13) show that $\kappa(p)$ has a singularity at $0.5$. Moreover, together with (6.12), we want to show that $\kappa(p)$ is not 
third differentiable.
We only need to discuss that  $i\geq 2$ and show (6.13) to hold for $i+1$.   Let $f_i(p) =d^i \kappa(p) /d^i p$.
We show that if  (6.13) does not  hold for $i+1$,  then there would be a contradiction.  We assume that there is $\epsilon >0$ such that 
for $p\in [ p_j-\epsilon(p_{j}-p_{j-1}), p_j]$,
$$| f_{i+1}(p)|\leq \epsilon  (C_5/8) L^{-2} (p_j)(0.5-p_j)^{-i-1}. \eqno{(6.14)}$$
By (6.14), we take  $x= p_j-\epsilon (p_{j}-p_{j-1})/2$
and $y= p_j-\epsilon (p_j-p_{j-1})$ such that
$$| f_{i+1}(y)|\leq  \epsilon (C_5/8) L^{-2} (p_j)(0.5-p_j)^{-i-1}\mbox{ and } | f_{i+1}(x)|\leq \epsilon  (C_5/8)  L^{-2} (p_j)(0.5-p_j)^{-i-1}.\eqno{(6.15)}$$
If we use  the mean value theorem for $f_{i+1}(p)$ at $y$ and $x$, then by (5.3) 
there exists $z\in  (p_j-\epsilon(p_{j}-p_{j-1}),  p_j-\epsilon(p_{j}-p_{j-1})/2)$ such that
$$|d^{i+2}\kappa(z) /dz^{i+2}|\leq (C_5/2)L^{-2} (p_j)(0.5-p_j)^{-i-2}.\eqno{(6.16)}$$
If  we use the mean value theorem again for $f_{i+2}(p)$ at $z$ and $p_j$, then by (5.1), (6.13), and (6.16),
there exists $\beta\in ( p_j-\epsilon (p_{j}-p_{j-1}), p_j)$ such that
$$|{d^{i+3}\kappa(\beta) \over d\beta^{i+3}}|\geq  ( C_5/2) L^{-2} (p_j)(0.5-p_j)^{-i-2}\epsilon^{-1} (p_{j}-p_{j-1})^{-1}\geq 
(C_5/2) M^{-1}L^{-2} (\beta)(0.5-\beta)^{-i-3}\epsilon^{-1}.\eqno{(6.17)}$$
By taking $\epsilon=(C_5/4)^{-1} M^{-1} C^i)^{-1}$ in (6.17) for the $C$ defined in Proposition 2, 
$$|d^{i+3}\kappa(\beta) /d\beta^{i+3}|\geq  2C^i L^{-2} (\beta)(0.5-\beta)^{-i-3} .\eqno{(6.18)}$$
On the other hand,  by Proposition 2,  
$$|d^{i+3}\kappa(\beta) /d\beta^{i+3}|\leq C^i L^{-2} (\beta)(0.5-\beta)^{-i-3}.\eqno{(6.19)}$$
(6.18) and (6.19) cannot hold together, so (6.14) cannot hold. Thus,  for $i\geq 2$, there exists $\alpha\in (p_{j}-\epsilon(p_{j}-p_{j-1}), p_j)$ such that   
$$|d^{i+1}\kappa(\alpha) /d\alpha ^{i+1}|\geq  \epsilon (C_5/8) L^{-2} (p_j)(0.5-p_j)^{-i-1}\geq \epsilon (C_5/8) L^{-2} (\alpha)(0.5-\alpha)^{-i-1}. \eqno{(6.20)}$$
Thus,  there exist $\alpha \in (p_{j-1}, p_j)$ and  $C_6=C_6(M, K)$ such that
$$|d^{i+1}\kappa(\alpha) /d\alpha ^{i+1}|\geq  C_6 L^{-2} (\alpha)(0.5-\alpha)^{-i-1}. \eqno{(6.21)}$$
We continue this way for $i-1$ times to show that    there is  $\alpha\in (p_{j-1}, p_j)$ such that
$$|d^{3}\kappa(\alpha) /d\alpha^{3}|\geq   C_7 L^{-2} (\alpha) (0.5-\alpha)^{-3} \mbox{ for some }C_7=C_7(M).\eqno{(6.22)}$$
Therefore, the lower bound of the Theorem follows from (6.12) and (6.22) for $p < 0.5$.
If $p > 0.5$, the lower bound of the Theorem follows from  (1.5), (6.12) and (6.22). $\Box$\\

\begin{center}{\large \bf References} \end{center}
Aizenman, M., Dulpantier, B. and Alharony, A. (1999). Path-crossing exponents and the external perimeter in 2D percolation. {\em Phy. Rev. Lett.}
{\bf 63}, 817--835.\\
 Dunford, N. and Schwartz, T. (1958). {\em Linear operators.}  {\bf 1},
{ Wiley-Intrscience}, New York.\\ 
Grimmett, G. (1981). On the differentiability of the number of clusters per vertex in percolation model. {\em J. Lond. Math. Soc.} (2) {\bf 23}, 372--384.\\
Grimmett, G. (1999). {\em Percolation}. { Springer-Verlag}, New York.\\
Higuchi, Y., Takei, M. and  Zhang, Y. (2012). Scaling relations for two-dimensional Ising percolation. {\em J. Stat. Phys.} {\bf 148}, 777--799.\\
Kesten, H. (1982). {\em Percolation Theory for Mathematicians.}  { Birkhauser}, Boston.\\
Kesten, H. (1987). Scaling relations for 2D-percolation. {\em Comm.
Math. Phys.} {\bf 109}, 109--156.\\
Kesten, H., Sidoravicius,  V. and Zhang, Y. (1998). Almost all words are seen in critical site percolation on the triangular lattice. {\em Electron. J. Probab.} {\bf 3}  no. 10, 75pp.\\
Lawler, G. F., Schramm, O. and  Werner, W. (2002) One-arm exponent for critical 2D percolation. {\em Electron. J. Probab}. {\bf 7}, 13pp.\\
 Morrow, G. J. and Zhang, Y. (2005). The sizes of the pioneering, lowest crossing and pivotal sites in critical percolation on the triangular lattice. {\em Ann.  Appl. Probab.} {\bf 5}, 1832--1886.\\
Smirnov, S. (2001). Critical percolation in the plane: conformal invariance, Cardy's formula, scaling limits. {\em  CR. Acad.  Sci. I-Math.} {\bf 333},  239--244.\\
Smirnov, S. and Werner, W. (2001). Critical exponent for two dimensional percolation. {\em Math.  Res. Lett. } {\bf 8}, 729--744.\\
Sykes, M. F., and Essam, J. W. (1964). Exact critical percolation probabilities for site and bond problems in two dimensions. {\em J. Math. Phys.}  {\bf 5}, 1117--1127.\\
Zhang, Y. (2011). A derivative formula for the free energy function {\em J. Stat. Phys.} {\bf 146}, 466--473.\\

Yu Zhang\\
Department of Mathematics\\
University of Colorado\\
Colorado Springs, CO 80933\\
yzhang3@uccs.edu

\end{document}